\documentclass{amsart}
\usepackage{amssymb, amsbsy, amsthm, amsmath, amstext, amsopn, verbatim,
multicol}
\usepackage[all]{xy}
\usepackage{amsfonts}
\usepackage{amscd}
\newtheorem{thm}{Theorem}[section]
\newtheorem{lemma}[thm]{Lemma}
\newtheorem{lemma-defn}[thm]{Lemma-Definition}
\newtheorem{sublemma}[thm]{Sublemma}
\newtheorem{corollary}[thm]{Corollary}
\newtheorem{prop}[thm]{Proposition}
\newtheorem{notation}[thm]{Notation}

\newtheorem{hypotheses}[thm]{Hypotheses}

\theoremstyle{definition}

\newtheorem{VC}[thm]{Vanishing Condition}

\newtheorem{defn}[thm]{Definition}

\newtheorem{example}[thm]{Example}

\newtheorem{construction}[thm]{Construction}

\theoremstyle{remark}

\newtheorem{remark}[thm]{Remark}

\newtheorem{fact}[thm]{Fact}

\begin{document}

\numberwithin{equation}{section}
\newcommand{\hs}{\mbox{\hspace{.4em}}}
\newcommand{\ds}{\displaystyle}
\newcommand{\bd}{\begin{displaymath}}
\newcommand{\ed}{\end{displaymath}}
\newcommand{\bcd}{\begin{CD}}
\newcommand{\ecd}{\end{CD}}

\newcommand{\proj}{\operatorname{Proj}}
\newcommand{\bproj}{\underline{\operatorname{Proj}}}
\newcommand{\spec}{\operatorname{Spec}}
\newcommand{\bspec}{\underline{\operatorname{Spec}}}
\newcommand{\pline}{{\mathbf P} ^1}
\newcommand{\pplane}{{\mathbf P}^2}
\newcommand{\coker}{{\operatorname{coker}}}
\newcommand{\ldb}{[[}
\newcommand{\rdb}{]]}

\newcommand{\Sym}{\operatorname{Sym}^{\bullet}}
\newcommand{\Symp}{\operatorname{Sym}}
\newcommand{\Pic}{\operatorname{Pic}}
\newcommand{\AAut}{\operatorname{Aut}}
\newcommand{\PAut}{\operatorname{PAut}}
\newcommand{\GKdim}{\operatorname{GKdim}}
\newcommand{\Rees}{{\cR}}
\newcommand{\KLpair}{\bK\hskip -3pt \bL}
\newcommand{\too}{\twoheadrightarrow}
\newcommand{\cA}{{\mathcal A}}
\newcommand{\cV}{{\mathcal V}}
\newcommand{\cM}{{\mathcal M}}
\newcommand{\bA}{{\mathbf A}}
\newcommand{\cB}{{\mathcal B}}
\newcommand{\cC}{{\mathcal C}}
\newcommand{\cD}{{\mathcal D}}
\newcommand{\D}{{\mathcal D}}
\newcommand{\cs}{{\mathbb C} ^*}
\newcommand{\bbC}{{\mathbb C}}
\newcommand{\boldc}{{\mathbf C}}
\newcommand{\cE}{{\mathcal E}}
\newcommand{\cF}{{\mathcal F}}
\newcommand{\cG}{{\mathcal G}}
\newcommand{\cH}{{\mathcal H}}
\newcommand{\cI}{{\mathcal I}}
\newcommand{\cK}{{\mathcal K}}
\newcommand{\bK}{{\mathbf K}}
\newcommand{\bL}{{\mathbf L}}
\newcommand{\cL}{{\mathcal L}}
\newcommand{\baL}{{\overline{\mathcal L}}}
\newcommand{\M}{{\mathcal M}}
\newcommand{\cJ}{{\mathcal J}}
\newcommand{\bM}{{\mathbf M}}
\newcommand{\cN}{{\mathcal N}}
\newcommand{\N}{{\mathcal N}}
\newcommand{\theo}{{\mathcal O}}
\newcommand{\cP}{{\mathcal P}}
\newcommand{\cR}{{\mathcal R}}
\newcommand{\cS}{{\mathcal S}}
\newcommand{\cT}{{\mathcal T}}
\newcommand{\boldp}{{\mathbf P}}
\newcommand{\boldq}{{\mathbf Q}}
\newcommand{\cQ}{{\mathcal Q}}
\newcommand{\cO}{{\mathcal O}}
\newcommand{\cU}{{\mathcal U}}
\newcommand{\cX}{{\mathcal X}}
\newcommand{\cW}{{\mathcal W}}
\newcommand{\boldz}{{\mathbf Z}}
\newcommand{\Qgr}{\operatorname{Qgr-}\hskip -2pt}
\newcommand{\qgr}{\operatorname{qgr-}\hskip -2pt}
\newcommand{\gr}{\operatorname{gr-}\hskip -2pt}
 \newcommand{\Ggr}{\operatorname{Gr-}\hskip -2pt}
\newcommand{\module}{\operatorname{mod-}\hskip -2pt}
\newcommand{\Module}{\operatorname{Mod-}\hskip -2pt}
\newcommand{\Tors}{\operatorname{Tors-}\hskip -2pt}
\newcommand{\tors}{\operatorname{tors-}\hskip -2pt}
\newcommand{\skl}{\operatorname{Skl}}
\newcommand{\coh}{\operatorname{coh}}
\newcommand{\QCh}{\operatorname{QCh}}
\newcommand{\End}{\operatorname{End}}
\newcommand{\uEnd}{\underline{\operatorname{End}}}
\newcommand{\Hom}{\operatorname{Hom}}
\newcommand{\uHom}{\underline{\operatorname{Hom}}}
\newcommand{\Ext}{\operatorname{Ext}}
\newcommand{\bExt}{\operatorname{\bf{Ext}}}
\newcommand{\Tor}{\operatorname{Tor}}
\newcommand{\uExt}{\underline{\operatorname{Ext}}}
\newcommand{\ASthree}{\underline{\operatorname{AS}}_3}
\newcommand{\inv}{^{-1}}
\newcommand{\airtilde}{\widetilde{\hspace{.5em}}}
\newcommand{\airhat}{\widehat{\hspace{.5em}}}
\newcommand{\nt}{^{\circ}}
\newcommand{\rank}{\operatorname{rank}}
\newcommand{\supp}{\operatorname{supp}}
\newcommand{\hd}{\operatorname{hd}}
\newcommand{\id}{\operatorname{id}}
\newcommand{\res}{\operatorname{res}}
\newcommand{\lrar}{\leadsto}
\newcommand{\im}{\operatorname{Im}}
\newcommand{\HH}{\operatorname{H}}
\newcommand{\TF}{\operatorname{TF}}
\newcommand{\Bun}{\operatorname{Bun}}
\newcommand{\Hilb}{\operatorname{Hilb}}
\newcommand{\nthord}{^{(n)}}
\newcommand{\Aut}{\underline{\operatorname{Aut}}}
\newcommand{\Gr}{\operatorname{Gr}}
\newcommand{\GR}{\operatorname{GR}}
\newcommand{\Fr}{\operatorname{Fr}}
\newcommand{\GL}{\operatorname{GL}}
\newcommand{\SL}{\operatorname{SL}}
\newcommand{\on}{\operatorname}
\def\Ext{\operatorname {Ext}}
\def\Hom{\operatorname {Hom}}
\def\bbZ{{\mathbb Z}}
\newcommand{\rref}{ref\ }

\title{Sklyanin Algebras and Hilbert Schemes of Points}
\author{T. A. Nevins and J. T. Stafford}
\address{Department of Mathematics, University of Michigan, Ann Arbor,
MI 48109, USA.} 
   \thanks{The 
first author was supported by an NSF Postdoctoral Research Fellowship and an
MSRI postdoctoral fellowship. The second author was supported in part by the
NSF through the  grants DMS-9801148 and DMS-0245320.}

\begin{abstract}
We construct projective moduli spaces for 
torsion-free sheaves on noncommutative projective planes.  These 
moduli spaces vary smoothly in the parameters describing 
the noncommutative  plane and have good properties
analogous to those of moduli spaces of sheaves over
the usual (commutative) projective plane $\mathbf P^2$.  

The generic noncommutative plane corresponds to the Sklyanin algebra
$S=\skl(E,\sigma)$ constructed from an   automorphism $\sigma$
of infinite order on an elliptic curve $E\subset \mathbf P^2$. In this case,
 the fine moduli space of line bundles over $S$ 
  with first Chern class zero and Euler characteristic $1-n$ provides a
 symplectic variety that is a deformation
  of the Hilbert scheme of $n$ points on $\mathbf P^2\smallsetminus E$.
  \end{abstract}
 \maketitle 
 \tableofcontents 
 
\clearpage
 \section{Introduction}\label{introduction} 
 
It has become apparent that vector bundles and
instantons on noncommutative spaces
play an important role in representation theory (for example, via
quiver varieties \cite{Nakajima, BGK, BGK2, Ginzburg, LeBruynbook}),
string theory \cite{CDS, Nekrasov-Schwarz, Seiberg-Witten, KKO},
and
integrable systems \cite{BW1, BW2}.  In this paper we make the first
systematic application of moduli-theoretic techniques from algebraic
geometry to a problem in this area, the study
of coherent sheaves on noncommutative   projective planes. 

The circle of ideas leading to this paper begins with the work of Cannings and
Holland on 
\emph{the first Weyl algebra} (or ring of differential operators  on the affine
line) $A_1 = \mathcal{D}(\mathbf A^1) = \mathbb{C}\{x,\partial\}/(\partial
x-x\partial-1)$.  \cite{CH} gives a (1-1) correspondence between
isomorphism classes of right ideals of $A_1$ and points of a certain infinite
dimensional Grassmannian.   Curiously, this space appears in a completely
different context: it is precisely the 
 {\em adelic Grassmannian} $\on{Gr}^{\on{ad}}$ that plays
an essential r\^ole in Wilson's study of rational solutions of the KP hierarchy
of soliton theory \cite{Wilson}. The space 
 $\on{Gr}^{\on{ad}}$ is related to a number of other interesting spaces.
  In particular it decomposes into a
disjoint union  $\on{Gr}^{\on{ad}} = \coprod_{n\geq 0}{\mathcal C}_n$ of 
certain quiver varieties ${\mathcal C}_n$, the
{\em completed Calogero-Moser phase spaces}. The varieties $\mathcal
C_n$ appear naturally in yet another way, as  deformations  of the Hilbert
scheme $(\mathbf A^2)^{[n]}$ of $n$ points in the affine plane (this deformation
is the twistor family of the hyperk\"ahler metric on 
$(\mathbf A^2)^{[n]}$).  

This suggests that there might be a natural way to  relate  right ideals
of $A_1$ more directly to both $\mathcal C_n$ and $(\mathbf A^2)^{[n]}$. This 
is done in \cite{Le,BW1,BGK,BGK2}, where  noncommutative projective geometry
is used to refine the work of \cite{CH}.  As will be
explained shortly, there is a noncommutative projective plane $\mathbf
P^2_{\hbar}$ and a ``restriction'' functor from the category of modules on
$\boldp^2_{\hbar}$ to the category of $A_1$-modules.  Under this map,
 isomorphism classes of  right ideals of $A_1$  are in  (1-1)
correspondence with isomorphism classes of line bundles
 on $\mathbf
P^2_{\hbar}$ that are ``trivial at infinity.'' These line bundles have  a
natural integer invariant  $c_2\geq 0$ analogous to the second Chern
 class and the line bundles with $c_2=n$ are
then in bijection with points of ${\mathcal C}_n$. This corresponds nicely with the 
results of \cite{CH} and \cite{Wilson}. Indeed, if one regards $\mathbf
P^2_{\hbar}$ as a deformation of $\boldp^2$, then 
$A_1$ provides a deformation of the ring of functions on
 $\mathbf A^2\subset  \mathbf P^2$ and  the line
bundles on $\mathbf P^2_{\hbar}$  with $c_2=n$ correspond naturally to 
points of a deformation of $(\mathbf A^2)^{[n]}$. 

The classification of line bundles on $\mathbf P^2_\hbar$
also has an analogue for  vector bundles of rank
greater than one  \cite{KKO}. In this case, the classification is related to
 influential work of Nekrasov and Schwarz \cite{Nekrasov-Schwarz}
 in string theory concerning instantons on a noncommutative $\mathbb R^4$. 

This story raises a number of problems:
\begin{enumerate}  
 \item The plane $\boldp^2_\hbar$ is  one of several
families of noncommutative planes. Construct moduli spaces that classify
vector bundles or even   torsion-free  coherent sheaves on all such
planes. Show, in particular,  that the classification results of
\cite{BW2,BGK,KKO} come from moduli space structures.

 \item Prove  that these  moduli spaces behave well in families. For example,
 when the noncommutative plane is a deformation of (the category of coherent
 sheaves on) $\mathbf P^2$, this should provide deformations of the
 Hilbert schemes of points $(\mathbf P^2)^{[n]}$.
\item The Hilbert scheme $(\mathbf A^2)^{[n]}$
 has an algebraic symplectic structure induced by  the hyperk\"ahler
metric. Construct Poisson or symplectic structures on the  analogous moduli spaces
 defined in (2) and  study the resulting Poisson and symplectic
geometry.
 \end{enumerate}
We solve all these problems in the present paper.
 For example, we obtain
 ``elliptic'' deformations of the Hilbert schemes of points $(\mathbf
 P^2)^{[n]}$ and $(\mathbf
 P^2\smallsetminus E)^{[n]}$ for any plane
  elliptic curve $E$. The deformations $(\boldp^2)^{[n]}$ and
 $(\boldp^2\smallsetminus E)^{[n]}$ do indeed
 carry Poisson and  holomorphic symplectic structures respectively.

\subsection{Line bundles on Quantum Planes} An introduction to the general
theory of noncommutative projective geometry can be found in  \cite{Stafford,
SVdB} and formal definitions are given in Section~\ref{section2}. Roughly
speaking, and in accordance with Grothendieck's  philosophy,\footnote{``In
short, as A. Grothendieck taught us,  to do geometry you really don't  need a
space, all you need is a category of sheaves on this would-be space''
\cite[p.83]{Ma}.}
 a noncommutative projective plane is an abelian category with
the fundamental  properties of the category $\coh(\mathbf P^2)$ of coherent
sheaves on the projective plane.  These
categories have been classified and they have a  surprisingly  rich geometry. 

 To be more precise, let $S=k\oplus \bigoplus_{i\geq 1}S_i$ be a connected
graded noetherian $k$-algebra over a field $k$. For simplicity we assume that
$\on{char} k=0$ in the introduction, although many of our results are proved 
for general fields. By analogy with the commutative
case one  regards $\qgr S$, the category of finitely generated  graded
$S$-modules modulo those of finite length, as the category of coherent sheaves
on the noncommutative (and nonexistent) 
projective variety  $\on{Proj}(S)$. As is
explained in \cite[Section~11]{SVdB}, the noncommutative analogues of
$\coh(\boldp^2)$ are exactly the categories  of the form $\qgr S$ where
 $S$ is an AS (Artin-Schelter) regular algebra with
the Hilbert  series $(1-t)^{-3}$ of a polynomial ring in three variables (to denote
which we
write $S\in \ASthree$; see
Definition~\ref{as-defn}). Accordingly, we will refer to $\qgr S$ for such an
algebra  $S$ as a {\em noncommutative} or {\em quantum projective plane}.

The algebras $S\in \ASthree$ have in turn been classified in terms of geometric
data \cite{ATV1}. When $\qgr S\not\simeq\coh(\mathbf P^2)$, the
classification is in terms of  commutative data $(E,\sigma)$, where
$E\hookrightarrow \mathbf P^2$ is a (possibly singular)  plane cubic curve and
$\sigma\in \on{Aut}(E)$ is a non-trivial automorphism.  Moreover, $S$ is then
determined by this data and so we write $S=S(E,\sigma)$. A key fact is that 
$\coh(E)\simeq \on{qgr}(S/gS)\subset \qgr S$ for  an element $g\in S_3$ that is
unique up to scalar multiplication
 and this inclusion  has a left adjoint of
``restriction to $E$.''  When $\qgr S\simeq\on{coh}(\boldp^2)$, there is no
canonical choice of $E$.  In order to incorporate the 
categories $\qgr S\simeq\coh(\boldp^2)$ into the same framework, we
identify such a category $\qgr S$ with $\coh(\boldp^2)$, 
choose any cubic curve $E\subset \boldp^2$, and set $\sigma=\operatorname{Id}$;
thus $S(E,\sigma)=k[x,y,z]$ in this case.

It follows from the classification that there is a rich
supply of noncommutative planes (see Remark~\ref{as-cases}). In particular the 
plane $\mathbf P^2_\hbar$ associated to the Weyl algebra equals $\qgr U$,
where  $U=\mathbb C\{x,y,z\}/(yx-xy-z^2, z\ \text{central})$ is an algebra
$S(E,\sigma)$  for which $E=\{z^3=0\}$ is the triple line at infinity. 
 Perhaps
the most interesting and subtle algebra in $\ASthree$  is the generic example: 
the  \emph{Sklyanin algebra}  $\skl(E,\sigma)=S(E,\sigma)$ which is determined
by a smooth elliptic curve $E$ and an  automorphism $\sigma$ given by
translation under the group law. One  can regard $\qgr
\skl$ as an elliptic deformation  of $\coh(\mathbf P^2)$ and so we call
$\qgr \skl$ an \emph{elliptic quantum plane}.

As has been mentioned, noncommutative projective planes have  all the basic
properties of $\coh(\mathbf P^2)$ and therefore admit natural definitions 
of vector bundles and torsion-free sheaves as well  as invariants like Euler
characteristics and Chern classes. In particular, \emph{a vector bundle} in 
$\qgr S$ is just the image of a reflexive graded $S$-module (see
Definition~\ref{vectorbundledefn} for a more homological definition in
 $\qgr S$). A \emph{line bundle} is then a vector
bundle of rank one. For any torsion-free  module $\cL\in \qgr S$ of rank one, 
a unique shift $\cL(n)$ of $\cL$ has \emph{first Chern class $c_1=0$} (see
page~\pageref{chern-defn}). \emph{The Euler characteristic} of 
$\cM\in \qgr S$ is defined just
as in the commutative case: $\chi(\cM)=\sum (-1)^j\dim_k \on{H}^j(\cM)$, where 
$\on{H}^j(\cM)= \on{Ext}^j_{\qgr S}(S,\cM)$.

 We can now describe our main results on the structure of rank 1 modules  in
$\qgr S$. Their significance will be described later in the
introduction. 
 
\begin{thm}[Theorem~\ref{connectedness1}]\label{mainthm1}
Let $\qgr S$, for $S=S(E,\sigma)$, be a noncommutative projective plane.
 Then
\begin{enumerate}
\item There  is a smooth,   projective, fine moduli space
$\M^{ss}_S(1,0,1-n)$ for 
 rank one torsion-free  modules in $\qgr S$
with $c_1=0$ and $\chi=1-n$. Moreover, $\dim \M^{ss}_S(1,0,1-n)=2n$.
\item $\M^{ss}_S(1,0,1-n)$ has a nonempty open 
subspace $(\boldp_S\smallsetminus E)^{[n]}$ parametrizing  modules whose 
restriction to $E$ is torsion-free. When $|\sigma|=\infty$, 
$(\mathbf{P}^2\smallsetminus E)^{[n]}$
parametrizes line bundles with $c_1=0$ and $\chi=1-n$.
\end{enumerate}
\end{thm}

\begin{remark} (1) When $\cM\in\qgr S$ is torsion-free of
 rank one and $c_1(\cM)=0$, 
the formula for 
$\chi(\cM)$ simplifies to  $\chi(\cM)=1-\dim_k\on{H}^1(\cM(-1))$
(see Corollary~\ref{monad6}).

(2) We emphasize that the moduli spaces in this theorem (and all
other results of this paper) are schemes in the usual sense; they are
constructed as  GIT quotients of subvarieties of a product of 
Grassmannians---see Theorem~\ref{semistability via GR}.\end{remark}

Most of the algebras in
$\ASthree$ (including the homogenized Weyl algebra $U$, the Sklyanin algebra
$\skl$ and the various quantum spaces defined by three generators $x_i$ which
$q$-commute $x_ix_j=q_{ij}x_jx_i$) occur in families that include the 
commutative polynomial ring. For these algebras, the moduli spaces behave well
in families (see Theorem~\ref{smooth family-intro}, below) 
and we can improve Theorem~\ref{mainthm1}
as follows:
\begin{thm}[Theorem~\ref{connectedness2}]\label{mainthm13}
Let $\cB$ be a smooth irreducible curve defined over $k$ and
let $S_\cB=S_{\cB}(E,\sigma)\in \ASthree$ be a flat family of algebras 
such that $S_p=k[x,y,z]$ for some point $p\in \cB$. Let $S=S_b$ for any
 $b\in \cB$. Then
\begin{enumerate}
\item $\M^{ss}_S(1,0,1-n)$ and $(\mathbf{P}_S\smallsetminus E)^{[n]}$
are  irreducible and hence connected.
\item
  $\M^{ss}_S(1,0,1-n)$ is a deformation of the Hilbert scheme 
$(\mathbf{P}^2)^{[n]}$, with its subvariety
$(\boldp_S\smallsetminus E)^{[n]}$ being a deformation of 
$(\mathbf{P}^2\smallsetminus E)^{[n]}$.
\end{enumerate}
\end{thm}
For the elliptic quantum planes, we get in addition:

\begin{thm}[Theorem~\ref{poisson
structure}]\label{mainthm11} Assume that $S=\skl(E,\sigma)$ is a Sklyanin 
algebra defined over $\mathbb C$.  Then
  $ \M^{ss}_S(1,0,1-n)$ admits a Poisson structure which 
restricts to a holomorphic symplectic structure on 
$(\boldp_S\smallsetminus E)^{[n]}$.
\end{thm}

When $S=\skl(E,\sigma)$ for an automorphism $\sigma$ of infinite order,   de
Naeghel and Van den Bergh \cite{dNvdB} have independently  used quiver
varieties to obtain a  geometric  description of the set of isomorphism
classes of line bundles in $\qgr S$.   One advantage of their approach is that
it  shows that, in stark contrast to the commutative case,
$(\boldp_S\smallsetminus E)^{[n]}$ is  affine.

Just as one can restrict modules from $\qgr S$ to $E$,  so 
there is   a natural functor of
 ``restriction to $\on{Proj}(S)\smallsetminus E$.'' 
Formally, this means
  inverting  the canonical  element $g\in S(E,\sigma)_3$ defining $E$
to obtain the localized
algebra $A(S) = S[g^{-1}]_0$ that can be regarded as a noncommutative
deformation of the ring of functions on $\boldp^2\smallsetminus E$.
The restriction map from $\qgr S$ to $A(S)\on{-mod}$ takes a module $M$
to $M[g^{-1}]_0$.
When $S$ is the algebra $U$,
$A(U)$ is just the Weyl algebra $A_1$ and so the next result gives a natural 
  analogue of the results of \cite{CH, BW1, BW2, BGK}
   discussed at the beginning of the introduction.

\begin{thm}[See Theorem~\ref{three-to-one}]\label{three-to-one-intro}
Let $S=\on{Skl}(E,\sigma)$. Then  the  closed points  
of \begin{equation}\label{disjoint union}
\coprod_{s\in{\mathbf Z}/3{\mathbf Z}}\;\; 
\coprod_{n\geq 0}\;\;(\boldp_S\smallsetminus E)^{[n]}
\end{equation}
are in bijective correspondence with the   isomorphism 
classes
of finitely generated rank 1 torsion-free 
$A(S)$-modules via the map $M\mapsto M[g^{-1}]_0$. 
\end{thm}
The appearance of ${\mathbf Z}/3{\mathbf Z}$ in
  Equation~\ref{disjoint union} corresponds to  the appearance of 
$\on{Pic}(\boldp^2\smallsetminus E) = {\mathbf Z}/3{\mathbf Z}$ 
in the commutative classification.  The bijection of 
Theorem~\ref{three-to-one-intro}
is {\em not} an isomorphism of moduli, nor should one hope for
it to be one---even in the commutative
case the moduli of sheaves on affine varieties are not well behaved.
One should note that the analogue of Theorem~\ref{three-to-one-intro} does not 
hold for all noncommutative planes: Proposition~\ref{eg-prop}
gives a counterexample for a ring of $q$-difference operators.

While Theorem~\ref{three-to-one-intro} holds for all values of the 
automorphism $\sigma$, there are really two distinct cases to the theorem. 
When $|\sigma|=\infty$, $A(S)$ is a simple hereditary ring and so the theorem 
classifies projective rank one modules, just as the work of Cannings-Holland
et al. classified projective $A_1$-modules. However, when $|\sigma|<\infty$,
$A(S)$ is an Azumaya algebra of dimension two and so the classification
includes rank $1$ torsion-free modules that are not projective.

\subsection{General Results and Higher Rank}\label{statement of results} The
moduli results from Theorems~\ref{mainthm1}, \ref{mainthm13}
 and \ref{mainthm11}  also have
analogues for modules of higher rank that  mimic the classical results for
vector bundles and torsion-free sheaves on $\mathbf{P}^2$. As in the
commutative  case, one has a natural notion of 
(semi)stable modules (see 
page \pageref{stability definition for modules}) and we prove results that are 
direct analogues of the commutative results as described, for example, in
\cite{Drezet-Le Potier, Le Potier} and \cite{OSS}.

\begin{thm}[Theorem~\ref{projective moduli spaces}]\label{firstthm}
 Let $\qgr S$ be a noncommutative projective plane and 
fix $r\geq 1$, $c_1\in{\mathbf Z}$, and $\chi\in{\mathbf Z}$. 
\begin{enumerate}
\item
There is a projective coarse moduli space
$\M^{ss}_S(r,c_1,\chi)$ for 
geometrically semistable torsion-free 
 modules in $\qgr S$ of
rank $r$, first Chern class $c_1$ and Euler characteristic $\chi$.  
 
\item
$\M^{ss}_S(r,c_1,\chi)$ has as an 
open subvariety the moduli space $\M_S^s(r,c_1,\chi)$ for  geometrically 
stable  modules in $\qgr S$.
\end{enumerate}
\end{thm}

As the next result shows, the moduli spaces $\M^{s}_S(r,c_1,\chi)$ 
also behave well in families.

\begin{thm}[Theorem~\ref{smooth family}]\label{smooth family-intro}
Suppose that $S=S_\cB$ is
a flat family of algebras in $\ASthree$
  parametrized by a   $k$-scheme $\cB$.
 Then there is a quasi-projective
 $\cB$-scheme $\M_{S}^s(r,c_1,\chi)\rightarrow \cB$ 
that is smooth over $\cB$,
and the fibre of which over $b\in\cB$ is
precisely $\M_{S_b}^s(r,c_1,\chi)$.  In particular, 
$\M^s_S(r,c_1,\chi)$ is smooth when $S\in \ASthree(k)$.
\end{thm}

In the case of the Weyl algebra, our methods give an easy proof 
that the parametrizations of 
Berest-Wilson \cite{BW1} and Kapustin-Kuznetsov-Orlov \cite{KKO} 
are really fine moduli spaces. The paper \cite{KKO}
constructs
a variety $V/\!\!/\on{GL}(n)$ together with a bijection between the
set of points of $V/\!\!/\on{GL}(n)$ and
the set of isomorphism classes of framed vector bundles in $\boldp^2_{\hbar}$
of rank $r$ and Euler characteristic $1-n$ (this generalizes the bijection 
constructed by \cite{BW1} for rank $1$).

\begin{prop}[Proposition~\ref{moduli for others}]\label{moduli for others-intro}
The variety $V/\!\!/\on{GL}(n)$ is 
 a fine moduli space for  framed vector bundles of rank $r$
and  Euler characteristic $\chi=1-n$ in $\boldp^2_{\hbar}$.
This isomorphism induces the bijections of \cite{BW1, KKO}.
\end{prop}

\subsection{Deformations and Poisson Structures} As we  have seen, many
of the noncommutative projective planes can be regarded as deformations of
$\coh(\mathbf P^2)$ and so, by  Theorem~\ref{mainthm13}, they determine
deformations of the Hilbert scheme of points on $\mathbf P^2$. We expand on
this observation in this subsection. For simplicity we just the discuss the
case of a Sklyanin algebra $S=\skl(E,\sigma)$ defined over $\mathbb C$. We
first note that Theorem~\ref{mainthm11} also generalizes to higher ranks:

\begin{thm}[Theorem~\ref{poisson structure}]\label{poisson-intro}
Let $S=\skl(E,\sigma)(\mathbb C)$.
\begin{enumerate}
\item The moduli space $\M_S^s(r,c_1,\chi)$ admits a Poisson structure.  
\item Fix a vector bundle $H$ on $E$ and let $\M_H$ denote the smooth
locus of the locally closed subscheme of $\M_S^s(r,c_1,\chi)$ parametrizing
modules $\cE$ that satisfy $\cE|_E \cong H$.  Then the Poisson structure
of $\M_S^s(r,c_1,\chi)$ restricts to a symplectic structure on $\M_H$.
\end{enumerate}
\end{thm}

The elliptic curve $E$ is the zero locus of a section $s$ of 
$\cO(3)=\Lambda^2T_{\mathbf P^2}$ that is unique up to scalar 
multiplication; upon restriction to $\boldp^2\smallsetminus E$, 
$s^{-1}$ is an algebraic symplectic  structure. 
The Poisson structure $s$ also induces a
Poisson structure on  the Hilbert scheme  $(\mathbf P^2)^{[n]}$ that
restricts to give
a symplectic structure on $(\mathbf P^2\smallsetminus E)^{[n]}$
(see \cite{Beauville}). 
This Poisson structure on $\mathbf P^2$  determines a  deformation of the polynomial
ring $k[x,y,z]$ that is precisely the Sklyanin algebra $\skl(E,\sigma)$, with
$\sigma$ corresponding to the deformation parameter
\cite[Appendix~D2]{Od}.  It is natural to  hope
that the noncommutative deformations of $\boldp^2$ induce deformations
 of  $(\boldp^2)^{[n]}$ and
$(\mathbf P^2\smallsetminus E)^{[n]}$ that also carry Poisson and
symplectic structures. 
The point of Theorems~\ref{mainthm13}(2) and \ref{mainthm11} is that 
$(\mathbf P^2_S)^{[n]}$ and  $(\mathbf P^2_S\smallsetminus E)^{[n]}$ provide
just such a deformation.

The relationship between the deformed Hilbert scheme of points $(\mathbf
P^2_S)^{[n]}$ and collections of points in $\mathbf P^2$ or $\qgr S$ is rather 
subtle (as was also the case with the original work on modules over the Weyl
algebra).  For simplicity, we will explain this when   $|\sigma|=\infty$ and we
regard $(\mathbf P^2_S)^{[n]}$ as parametrizing  ideal sheaves $\cI \subset
\cO_{\mathbf P^2}$ for which  $ \cO_{\mathbf P^2}/\cI $ has length $n$. Now
the  only  simple objects in  $\qgr S$ are the modules corresponding to points
on $E$.\footnote{From the point of view of quantization, this corresponds to the fact
that the Poisson structure on $\boldp^2$ is nondegenerate over
$\boldp^2\smallsetminus E$.}  By mimicking the commutative procedure, it is
therefore easy to find the modules in $(\mathbf P^2_S)^{[n]}$ corresponding to 
collections of points on $E$, but the modules parametrized by 
$(\boldp^2_S\smallsetminus E)^{[n]}$ are necessarily  more subtle (compare
Corollary~\ref{stable in rank 1} with Proposition~\ref{existence}). In
essence,  the proof of  Theorem~\ref{mainthm1}(2) shows that  the ideal sheaves
parametrized by $(\boldp^2\smallsetminus E)^{[n]}$ can be  constructed as the
cohomology of certain complexes (monads) that  deform to complexes in $\qgr
S$.  The  cohomology of such a deformed complex cannot correspond to a
collection of points in $\qgr S$ simply because these points do not exist. It
is therefore forced to be an interesting line bundle.  This subtlety is also
reflected  in Theorem~\ref{three-to-one-intro} since one is parametrizing
projective $A(S)$-modules rather than annihilators of finite dimensional
modules.

The Poisson structure of  $\M_S^s(r,c_1,\chi)$ in
Theorem~\ref{poisson-intro}
 is related to work of Feigin-Odesskii \cite{FO2} on their
generalizations $Q_{n,k}(E,\sigma)$ of higher dimensional Sklyanin 
algebras.
The classical limit of the  $Q_{n,k}(E,\sigma)$ induces a Poisson structure  on
the moduli spaces of certain vector bundles over  $E$ \cite[Theorem~1]{FO2}.
This structure is in turn a special case of  a Poisson pairing on the moduli
space of  $P$-bundles on $E$, where $P$ is a parabolic subgroup  of a reductive
group.  Now, our moduli spaces may also be identified with moduli spaces of
filtered $E$-vector bundles and  the method we use to obtain our Poisson
structure on $\M_s^s(r,c_1,\chi)$ is to use Polishchuk's generalization
\cite{Polishchuk} of the Feigin-Odesskii construction (see 
Subsection~\ref{construction of poisson structure}). Given this connection with
\cite{FO2} and the $Q_{n,k}(E,\sigma)$, it would be interesting to know whether
our Poisson structure  induces similarly interesting noncommutative
deformations  of $\M_s^s(r,c_1,\chi)$.

Finally, it would be interesting to relate our  moduli spaces to integrable
systems. We  make a small step in this direction in 
Subsection~\ref{equivariant higgs section} by explaining how the symplectic
leaves of $\M^s_S(r,c_1,\chi)$ may be related to parameter spaces for Higgs
bundles with values in the centrally extended current group of  \cite{EF}.

\subsection{Methods} 
One key element of the  paper  is the thoroughgoing use of
generalizations 
of the  cohomological tools from commutative algebraic
geometry, primarily Cohomology and Base Change.
These may be found in a form sufficient for our purposes 
in Section~\ref{section-base}.  In particular, it is through these methods that
we are able to prove formal moduli  results.

Our treatment of semistability and linearization of the group
actions in Sections~\ref{section3} and \ref{section5} follows
the outline of Drezet and Le Potier \cite{Drezet-Le Potier, Le Potier}
for $\mathbf P^2$.
 Their techniques require some 
modification, however,  because of the shortage of ``points'' on our
noncommutative surfaces.  
This proof has several steps. First, we construct a version of the Beilinson
spectral sequence. The version  used here is different than that appearing in
\cite{Le, BW1, KKO, BGK}, since it uses the \v{C}ech complex from \cite{VW2} to
avoid the problem that  tensor products of modules over a noncommutative ring 
are no longer modules. This spectral sequence is then used to show that the
moduli space of semistable modules in $\qgr S$ is equivalent to that for 
semistable monads and then for semistable Kronecker complexes.   Here, \emph{a
Kronecker complex} is a complex
  $S(-1)^a \xrightarrow{\alpha} S^b\xrightarrow{\beta}
S(1)^c$ in $\qgr S$.  This complex is \emph{a monad} if $\alpha$ is injective and
$\beta$ is surjective. The semistable Kronecker complexes can be described by
purely commutative  data and are easy to analyse by  standard techniques of GIT
quotients. An important fact is that this all works in families, which is why
we are able to  construct our moduli spaces.

The construction of a Poisson structure and the relation with Higgs bundles
are  analogs of results of  \cite{Polishchuk} and \cite{GM}, respectively.

Since we hope that the paper may be of interest to readers of varied
backgrounds, we have included some details in proofs that we imagine
will be useful for some readers but unnecessary for others.

\subsection{Further Directions}
 
A natural question  that is not answered in the present work
concerns the metric geometry of our
 deformations of   Hilbert schemes.  The
plane cubic complement $\boldp^2\smallsetminus E$ admits a complete hyperk\"ahler
metric (we are grateful to Tony Pantev and Kevin Corlette for explaining this
to us) and one imagines that the Hilbert schemes of points also admit
such metrics; if this is true, one would like to know the
relationship between our deformations and the twistor family for the
hyperk\"ahler metric on the Hilbert scheme.

In a more speculative direction, one would like to have an interpretation,
 parallel to that in \cite{BW1,BW2, BGK, BGK2}, of the moduli spaces
$\M^{ss}_S(r,c_1,\chi)$ in terms of integrable systems.  
Section~\ref{current groups} makes a start in this direction, but this is
certainly an important direction for further work to which we hope to
return.

As was remarked earlier, the Weyl algebra and its homogenization
have been applied to questions in string theory \cite{KKO,Nekrasov-Schwarz}
and so it would be interesting to understand whether the other quantum planes 
have applications in this direction \cite[Section~VI.B]{DN}.
The Sklyanin algebra  
appears in general marginal deformations of $N=4$ super Yang-Mills theory
(see \cite[Equations~2.5--2.7 and Section~4.6.1]{BJL}) 
and, in the terminology of that paper,
Theorems~\ref{mainthm1} and \ref{mainthm13}
 can be interpreted as a parametrization of 
space-filling D-branes  (see \cite[Section~6]{BL}).

\subsection{Acknowledgments}
The authors are grateful to David Ben-Zvi, Dan Burns, 
Kevin Corlette, Ian
Grojnowski, and Tony Pantev for many helpful conversations.  
In particular, Tony Pantev  introduced the authors to
noncommutative deformations in the 
context of moduli spaces.

\section{Background Material}\label{section2}

In this section we collect the basic definitions and 
 results from the literature that will be used throughout the paper.

Fix a   base field $k$ and a commutative $k$-algebra $C$.
Noncommutative projective geometry 
\cite{Stafford, SVdB} is concerned with the study of connected graded (cg)
$k$-algebras  or, more generally cg $C$-algebras, where a graded $C$-algebra
$S=\bigoplus_{i\geq 0} S_i$ is called \emph{connected graded}\label{cg-defn} if
$S_0 =C$ is a central subalgebra. The $n$th {\it Veronese
ring}\label{veronese-defn} is defined to be the ring $S^{(n)}=\bigoplus_{i\geq
0} S_{ni}$, graded by $S^{(n)}_i=S_{ni}$. Write $\Module S$ for the category of
right $S$-modules and $\Ggr S$ for the category of graded right $S$-modules,
with homomorphisms $\Hom(M,N)=\Hom_S(M,N)$  being graded homomorphisms 
 of degree zero.  Given $M=\bigoplus_{i\in \mathbb
Z}M_i$, the {\it shift} $M(n)$\label{shift-defn} of $M$ is the graded module 
$M(n)=\bigoplus M(n)_i$ defined by $M(n)_i=M_{i+n}$ for all $i$. The other hom
group is  $\underline{\Hom}_S(M,N) = \bigoplus_{r\in \mathbb Z} \Hom(M,N(r))$,
with derived functors $\underline{\Ext}^j(M,N)$.\label{ext-defn}
 If $M$ is
finitely generated, then  $\underline{\Ext}^j(M,N)=\Ext^j_{\Module S}(M,N).$ 
Similar conventions apply to  Tor groups.

A module   $M=\bigoplus_{i\in \mathbb Z}M_i\in \Ggr S$ is called
\emph{right bounded}\label{bounded-defn}
if $M_i=0$ for  $i\gg 0$. The full Serre subcategory of $\Ggr S$\label{gr-defn}
generated by the right bounded  modules is denoted  $\Tors S$ with quotient
category  $\Qgr S=\Ggr S/\Tors S$.  If $S$ happens to be noetherian (which will
almost always be the case in this paper)  write $\module S$, $\gr S$ and $\qgr
S = \gr S/\tors S$ for the subcategories of noetherian objects in  these three
categories. Observe that $\tors S$ is just the category of finite-dimensional
graded modules. Similar definitions apply for left modules and we write the
corresponding categories as $S$-${\mathrm{gr}}$, etc.

We now turn to the definitions for 
the  algebras of interest to us: Sklyanin algebras and, more
generally, Artin-Schelter regular rings. 
Fix a (smooth projective) elliptic  curve $\iota: E \hookrightarrow
 {\boldp}^2$ with corresponding 
line bundle $\cL=\iota^*({\mathcal O}_{{\boldp}^2}(1))$ 
of degree $3$. Fix an automorphism
$\sigma\in {\mathrm{Aut}}(E)$ given by translation under the group law and 
denote the graph of $\sigma$ by $\Gamma_{\sigma}\subset 
E\times E$. 
 If $V= \HH^0(E,\cL)$,  there is a $3$-dimensional space 
\begin{equation}\label{relations}
\cR = \HH^0\big(E\times E, (\cL\boxtimes\cL)(-\Gamma_{\sigma})\big)\subset
 \HH^0\big(E\times E, \cL\boxtimes \cL\big) = V\otimes V.
\end{equation}

\begin{defn}\label{skl-def}
The {\em $3$-dimensional Sklyanin algebra}  is the algebra
$$S=\skl = \skl(E,{\mathcal L},\sigma) = T(V)/(\cR),$$ where $T(V)$
denotes the tensor algebra on $V$. 
\end{defn}
\noindent

 When $\sigma$ is
 the identity, $\cR=\Lambda^2 V$ and so 
 $\skl(E,{\mathcal L},\mathrm{Id})$ is just
 the polynomial ring $k[x,y,z]$. One can therefore regard
  $\skl(E,{\mathcal
 L},\sigma)$ as a deformation of that polynomial ring and this deformation
is flat \cite{TVdB}. 
 One may also write $\skl(E,{\mathcal L},\sigma)$ as the $k$-algebra with 
 generators $x_1,x_2,x_3$ and relations:
\begin{equation}\label{skly-def2}
ax_ix_{i+1} + bx_{i+1}x_i + cx_{i+2}^2=0\qquad i=1,2,3 \mod 3,
\end{equation}
where $a,b,c\in k^*$ are any scalars satisfying 
$(3abc)^3\neq (a^3+b^3+c^3)^3$ (see the introduction to \cite{ATV1}).

Basic properties of $S=\skl(E,{\mathcal L},\sigma)$
 can be found in \cite{ATV1, ATV2} and  are summarized
in \cite[Section~8]{SVdB}. In particular, it is 
one of the most important examples of an AS regular ring defined as follows:

\begin{defn} \label{as-defn} An \emph{Artin-Schelter (AS) regular ring of
dimension three} is a connected graded algebra $S$ satisfying 
\begin{itemize} \item
$S$ has global homological dimension $3$, written ${\mathrm{gldim}}\; S=3$;
\item $S$ has  Gelfand-Kirillov dimension $3$, written $ \operatorname{GKdim} S
= 3$; and   
\item  $S$ satisfies the \emph{AS Gorenstein condition}: 
$\underline{\Ext}^i(k,S)=\delta_{i,3} k(\ell)$, for some $\ell$.  \end{itemize}
Let $\ASthree=\ASthree(k)$\label{ASthree-defn} denote the class of AS regular
$k$-algebras of dimension three with Hilbert series $(1-t)^{-3}$ and with
$\ell=-3$. \end{defn}

The algebras in $\ASthree$
have been classified in \cite{ATV1} and we will list the generic 
examples at the end of this section (see  
Remark~\ref{as-cases}).
 In particular, $\skl\in \ASthree$. As is
justified in \cite[Section~11]{SVdB},  the category $\qgr S$ for  $S\in
\ASthree$  can be regarded as a (nontrivial) noncommutative projective plane 
and most of our results work for $\qgr S$ for any 
$S\in \ASthree$. 

There are two basic techniques for understanding $S$ and $\Qgr S$
for $S\in \ASthree$. 
First, assume that $S\in \ASthree$ but that $\Qgr S$ is \emph{not} 
equivalent to $\operatorname{Qcoh}(\boldp^2)$,\label{qcoh}
the category of quasi-coherent
sheaves on $\boldp^2$. 
By \cite[Theorem~II]{ATV1} there exists, up to a scalar multiple,
 a canonical normal element $g\in S_3$ (thus, $gS=Sg$).  
The factor ring $S/gS$ is isomorphic to the
twisted homogeneous coordinate ring $B=B(E,\cL,\sigma)$\label{tw-defn} for a cubic curve
 $\iota: E\hookrightarrow \boldp^2$, the line bundle
 $\cL=\iota^*\cO_{\boldp^2}(1)$ of degree three
  and an automorphism $\sigma$ of the scheme $E$.
   The details of this construction can be found in 
\cite{AVdB} or  \cite[Section~3]{SVdB}, but the essential
properties are  the following:  write $\cN^\tau$ for the
 pull-back of a sheaf $\cN$ along an automorphism $\tau$ of $E$ and set
$
\cL_n=\cL\otimes\cL^{\sigma}\otimes\cdots\otimes\cL^{\sigma^{n-1}}.$
 Then
$B=\bigoplus_{n\geq 0} B_n$, where
$B_n=\HH^0(E,\cL_n) $.
The algebra $S$ is  determined by the data $(E,\mathcal{L},\sigma)$ and 
so we write $S=S(E,\mathcal{L},\sigma)$. When $S=\skl$, $E$ and $\sigma$
coincide with the objects used in Definition~\ref{skl-def}. 
In fact, $\cL$ is superfluous to the
classification, both  for the Sklyanin algebra and  in general.  
However,  this sheaf appears in many of our results and so we have
included it in the notation.

We are mostly interested in $\Qgr S$ rather than $S$ itself, 
and it is useful to be able to use the curve $E$. However 
there is no canonical 
embedded curve when  $\Qgr S\simeq \on{Qcoh}(\mathbf{P}^2)$.  
Thus, in order to have a uniform approach, we let 
$\ASthree'$\label{ASthree-other} denote the algebras 
 $S\in \ASthree$ for which $\Qgr S$ is \emph{not} 
equivalent to $\operatorname{Qcoh}(\boldp^2)$, together 
with $S=k[x,y,z]$. For $S=k[x,y,z]$ we fix any nonzero 
homogeneous element $g\in S_3$,
 set $E=\on{Proj}(S/gS)$ and, somewhat arbitrarily, write $S=S(E,\cL,\on{Id})$,
 where $\cL $ is the restriction of $\theo_{\mathbf P^2}(1)$ to $E$.
  Thus, for each $S\in \ASthree$ there exists 
$S'\in \ASthree'$ with $\Qgr S\simeq \Qgr S'$.

Let $S=S(E,\cL,\sigma)\in \ASthree'$. 
The factor ring $B=B(E,\cL,\sigma)$ is noetherian. More significantly,
 $\on{Qcoh}(E)\simeq \Qgr B$
via the map $$\xi: \cF \mapsto \bigoplus_{i\geq 0} 
 \HH^0(E, \cF\otimes_{\cO_X}\cL_n)\qquad\text{for}\quad\cF \in \on{Qcoh}(E).$$
This induces a map 
$\rho: \Qgr S\to \Qgr S/gS = \Qgr B \simeq  \operatorname{Qcoh}(E)$ and,
 mimicking geometric notation, we write 
\begin{equation}\label{defn-restriction}
 \cM|_{E}=\rho(\cM)\in  \operatorname{Qcoh}(E)\qquad \mathrm{for}\quad
 \cM\in \Qgr S.
 \end{equation}
When $S=\skl$,  the element $g$ is   central
 and $B$ is   a domain.
For other $S\in \ASthree$,  it can happen
that $g$ is not central or that $E$ is not integral. 
In the latter case $B$ will not
be a domain. 

The second important technique in the study of $\qgr S$, for $S\in \ASthree$,
 is to use cohomological techniques modelled on the classical case.
  This works well is because  $S$ satisfies the $\chi$ condition of 
\cite[Definition~3.7]{AZ}, defined as follows: 
 a cg $C$-algebra $R$ satisfies $\chi$\label{chi-defn} if 
$\Ext_{\Module R}^j(R/R_{\geq n}, M) $ is a finitely generated $C$-module for 
all   $M\in \Module R$, all $j\geq 0$, and all $n\gg 0$.
Write  $\pi$    for  the natural  projections 
 $\Ggr S \to \Qgr S$ and $\gr S\rightarrow \qgr S$ and set $\pi(S)=\theo_S$.
 One
may pass back from $\Qgr S$ to $\Ggr S$ by the ``global sections'' functor:
\begin{equation}\label{gamma-fn}
\Gamma^*(\cM) = \bigoplus_{n\geq 0}{\Hom}_{\Qgr S}(\theo, \cM(n)),
\qquad{\mathrm{for} }\ \cM\in\Qgr S.
\end{equation}
One then has an adjoint pair $(\pi, \Gamma^*)$. 
We will occasionally use the module
$\Gamma(\cM) = \bigoplus_{n\in \mathbb Z}{\Hom}_{\Qgr S}(\theo, \cM(n))$.
This module has the disadvantage that   
$\Gamma(\pi(M))$ need not be finitely generated, but 
 $\Gamma(\cM)/\Gamma^*(\cM)$ is at least right bounded.

As in the commutative situation, we will write the derived functors of 
$\HH^0(\cM)={\Hom}_{\Qgr S}(\theo, \cM)$,\label{Hi-defn} for  
$\cM\in\Qgr S$, as 
$$\HH^i(\cM)=\HH^i(\Qgr S, \cM)=\Ext^i_{\Qgr S}(\cO,\cM).$$ 
One also has the analogous objects for $\qgr B$ and 
one should note that 
there is no confusion in the notation since 
$\HH^i(\Qgr S, \cM)=\HH^i(E, \cM)$ for  $i\geq 0$ and  $\cM\in \qgr B$
(see \cite[Theorem~8.3]{AZ}).

The significance of the $\chi$ condition is that the functors $\HH^i$
 are particularly well-behaved. 
In particular,  \cite[Theorems~8.1 and
4.5 and Corollary, p.253]{AZ} combine to show:
\begin{lemma}\label{chi-lemma} Let $S\in \ASthree$. Then:
\begin{enumerate}
\item\label{chi-one} $\pi\Gamma^*(\cM)=\cM$ for any $\cM\in \Qgr S$.
\item \label{chi-two} If $M\in \gr S$, then
$\Gamma^*(\pi M)=M$ up to a finite dimensional vector space.
\item \label{chi-three} More precisely, if 
$M=\bigoplus_{i\geq 0}M_i\in \gr S$ and $F$ is the
maximal finite-dimensional submodule of $M$, then $\Gamma^*\pi M$ is the
maximal positively graded essential extension of $M/F$ by a 
finite-dimensional module.
\item
$\Gamma^*(\pi(S))=S=\Gamma(\pi(S))$ and 
$\Gamma^*(\pi(S/gS))=S/gS=\Gamma(\pi(S/gS))$ whenever $S\in \ASthree'$.
\item $\HH^1(\theo_S(j))=0$ for all $j\in \mathbb Z$
and $\HH^2(\theo_S(k))=0$ for all $k\geq -2$.
\item $S$ has \emph{cohomological dimension $2$} in the sense that 
$\HH^i(\cM)=0$ for all $\cM\in \Qgr S$ and $i>2$.\label{cohom-defn}
\qed
\end{enumerate}
\end{lemma}

There is also an analogue of Serre duality:
 \begin{prop}[Serre Duality]\label{serre-duality}
Let $S\in \ASthree$ and $\M, \N\in \qgr S$. Then
\bd
\Ext^i(\M,\N) \cong \Ext^{2-i}(\N,\M(-3))^*.
\ed
\end{prop}
\begin{proof} By 
\cite[Theorem~2.3]{YZ} and \cite[Theorem~8.1(3)]{AZ} one has 
  $$\Ext^i(\cM,\cO(-3))\cong \HH^{2-i}(\cM)^*=\Ext^{2-i}(\cO, \cM)^*$$
  for all $\cM\in \qgr S$. 
This version of Serre duality implies that 
$\Ext^2(\M,\theo(-n))\neq 0$ for $n\gg 0$, 
and  $\Ext^i(\M, \theo(-n)) = 0$ for $n\gg 0$ and $i<2$.  
Also, since $S$ has finite global homological dimension, 
Lemma~\ref{chi-lemma}(6) implies that $\Ext^j(\cM,\cQ)=0$ 
for all $\cQ\in \Qgr S$ and $i>2$.  Set
$
\cL=\pi\underline{\Ext}^2(\M,\theo(-3))^* 
  \cong \M(-3)$. 
Then all the hypotheses of  \cite[Theorem~2.2]{YZ} are satisfied and that result
implies that 
$
\Ext^i(\M,\cN)\cong \Ext^2(\cN,\cL)^* \cong \Ext^{2-i}(\cN,\cM(-3))^*.
$
\end{proof}

The next two lemmas give some standard facts about the relationship 
between modules over $S$ and $E$ for which we could not find a 
convenient reference.

\begin{lemma}\label{shifts}  Let $S=S(E,\cL,\sigma)\in \ASthree'$.
  Suppose that 
 $\cM$ is a module in $\qgr S$
and that $j\in {\mathbb Z}$. Then 
$\cM(j)|_E =
(\cM|_E)^{\sigma^{-j}}\otimes_{\cO_E} \cL_{j}^{\tau_{j}}$,
where $\tau_j = \sigma^{-j}$ if $j\geq 0$ 
but $\cL_j=\cL_{-j}^{-1}$ and $\tau_{j} = 1$ 
if $j<0$.
\end{lemma}

\begin{proof} 
Let $\cF\in \coh(E)$, with $\xi(\cF)=\bigoplus_{j\geq 0} \HH^0(\cF\otimes
\cL_j)\in \gr B$.  
By \cite[(3.1)]{SVdB} the shift functor in $\gr B$
 can be reinterpreted in terms of $\theo_E$-modules as:
$$
 \xi(\cF)(n) = \xi\left(
 \sigma^n_*(\cF\otimes \cL_n)\right) = 
 \xi\left(\cF^{\sigma^{-n}}\otimes 
\cL_n^{\sigma^{-n}}\right)\qquad {\rm for\ } n\geq 0.
$$
A simple computation then shows that 
$
\xi(\cF)(-n) =\xi\left(\cF^{\sigma^{n}}\otimes 
\cL_n^{-1}\right)$ for $  n\geq  0.$
Since restriction to $E$ commutes with the shift functor, the lemma
 follows.
  \end{proof}

Let $S\in \ASthree$. A  module  $M\in \Ggr S$ is called {\it torsion-free},
respectively {\it torsion},\label{tf-defn} if no element, respectively every
element, of $M$ is killed by a nonzero element of $S$. A module $\cM\in \Qgr S$
is {\it torsion-free} if $\cM=\pi(M)$ for some torsion-free module $M\in \Gr S$.
By \cite[S2, p.252]{AZ} this is equivalent to 
$\Gamma^*(\cM)$ being torsion-free. 
The term ``torsion'' is also used in the literature  to mean elements of $\Tors
S$, but that will never be the case in this paper. 
A torsion-free module $\cM$ in ($\qgr S$ or $\gr S$) has {\it rank
$m$}\label{rank-defn} if  $\cM$ has Goldie rank $m$; that is, $\M$ contains a
direct sum of $m$,  but not $m+1$,  nonzero submodules. Note that $\cM\in \qgr
S$ is torsion-free of rank one  if and only if $\Gamma^*(\cM)$ is isomorphic
to a shift of  a right ideal of $S$.

\begin{lemma}\label{restriction of torsion-free} 
Let $S\in \ASthree'$ and suppose that
$
0\rightarrow \cF_1\rightarrow \cF_2 \rightarrow \cF_3 \rightarrow 0
$
is a short exact sequence  in $\qgr S$ such that $\cF_3$ is torsion-free.
 Then the restriction of 
this sequence to $E$ remains exact.
\end{lemma}
\begin{proof}
Since $\Gamma^*$ is left exact,
we have an exact sequence of   $S$-modules
\begin{equation}\label{restr22}
0\rightarrow F_1\rightarrow F_2 \rightarrow 
F_3 \rightarrow 0,
\end{equation}
where $F_i=\Gamma^*(\cF_i)$ for $i=1,2$ 
and $F_3$ is the image of $\Gamma^*(\cF_2)$ in $\Gamma^*(\cF_3).$ 
By \cite[S2, p.252]{AZ} $F_3$ is torsion-free and, by \cite[Corollary, p.253]{AZ}, 
$\Gamma^*(\cF_3)/F_3$ is finite dimensional. 
As $F_3$ is torsion-free,  $F_3\otimes gS \hookrightarrow F_3\otimes
S$  and so  $\underline{\mathrm{Tor}}^1(F_3,S/gS)=0$.
Therefore, \eqref{restr22} induces the exact
sequence:
$$0\rightarrow F_1/F_1g\rightarrow F_2/F_2g \rightarrow 
F_3/F_3g \rightarrow 0.$$
Since $\cF_i|_E$ is the image of  $F_i/F_ig$ in $\coh(E)$
 for $1\leq i\leq 3$,
this proves the lemma.
\end{proof}

There is one further algebra that will be considered in this paper.
 Let $S\in \ASthree'$.
Write $S[g^{-1}]$ for the localization of $S$ at the powers of $g$.  
Since $g$ is 
homogeneous of degree $3$,  $S[g^{-1}]$ is still ${\mathbb Z}$-graded.
Let $$A=A(S)=S[g^{-1}]_0\label{A(S)-defn}$$ denote the $0$th graded 
piece of $S[g^{-1}]$. 
Just as $S$
can be thought of as the homogeneous coordinate ring
 of the noncommutative projective plane
$\qgr S$,  
so the algebra $A$ can be
regarded as  the ``coordinate ring of the quantum affine 
variety $\qgr S \smallsetminus E$,''
and one of the aims of this paper is to understand the right ideals of
this ring in terms of  modules in $\qgr S$ and the geometry of 
$\boldp^2\smallsetminus E$.

The next two results collect some  basic facts about
$A=A(S)$. 
 Let $\Lambda_1 =
S_3g^{-1}\subset A$ and note that 
$\Lambda =\{\Lambda_i=\Lambda_1^i\}$ provides a filtration of $A$. 
Since $g$ is a homogeneous normal element,
$S_ig=gS_i$ for all $i$ and so
$\Lambda_i = (S_3g^{-1})^i = S_3^ig^{-i} = S_{3i}g^{-i}.$

\begin{lemma}\label{firstlemma}
Let $S=S(E,\mathcal{L},\sigma)\in \ASthree'$ and set $A=A(S)$. Then:
\begin{enumerate}
\item\label{firstlemmapt1} $A$ is generated by
 $\Lambda_1$ as a $k$-algebra. 
\item $R=S[g^{-1}]$ is strongly graded in the
 sense that  $R_nR_m=R_{n+m}$ for all $m,n$. Moreover,
  $R$ is  projective as a left or right  $A$-module.
\item\label{firstlemmapt3} If $|\sigma | =\infty$ (as an automorphism of 
the scheme $E$) then  $A$ is a simple
hereditary ring.
\item\label{firstlemmapt4}  If $|\sigma|<\infty$ then $A$ is an Azumaya algebra
of global dimension two.
\item\label{firstlemmapt5} Set
 $\cR(A) = \bigoplus \Lambda_1^i g^i \subset S[g^{-1}]$. 
 Then $\cR(A)= S^{(3)}$.
\end{enumerate}
\end{lemma}
\begin{proof} (1) 
By construction, 
$ A= \sum_{i\in\boldz}S_{3i}g^{-i} =
 \sum_{i\in\boldz}\left(S_3g\inv\right)^i=\sum_{i\geq 0}\Lambda_1^i$.

(2) By  \cite[Proposition~7.4]{ATV2}, $R$ is strongly graded.
 Thus  $R_iR_{-i}=R_0=A$ for all $i$ and  the dual basis lemma implies that
 each $R_i$ is  a projective $A$-module.
 
 (3)
 $A$ is simple by \cite[Theorem~I]{ATV2}.
  By \cite[Theorem~A.I.3.4]{NV}
 $\module A\simeq\gr R$ via the map 
   $M\mapsto M\otimes_AR$, for $M\in \module A$. Thus, it suffices to
  show that each  $N\in \gr R$ has homological dimension  
  $\on{hd}(N)\leq 1$.
 This is proved in \cite[Proposition~2.18]{Aj}
   for the Sklyanin algebra and the proof works in general.

   (4)  $A$ is Azumaya by \cite[Theorem~I]{ATV2}.
 Since $S$ has global dimension three, it follows that graded $R$-modules have
 homological dimension at most two. Thus, as in the proof of (3), 
${\mathrm{gldim}}(A)\leq 2$.
The other inequality follows from
 \cite[Corollary~6.2.8 and Proposition~13.10.6]{MR}.
 
 (5) By construction, $\cR(A)$ is a graded ring with 
 $\cR(A)_i=\Lambda_ig^i=S_{3i}g^{-i}g^i=S_{3i}$ for all $i$. Therefore, 
$\cR(A)=S^{(3)}$. 
\end{proof}

\begin{corollary}\label{firstcor} Let $S\in \ASthree'$ with factor ring
$B=S/gS$ and set $A=A(S)$.
Then:
\begin{enumerate}
\item If $g$ is central in $S$, then $S^{(3)}$ is naturally isomorphic to the
Rees ring of $A$ and $B^{(3)}\cong \operatorname{gr}_\Lambda
A=\bigoplus\Lambda_i/\Lambda_{i-1}$.
\item If $E$ is an integral curve then  $\operatorname{gr}_\Lambda
A$  is a domain.
\end{enumerate}
\end{corollary}

\begin{proof} (1) The Rees ring of $A$ is defined to be the graded algebra
$\Rees=\bigoplus \Lambda_it^i $, regarded as a subring of the polynomial ring
$A[t]$.  Since $A[t]\cong A[g]\subset S[g^{-1}]$, clearly the map $t\mapsto g$ 
induces an isomorphism $\Rees\cong \cR(A)\cong S^{(3)}$.  The identity
$B^{(3)}\cong \operatorname{gr}_\Lambda A$
then follows from the observation that  $\operatorname{gr}_\Lambda A\cong
\Rees/t\Rees$.

(2) When $g$ is not central, $\cR(A)$ is isomorphic to a {\it twist} of $\Rees$,
in the   sense of \cite[Section~8]{ATV2}.
 In more detail, conjugation by $g$ defines an
automorphism $\tau $ on $S$ and hence on $A$; thus $\tau(a)=gag^{-1}$.
We may then define a new multiplication on $\Rees$ by 
$a\circ b=a\tau^r(b)$, for $a\in \Rees_r$ and $b\in \Rees_s$. 
It is readily checked that the map $t\mapsto g$ induces an isomorphism 
$(\Rees,\circ)\cong \cR(A)$
and so  $(\Rees/t\Rees,\circ)\cong \cR(A)/g\cR(A)=B^{(3)}$. 
If $E$ is integral then $B$ is a domain \cite{AVdB} and hence so is 
$B^{(3)}=S^{(3)}/gS^{(3)}$. This immediately  implies that 
$\Rees/g\Rees$ is  a domain.
\end{proof}

We end this section by describing  some of the algebras   appearing
in the classification of $\ASthree'$.

\begin{remark}\label{as-cases} According to the classification in 
\cite[(4.13)]{ATV1}, the algebras in $\ASthree$ with 
$\qgr S\not\simeq\coh(\mathbf P^2)$ break into 
seven classes.
For the reader's convenience, we note the generic example of 
the curve $E$, the automorphism $\sigma$, the number of essential parameters
$\mathbf p$ defining the family and 
the defining relations $\{f_i\}$ for a typical algebra
 $R=k\{x,y,z\}/(f_1,f_2,f_3)$ in each class. The letters $p,q,r$ will 
 always denote elements of $k^*$. 
\begin{enumerate}
\item[] {\bf Type $\mathbf A$}: $E$ is elliptic with  $\sigma$ translation,
 $\mathbf p = 2$
and $R=\skl(E,\cL,\sigma)$.
\item[] {\bf Type $\mathbf{B}$}:  $E$ is elliptic with $\sigma$ being
multiplication by $(-1)$ and $\mathbf p = 1$.
Take $\{f_i\}=\{ xy+yx-z^2-y^2, x^2+y^2+(1-p)z^2, zx-xz+pzy-pyz\}$.
\item[] {\bf Type $\mathbf{S_1}$}:  $E$ is a triangle with each component
 stabilized by $\sigma$ and $\mathbf p = 3$. 
 Take $\{f_i\}=\{xy+pyx,yz+qzy,zx+rxz\}$, with $pqr\not=-1$.
\item[] {\bf Type $\mathbf{S_1'}$}:  $E$ is the union of a line and a conic, 
with each component stabilized 
by $\sigma$ and $\mathbf p = 2$. Take $\{f_i\}=\{xy+pyx+z^2,yz+qzy,zx+qxz\}$.
\item[] {\bf Type $\mathbf{S_2}$}:  $E$ is  a triangle with two components
interchanged by $\sigma$ and $\mathbf p = 1$. Take 
$\{f_i\}=\{x^2+y^2,yz+qzx,xz+qzy\}$.
\item[] {\bf Types $\mathbf{E,H}$}:  $E$ is elliptic with
complex multiplication $\sigma$. Here,  $\mathbf p = 0$.
\end{enumerate}
If one takes the example $R$ of  type $S_1$, but with $pqr=-1$, then  one
obtains an algebra $R$ with $\qgr R\cong \coh(\mathbf P^2)$. 
\end{remark}

It is useful to  think of $A=A(\skl)$ as an ``elliptic'' analogue of the first
Weyl algebra $A_1=k\{u,v\}/(uv-vu-1)$,\label{weyl-defn}
 as this will help illustrate the
relationship between the results of this paper and those of  Berest-Wilson
\cite{BW1, BW2} and others.  The  AS regular algebra associated to $A_1$ is 
\begin{equation}\label{U-defn}
U=k\{x,y,z\}/(xy-yx-z^2, z\text{ central}).
\end{equation}
It is easy to check that 
$U\in \ASthree'$ (of Type $S_1'$), with $g=z^3$ and 
 $A_1\cong A(U)$.

There is one significant difference between $A(\skl)$ and the first Weyl algebra
$A_1$. The group of $k$-algebra automorphisms  $\AAut_k(A_1)$  is large
\cite[Th\'eor\`eme~8.10]{Di} and this group plays an
important role in the work of Berest and Wilson \cite{BW1}.
However, the opposite is true of the Sklyanin algebra.
 
\begin{prop}\label{aut-a} Let $S=\skl$ and $A=A(S)$. Write
$\PAut_{gr}(S^{(3)})$ for the group of graded 
$k$-algebra automorphisms $\theta$ of 
$S^{(3)}$  satisfying  $\theta(g)=g$.

If $|\sigma|=\infty$ then $\AAut_k(A)\cong \PAut_{gr}(S^{(3)})$ is a finite
group. \end{prop}

Since this result is not used in the paper, we will omit the proof.

\section{From $A$-Modules to $S$-Modules and Back}\label{sectionthereandback}

In the commutative case, given a line bundle on 
$\boldp^2$, one can restrict it to $\boldp^2\smallsetminus E$
and conversely, given a  line bundle on  
 $\boldp^2\smallsetminus E$ one can extend it to one on $\boldp^2$. 
The same ideas work for $\qgr S$ for $S\in \ASthree'$. The details are
given in this section, which will in particular prove
that there is a (3-1) correspondence between projective right ideals of
$A(\skl)$ and line bundles in $\qgr \skl$ modulo shifts. We also prove analogous
results for rank one torsion-free modules.
 This is one of a few places in this paper
where the results are distinctly different for different AS regular algebras. In
particular,  for arbitrary $S$ one does not get a finite-to-one correspondence 
(see Proposition~\ref{eg-prop}).

Let $S\in \ASthree'$ and write $A=A(S)$. If $\cM\in \qgr S$, set 
$M=\Gamma^*(\cM)$ and $\cM|_A =M|_A= M[g^{-1}]_0$.\label{M-to-A-defn}
 Note that if $\cM=\pi(N)$
for some other  finitely generated $S$-module  $N$, then $M$ and $N$ only
differ by a finite-dimensional $S$-module and so $M[g^{-1}]\cong N[g^{-1}]$. 
There is a natural filtration  $\Phi^{M}$ on $M|_A$ defined by
$\Phi^{M}_i=M_{3i}g^{-i}$. In the other
direction, let $P=\bigcup \Theta_i$ be a filtered right $A$-module. Then
$\cR(P)=\cR_\Theta(P)=\bigoplus_{i\geq 0}\Theta_ig^i$ is naturally a right
  $S^{(3)}$-module under the identification of Lemma~\ref{firstlemma}(5).
Thus, we may form $\cS_\Theta(P)=\cR(P)\otimes_{S^{(3)}}S$\label{cS-defn} and
$\cV_\Theta(P)=\pi\cS(P)$. Set $R=S[g^{-1}]$ and give each graded piece
$R_n=\sum S_{n+3j}g^{-j}$ the  $A$-module filtration defined by
$\Psi^{R_n}_j=S_{n+3j}g^{-j}$.

\begin{lemma}\label{S-to-A}
Let $S\in \ASthree'$ and $\cM\in \qgr S$ be a torsion-free module. If
$M=\Gamma^*(\cM)$ then:
\begin{enumerate}
\item $\cM \cong \cV(\cM|_A)$, where $\cM|_A$ is given the induced filtration
$\Phi^M$.
\item There is a filtered isomorphism
 $\cM(\ell)|_A \cong (\cM|_A)\otimes _AR_\ell$ for $\ell \in \mathbb Z$.
 \item If $\cN=\cM g$  and $N=\Gamma^*(\cN)$ then $\cN|_A\cong \cM|_A$ and 
 $\Phi^N_i=\Phi^M_{i-1}$ for $i\gg 0$.
\item Conversely, if $P=\bigcup \Theta_i$ is a filtered $A$-module 
then $\cS_\Theta(P)|_A=P$.
 \end{enumerate}
 \end{lemma}
 
 \begin{proof} (1) Since $M$ is torsion-free,
  $\cR(\cM|_A) = \bigoplus_{i\geq 0} M_{3i}g^{-1}g^i \cong M^{(3)}$.
   The natural map 
 $M^{(3)}\otimes_{S^{(3)}}S \to M$ then has bounded kernel and cokernel
 \cite[Theorem~4.4]{Ve} and so $\cV(\cM|_A)\cong \cM$.

(2) As $M|_A$ is torsion-free and $R$ is a projective left $A$-module
(Lemma~\ref{firstlemma}), the natural map 
$(M|_A)\otimes_AR_\ell \to (M|_A)R_\ell $ is an isomorphism and 
we can use products in place of
tensor products. The filtration on $(\cM|_A)R_\ell$   
is then the  natural one defined by
 $\Theta_j=\sum_k\Phi^{M}_{j-k}\Psi^{R_\ell}_{k}$. 
Since  $gS_r=S_rg$,
one obtains $$\Theta_j = \sum_k M_{3j-3k}g^{j-k}S_{3k+\ell}g^{-k}
= \sum_k M_{3j-3k}S_{3k+\ell}g^{-j}
=M_{3j+\ell}g^{-j},$$ for any $j$ such that 
$M$ is generated in degrees $\leq j$.
On the other hand, the filtration on $M(\ell)|_A $
is given by $\Theta_j' = M(\ell)_{3j}g^{-j}  
= M_{\ell + 3j}g^{-j}$. As these vector spaces agree for 
$j\gg 0$, the $A$-modules are also equal.

(3) This is a straightforward computation. 

(4) Note that $\cS_{\Theta}(P)_{3j} = \cR_{\Theta}(P)_{j} = \Theta_jg^j$
and so $\cS(P)|_A = \sum \Theta_jg^jg^{-j}=P$.
\end{proof}

In order to get a way of relating unfiltered 
$A$-modules with modules in $\Qgr S$
we need to relate different filtrations on $A$-modules. 
This is given by the following mild generalization of  \cite[Lemma~10.1]{BW2}.

\begin{lemma}\label{canonical}
  Let $A=\bigcup_{i\geq 0} \Lambda_i$ be any filtered $k$-algebra
 (thus the $\Lambda_i$ are $k$-subspaces of $A$) such that
 both $A$ and $\operatorname{gr} A=\bigoplus
\Lambda_i/\Lambda_{i-1}$ are Ore domains. Let $M=\bigcup \Theta_i$ be a
filtered right ideal of $A$
 such that either   $\dim_k\Theta_i<\infty$ for all $i$ or that
 $\on{gr}_\Theta M=\bigoplus \Theta_i/\Theta_{i-1}$
 is a finitely generated $\on{gr} A$-module.

 Assume, for some $r_0\in{\mathbb N}$, that $(\on{gr}_\Theta M)_{\geq
r_0}$
  is a
torsion-free
 right $\operatorname{gr}(A)$-module.
 Then there exists $j\in\boldz$ such that
 $\Theta_i = \Lambda_{i+j}\cap M$ for all $i\gg 0$.
\end{lemma}

\begin{proof}
We have two filtrations on $M$, the given filtration
$\Theta$ and $\{\Gamma_i=M\cap \Lambda_i\}$.
Write $\Theta^\circ_i=
\Theta_i\smallsetminus \Theta_{i-1}$, and similarly for $\Gamma$ and
$\Lambda$.
Fix
 $m\in \Theta^\circ_r$ and let $n\in \Theta^\circ_{p}$
 for some  $r,p\geq r_0$.  If
 $m\in \Gamma^\circ_s$ and $n\in \Gamma^\circ_q$,
we claim that $r-s=p-q$.

The  assumption that  $(\operatorname{gr}_\Theta M)_{\geq r_0}$
is torsion-free is equivalent to the assertion that $ma\in
\Theta^\circ_{r+t}$
 for all $a\in \Lambda^\circ_t$
with  $t\geq 0$.
As $A$ is Ore, there exist $a,b\in A\smallsetminus\{0\}$, say
$a\in \Lambda^\circ_u$ and $b\in \Lambda^\circ_v$,
such that $ma=nb\not=0$. Thus
$ma\in \Theta^\circ_{r+u}$
while $nb\in \Theta^\circ_{p+v}$.
So $r+u=p+v.$ Since $\operatorname{gr}_\Gamma M \subseteq
\operatorname{gr} A$ is a torsion-free module, the
same argument implies that
$ma\in \Gamma^\circ_{s+u}$
while $nb\in \Gamma^\circ_{q+v}$.
Thus $s+u=q+v.$
Subtracting the two  equations proves the claim that $r-s=p-q$.

As $n$ is arbitrary,  the claim implies that
$\Theta^\circ_p\subseteq \Gamma^\circ_{p-(r-s)}$ for all $p\geq r_0$.
Conversely, either assumption on the $\Theta_j$
implies that  $\Theta_{r_0}\subseteq \Gamma_{p_0}$ for some $p_0>0$.
  Pick   $n' \in
\Gamma^\circ_{q'}$ for some $q'> p_0$ and suppose that
$n'\in \Theta^\circ_{p'}$. Necessarily, $p'\geq r_0$.
Thus,  applying the claim to $n=n'$
shows that
$r-s=p'-q'$ and hence that
$\Gamma^\circ_{q'}\subseteq \Theta^\circ_{q'+(r-s)}$.
Thus, $\Theta^\circ_p= \Gamma^\circ_{p-(r-s)}$
for all $p\gg 0$.
\end{proof}

\begin{notation}\label{canonical-defn}
We will call two filtrations $\Gamma$ and $\Gamma'$ on an $A$-module $M$ {\em
equivalent} if $\Gamma_i=\Gamma'_i$ for all $i\gg 0$.  If $M$ is a finitely
generated torsion-free
rank $1$ right $A$-module, then  we may choose an embedding
 $\psi: M\hookrightarrow A$  as
a right  ideal of $A$ and hence obtain a filtration $\Lambda(M,\psi)$ defined
by $\Lambda(M,\psi)_n = \psi^{-1}(\Lambda_n\cap \psi(M))$.  By the last lemma,
this filtration is unique up to shift and equivalence.  We  call it a {\em
canonical filtration}.
\end{notation}

Let $S=S(E,\cL,\sigma)\in \ASthree'$ with $A=A(S)$ and  $B=B(E,\cL,\sigma)=S/gS$. 
We next want to determine the modules in $\qgr S$ that appear as 
$\cV_\Theta(P)$ for some $P\in\text{Mod-}A$.
Write $\cV=\cV_S$\label{VS-defn}
 for the set of isomorphism classes of rank one, torsion-free
modules $\cM$  in $\qgr S$ for which  $\cM|_E=\cM/\cM g$ is a vector bundle on
$E$. A technical but more convenient description  of $\cV_S$
follows from the next lemma.

\begin{lemma}\label{tf-vble} Set $S=S(E,\cL,\sigma)$ and
let $\M\in \qgr S$ be a torsion-free module.
 Then   $\M|_E$ has
no simple submodules if and only if $\M|_E$ is a vector bundle.
\end{lemma}

\begin{proof} Assume that $\M\in \qgr S$ is a torsion-free module such that 
   $\M|_E$ has no simple submodules. 
We  claim that $\cM$ has a resolution $0\to \cP\to \cQ
\to \cM\to 0$ where $\cP$ and $\cQ$ are sums of shifts of $\theo_S$. 
To see this, note that  $M=\Gamma(\cM)$ is torsion-free and so 
embeds into a direct sum $F$  of shifts of $S$.  By
Lemma~\ref{chi-lemma}(3), $F/M$ has no finite dimensional submodules and so 
\cite[Proposition~2.46(i)]{ATV2} implies that  $\on{hd}(F/M)\leq 2$.
Thus, $\on{hd}(M)\leq 1$ and so $M$ has a graded free resolution of length 
one.  Now apply $\pi$ to this resolution.

By Lemma~\ref{restriction of torsion-free}, the complex $0\to \cP|_E \to
\cQ|_E\to \cM|_E\to 0$ is also exact and so, by Lemma~\ref{shifts},
 $\on{hd}_{E}(\cM|_E)<\infty$.
 Since $E$ is a Gorenstein curve and $\cM|_E$ has no simple subobjects,
 $\cM|_E$ locally has depth at least one. Thus the 
 Auslander-Buchsbaum formula \cite[Theorem~19.9]{Ei}
 proves that $\on{hd}_{E}(\cM|_E)=0$, locally and hence
 globally.  Thus $\cM|_E$ is a vector bundle. The other direction is trivial.
\end{proof}

\begin{lemma}\label{rees} Let $S\in \ASthree'$ and set $A=A(S)$. 
 Suppose that $P$ is a right ideal of $A$ with the induced filtration 
 $\Gamma_i=P\cap \Lambda_i$.    Then $\cV(P)\in \cV$.
\end{lemma}

\begin{proof}   Let $Q=\cR_\Gamma(P)$. By
Lemma~\ref{firstlemma}, $Q$ is a right ideal of $S^{(3)}$  and hence
$Q S$ is a right ideal of $S$.  
Since $\left(Q \otimes_{S^{(3)}}S\right)^{(3)} =
Q= \left(Q S\right)^{(3)},$
 \cite[Theorem~4.4]{Ve} implies that $\cV(P)=\pi(Q\otimes S)=\pi(Q S)$ is
  torsion-free. Now consider $\cV(P)|_E$.
Then
$$ 
Q/Q g\ = \bigoplus \frac{(P\cap \Lambda_i)g^i}{(P\cap \Lambda_{i-1})g^i}
\cong \bigoplus\frac{(P\cap \Lambda_i+\Lambda_{i-1})g^i}{\Lambda_{i-1}g^i} 
 \hookrightarrow S^{(3)}/gS^{(3)} = B^{(3)}.
$$
Once again, \cite{Ve} shows that
 $\cS(P)/\cS(P)g  \cong (Q/Q g)B\hookrightarrow B$ in $\qgr B$.
 Thus, $\cV(P)|_E=\pi(\cS(P)/\cS(P)g)$ has no simple subobjects.
\end{proof}

We let $\sim$ denote the equivalence relation on $\cV$ defined by $\cM\sim
\cM g$. Let $\cP=\cP_A$\label{P-defn} denote the set of isomorphism 
classes of finitely generated,
torsion-free rank one $A$-modules. 

\begin{prop}\label{extension thm for t.f. sheaves} 
Let $S\in \ASthree'$ and write  $A=A(S)$.
Then:
\begin{enumerate}
\item The map $  \cM\mapsto \cM|_A$ induces a surjection
$\psi: (\cV/\hskip -3pt\sim) \to \cP$. 
\item If $E$ is integral, then $ \psi$ is a bijection.
\end{enumerate}
\end{prop}

\begin{proof} (1) If $\cM\in \cV$, then $M=\Gamma^*(\cM)$ is torsion-free
and hence so is $\cM|_A$. By Lemma~\ref{S-to-A},
$\cM|_A \cong (\cM g)|_A$ and so $\psi $ is defined. 
If $P\in\cP$, pick an embedding $P\hookrightarrow A$ and take the induced
filtration $\Theta_i=P\cap \Lambda_i$. 
Then Lemma~\ref{rees} ensures that $\cV(P)\in \cV$ and so $\psi$ is
surjective.

(2) In this case, Corollary~\ref{firstcor} implies that 
$\operatorname{gr}_\Lambda A$ is a domain and hence Lemma~\ref{canonical}
ensures that the filtration $\Theta$ in part (1) is unique up to shift and 
equivalence. Given two equivalent filtrations, say $\Theta$ and $\Theta'$ of
$P$, then $\cS_\Theta(P)$ and 
$\cS_{\Theta'}(P)$ only differ in a finite number of graded pieces 
and so $\cV_\Theta(P)=\cV_{\Theta'}(P)$.
On the other hand, suppose that 
$\Theta'_i=\Theta_{i-1}$, for $i\gg 0$.  Then, in high degree,
$\cR_{\Theta'}(P)=\sum \Theta_{i-1}g^i=\cR_{\Theta}(P)g$.
Thus, $\cV_{\Theta}(P)= \cV_{\Theta'}(P) g$
and so $\cV_{\Theta}(P)\sim  \cV_{\Theta'}(P)$.
Thus, $\psi $ is a bijection.
\end{proof}

We want to refine Proposition~\ref{extension thm for t.f. sheaves}(2) in two ways.
In  Theorem~\ref{three-to-one} we will 
use  an analogue of the first Chern class to 
chose a canonical member of each equivalence class of $\sim$.
Then we will show that under this
equivalence the projective $A(S)$-modules correspond precisely to the 
line bundles in $\qgr S$ (see Corollary~\ref{vble-corr}).

Let  $\cM\in \qgr S$
for $S\in \ASthree$.  The {\it first Chern class}\label{chern-defn}
$c_1(\cM)$ of $\cM$ is defined by the properties that it should be additive on
short exact sequences and satisfy $c_1(\theo_S(j))=j$. As the next lemma shows, 
this  uniquely determines $c_1(\M)$: one simply takes a resolution of
$\M$ by shifts of  $\theo_S$ and applies additivity. 
When $\qgr S\simeq \coh (\boldp^2)$, it is easy to see
that  this definition of $c_1$  coincides with the usual commutative one. 

We first describe some of the basic properties of $c_1(\cM)$.

\begin{lemma}\label{hilbert polynomial'}
Let $\M\in \qgr S$ for $S\in \ASthree$.
\begin{enumerate}
\item There is a unique function $\M\mapsto c_1(\M)$ with the given 
properties.
\item   $c_1(\M(s)) = c_1(\M)+s\cdot \on{rk}(\cM)$
for any $s\in \mathbb Z$.

\item Assume that $S=S(E,\cL,\sigma)\in \ASthree'$ for a smooth curve $E$
and that $\cM$ is torsion-free.
 Then  $c_1(\cM) = \frac{1}{3}\deg
(\cM|_E)$. 
\end{enumerate} 
\end{lemma}

\begin{proof} 
(1) Let $M\in \gr S$. 
For  $q\geq 0$, write 
$ \on{Tor}^{S}_q(M,k)=\bigoplus_{j} k(\ell_{qj})$,
 as a graded  $k$-module and define $c_1'(M)=\sum_{q,j} (-1)^q \ell_{qj}$.
Clearly, $c_1'(-)$ is additive on short exact sequences and is well defined.
Thus it suffices to prove that 
$c_1'(M)=c_1(\pi(M))$. 

If 
$M$ has a graded projective resolution 
$   P^\bullet\to M\to 0$,
with $P^q=\bigoplus_j S(m_{qj})$, then, as in  \cite[(2.8)]{ATV1},
$c_1'(M)=\sum_{q,j} (-1)^q m_{qj}$.
On the other hand, \cite[(2.15)]{ATV1} implies that $c_1'(k)=0$ and hence that 
$c_1'(F)=0$ for any finite dimensional graded $S$-module $F$. 
Thus, Lemma~\ref{chi-lemma}(3) implies that 
$c_1'(M)=c_1'(\Gamma^*(\pi(M))$, for any $M\in \gr S$.
Thus, defining  $c_1(\pi(M))=c_1'(M)$ does indeed give a unique, well-defined
function determined by the properties that it is additive
 on exact sequences and satisfies $c_1(\theo_S(m))=m$ for $m\in \mathbb Z$.

(2) By additivity of $c_1$ and rank, 
it suffices to prove this for $\M=\theo_S(t)$, for which it is obvious.

(3) 
If $P^\bullet \to \cM\to 0$ is a resolution of $\cM$ by shifts of $\theo_S$,
then Lemma~\ref{restriction of torsion-free} implies that 
$P|_{E}^\bullet\to \cM|_E\to 0$ is a resolution of $\cM|_E$ by shifts of
$\theo_E$.  By Lemma~\ref{shifts}, $\deg \theo(n)|_E = 3n$, for any $n$ and so 
the result follows from additivity.
\end{proof}

\begin{thm}\label{three-to-one}
Let $S=S(E,\cL,\sigma)\in \ASthree'$ be such that $E$ is integral. 
Then:
\begin{enumerate}
\item[(1)] There is a bijection between $\cP$ and  
 $\big\{\cM\in \cV : 0\leq c_1(\cM)\leq 2\big\}$.
\item[(2)] There exists a (3-1) correspondence between modules $P\in\cP$ 
and modules $\cM\in \cV$ satisfying
$c_1(\cM)=0$.
\end{enumerate}
\end{thm}
\begin{remark}
(1) By Remark~\ref{as-cases}, this theorem applies to algebras of types 
{\bf A, B, E} and {\bf H}. We give a brief discussion of the other cases
at the end of the section.
 
 (2) This result proves Theorem~\ref{three-to-one-intro} from the
 introduction, modulo  a proof of Theorem~\ref{mainthm1}(2). 
 \end{remark}

\begin{proof}  (1) Suppose that $M\in \gr S$ has a projective resolution 
$P^\bullet\to M\to 0$. 
Since $S(n)g=gS(n)\cong S(n-3)$, the module $Mg$ has a projective resolution 
$P^\bullet g \cong P^\bullet(-3)$. Thus, by Lemma~\ref{hilbert polynomial'}(2),
$c_1(\cM g)=c_1(\cM)-3$
 for any torsion-free, rank one module $\cM\in \qgr S$. 
Now apply Proposition~\ref{extension thm for t.f. sheaves}(2).

(2) Let $\cM\in \cV$ be  such that $0\leq \cM\leq 2$.
By Lemma~\ref{hilbert polynomial'}(2), again,
there exists a unique $r\in \{0,1,2\}$ such that 
 $c_1(\cM(-r))=0$. By Lemma~\ref{S-to-A}(2),
 $\cM(-r)|_A\cong (\cM|_A)\otimes R_{-r}$, giving the (3-1) equivalence. 
 \end{proof}

A natural question raised by  Theorem~\ref{three-to-one}
is whether one can identify  projective right ideals of $A$ in terms of $\qgr
S$-modules. As we show, they correspond precisely to   line bundles  in
$\qgr S$, as defined below.

If $S\in \ASthree$ and $\cM\in \qgr S$, write
\begin{equation}\label{m-vee}
\uHom_{\qgr S}(\cM,\theo) =
\bigoplus_{n\in \mathbb Z}{\Hom}_{\qgr S}(\cM,\theo(n)) 
\in S{\on{-gr}}.
\end{equation}
One can think of $\pi\uHom$  as  sheaf $\Hom$ on $\qgr S$. 
 The right derived functors of $\pi\uHom$
will be denoted by $\pi\underline{\Ext}$ and these are again objects in 
$S\on{-qgr}$. A more complete discussion of these
concepts can, for example, be found in \cite[Section~5.3]{KKO}.
In particular, it is observed there that  
$\pi\uHom_{\qgr S}(\cM,\theo) \cong \pi\uHom_{\gr S}(\Gamma^*(M),S)$,
and so this notation is consistent with the earlier definitions of
$\uHom$ and  $\uExt$.

\begin{defn}\label{vectorbundledefn}(\cite[Definition~5.4]{KKO}) 
Let $S\in \ASthree$. An 
object  $\cF\in\qgr S$  is called a {\it vector bundle} if 
$\pi\underline{\Ext}^i(\cF,\cO)=0$ for all $i>0$. Such an object is necessarily
torsion-free. 
If $\cF$ is a vector bundle of rank one, then $\cF$ is called 
 a {\it line bundle}.\label{line-ble-defn}
 \end{defn}

  This definition of vector bundle is only appropriate for rings of  finite
homological dimension that satisfy the $\chi$ condition. In particular, the
analogous definition for  $\qgr S/gS$ need not correspond to vector bundles in
$\coh(E)$ and so we do not talk about vector bundles in $\qgr S/gS$.

\begin{lemma}\label{vble-reflexive} Let $S=S(E,\cL,\sigma)\in \ASthree'$
 and $\cM\in \qgr S$. 
Then:
\begin{enumerate} 
\item  $\cM$ is a vector bundle if and only if $\Gamma(\cM)$ is a reflexive 
$S$-module. 
\item If $\cM$ is a vector bundle, then $\cM|_E$  is a vector bundle over $E$.
\end{enumerate}
\end{lemma}

\begin{proof} (1) This follows from the Auslander-Gorenstein and
Cohen-Macaulay
conditions (see \cite[Note, p.352]{ATV2}). In more detail, 
let $\cM\in \qgr S$ be a torsion-free module. Since $\cM$ embeds 
in a direct sum of shifts of 
$\theo_S$, Lemma~\ref{chi-lemma}(4) implies that  $M=\Gamma(\cM)$
is finitely generated.
  Set $M^*=\underline{\Hom}(M,S)$.
By \cite[Theorem~4.1]{ATV2}, $ 
\dim_k\underline{\Ext}^j(M,S)<\infty $  for $j>1$ and  $ \GKdim\,
\underline{\Ext}^1(M,S)\leq 1$.  By \cite[(4.4)]{ATV2} 
the canonical map 
$M\to M^{**} $ has cokernel 
$Q\subseteq \underline{\Ext}^2(\underline{\Ext}^1(M,S),S)$. 
By Lemma~\ref{chi-lemma}(3), $M$ has no finite dimensional extensions, 
so $M$ is reflexive if and only if $\dim_kQ<\infty$.
By \cite[Theorem~4.1(iii)]{ATV2} this happens  
 if and only if  $\dim_k\underline{\Ext}^1(M,S)<\infty$.

By \cite[Equation~14, p. 406]{KKO},
 $\pi\underline{\Ext}^j(\cM,\cO)= 
\pi(\underline{\Ext}^j(\Gamma^*(\cM),S) )$. 
Thus, $\cM$ is a vector bundle if and only if  
$\dim_k\underline{\Ext}^j(\Gamma^*(\cM),S)<\infty$ for all $j>0$. 
By the last paragraph this is equivalent to $M$
being reflexive. 

(2) Let $M= \Gamma(\cM)$ and $N=M^*$. 
Since $M$ is reflexive and hence torsion-free, it is easy to see that 
$M/Mg \hookrightarrow \Hom_S(N,S/gS)
=\Hom_{S/gS}(N/gN,S/gS).$
Therefore, $\cM|_E = \pi(M/Mg) \subset\pi(\Hom_{S/gS}(N/gN,S/gS))$ 
which certainly has no simple submodules. Now apply Lemma~\ref{tf-vble}.
\end{proof}

\begin{corollary}\label{vble-corr}
Let $S=S(E,\cL,\sigma)\in \ASthree'$ with $A=A(S)$  and let 
$\cM\in \qgr S$ be
 such that $\cM|_E$ is a vector bundle.
 Then:
\begin{enumerate}
\item $\cM$ is a vector bundle if and only if $\cM|_A$ is a
projective $A$-module. 
\item If $|\sigma|=\infty$ as an automorphism of the scheme $E$,
 then $\cM$   is automatically a vector bundle.
\end{enumerate}
\end{corollary}

\begin{proof} (1) Set  $M=\Gamma(\cM)$, with double dual $M^{**}$.
By \cite[Corollary~4.2(iv)]{ATV2}, $\GKdim(M^{**}/M)\leq 1$.
 Suppose, first, that there 
is a proper essential extension  $0\to M\to N\to H\to 0$ with
$\GKdim(H)\leq 1$  and 
$Hg=0$. Since $M$ has no extensions by finite dimensional modules, 
$\GKdim H = 1$. 
Consider the  exact sequence 
$$ \Tor^1_S(H,S/Sg)\
\buildrel{\theta}\over{\longrightarrow}\ M/Mg
\longrightarrow N/Ng\longrightarrow H/Hg\longrightarrow 0.$$
Now, $\GKdim \Tor^1_S(H,S/Sg)\leq \GKdim H \leq 1 $ by
\cite[Proposition 2.29(iv)]{ATV2}.
 By Lemma~\ref{tf-vble}, this forces $\on{Im}(\theta)=0$ and hence 
$Mg=M\cap Ng$. On the other hand, 
$Ng\subseteq M$ since $Hg=0$. Thus, $Mg=Ng$. Since $M$ and $N$ are torsion-free,
this implies that $M=N$. Thus, no such extension exists.

By Lemma~\ref{vble-reflexive},
 $\cM$ is  a vector bundle if and only if 
$M$ is   reflexive. 
By the last paragraph this happens   if and only if 
$M[g^{-1}]$ is  reflexive. 
By the equivalence of categories $\module A \simeq \gr S[g^{-1}]$
\cite [Theorem~A.I.3.4]{NV}, this happens if and only if 
$\cM|_A=M[g^{-1}]_0$ is   reflexive. Finally, Lemma~\ref{firstlemma}
implies that 
$\on{gldim}(A)\leq 2$ and so reflexive $A$-modules 
are the same as projective $A$-modules. 

(2)  This is immediate from part (1) combined with Lemma~\ref{firstlemma}(3).
\end{proof}

\begin{remark}\label{vble-remark} Despite the fact that it holds for all values
of $|\sigma|$, there is actually a striking dichotomy in 
Theorem~\ref{three-to-one}: If
$|\sigma|=\infty$ then   Corollary~\ref{vble-corr}(2) applies. In contrast,
when $|\sigma|<\infty$, $A$ has global (and Krull) dimension two. It follows
that $A$ has many non-projective right ideals and hence that $\cV$ contains
many modules $\cM$ that are not line bundles. \end{remark}

The reason for demanding that $E$ be integral in 
 Proposition~\ref{extension thm for t.f. sheaves}(2)
and Theorem~\ref{three-to-one} is so that 
Lemma~\ref{canonical} can be applied 
to restrict the number of possible filtrations on an
$A$-module $P$. If $E$ is not integral, then that lemma can fail and so 
  the map $\psi$  in
 Proposition~\ref{extension thm for t.f. sheaves}(1)
need not  be a bijection.

 An example where $\psi$ has infinite  
fibres is given by the ring
\begin{equation}\label{eg-eg}
T=k\{x,y,z\}/(yx-pxy,xz-zx,zy-yz),
\end{equation}
 where $p\in k$ is
transcendental over the prime subfield. 
In the notation of Remark~\ref{as-cases}, 
$T$ is an   AS-regular ring of Type $\mathbf{S}_1$ with 
 $g=xyz$.   Moreover,  $T$ is an Ore extension
$T=R[x,\tau]$, where $\tau$ is an automorphism 
of the polynomial ring $R=k[y,z]$;
 thus multiplication  is defined by $xr=\tau(r)x$ for $r\in R$.   
We   regard elements of $T$ as polynomials in $x$ with
coefficients in $R$ and  write  $\deg_x$ for the corresponding degree
function.

A method for constructing non-free projective modules over Ore extensions is
given in \cite[Theorem~1.2]{St2} and we use a similar technique to build
reflexive  $T$-modules.

\begin{lemma}\label{eg-comp}
Define $T$ by \eqref{eg-eg} and set 
 $M =\big\{t\in T : (z+y)t\in (x+z)T\big\}$.
 Then:
 \begin{enumerate}
 \item $M$ contains 
 $\alpha=x^2+x(1+p)z+pz^2$ and $\beta= x(z+y)+z(z+py)$.
 
 \item $M$ contains no polynomial of the form
  $\gamma=xz^i+r$ for $r\in R$ and $i\geq 0$.
 
 \item  $M[g^{-1}]_0$ is a non-cyclic projective right ideal of $A(T)$.
 \end{enumerate}
 \end{lemma}
 
 \begin{proof}
 (1) Use the two identities
$(y+z)\alpha = (x+z)\left(x(p^2y+z)+pz(y+z)\right)
$ and 
$(y+z)\beta = (x+z)(y+z)(py+z).$

  (2) Suppose that $  xz^i+r\in M$,
for some $i\geq 0$ and $r\in R$.  Then 
$$(z+y)(xz^i+r) = x(z^{i+1}+pyz^i) +
(z+y)r =(x+z)(z^{i+1}+pyz^i) + w,$$ for $w=(z+y)r-(z^{i+1}+pyz^i).$ 
Since $\deg_xw=0$, the definition of $M$  forces
$w=0$. But the equation $w=0$  is  impossible in the polynomial ring $k[y,z]$.

 (3) This follows easily from \cite[Theorem~1.2]{St2}.
\end{proof}

Extend $\tau$ to an automorphism of $T$ by defining 
$\tau(t) = xtx^{-1}$, for $t\in T$.

\begin{prop}\label{eg-prop} Let $M$ be
the $T$-module defined in Lemma~\ref{eg-comp}. Then:
\begin{enumerate}
\item  $\pi(M)$ is a line bundle
such that $\pi(M)\not\cong\pi(M^{\tau^i}[j])$, for 
$(i,j)\not=(0,0)$.

\item However, $M[g^{-1}]_0 \cong N[g^{-1}]_0$, where $N=M^{\tau^i}[j]$,
for any $i,j\in \mathbb{Z}$. 

\item Thus, the map $\psi$ defined in 
Proposition~\ref{extension thm for t.f. sheaves}(1)
has some infinite fibres. 
\end{enumerate}
\end{prop}

\begin{proof} (1) Since $M$ is reflexive,
Lemma~\ref{vble-reflexive} implies that $\pi(M)$ is a line bundle. 

 Suppose there exists an  isomorphism 
$\theta: \pi(M)\to \pi(N)$, where $N=\pi(M^{\tau^i}[j])$, for some $i,j$.
Then 
$\theta$ induces a homomorphism $\theta$ from 
$M$ to the injective hull of $N$ such that
$\theta(M_{\geq n})= N_{\geq n}$, for all $n\gg 0$.
However, by construction, $M$ and  $N$ are reflexive and so 
have no essential extensions by finite dimensional $T$-modules.
  Thus, $\theta$   induces an isomorphism 
from $M$ to $N$. Since $\tau$ is a graded isomorphism 
of $T$, this implies that $M$ and $N$  have the same Hilbert series. As $M$
 and $M^{\tau^i}[j]$ have different Hilbert series for $j\not=0$, this forces 
$j=0$. 

We may write
 $f_1M=f_2N$ for some 
$f_i\in T_d$ and  some $d$. By \cite{RSS}, elements in $T$
that are monic when 
regarded as polynomials in $x$ form an Ore set $\mathcal{C}$.
Localizing at $\mathcal{C}$ and using Lemma~\ref{eg-comp}(1) 
gives the identity
 $M_{\mathcal{C}} =T_{\mathcal{C}} =
N_{\mathcal{C}}$. Thus $f_1T _{\mathcal{C}} =f_2T _{\mathcal{C}} $. 
The only units
in $T _{\mathcal{C}} $ are of the form $\lambda_1g_1^{-1}g_2$, where $\lambda\in
k^*$ and $g_i\in \mathcal{C}$. Thus, we may assume that the $f_i\in
\mathcal{C}$.

Consider elements of least degree in $x$ in $f_1M=f_2N$. 
By Lemma~\ref{eg-comp}, this degree is $(d+1)$ and 
$f_2N\ni f_2\beta^{\tau^i} = x^{d+1}(z+p^{-i}y) +$ lower order terms. 
But, $f_1M\ni f_1\beta$, with leading term $x^{d+1}(z+y) $. 
Thus $ f_2\beta^{\tau^i}- p^{-i} f_1\beta \in f_1M$ has
 leading term $x^{d+1}z(1-p^{-i})$. Since 
$p^i\not=1$, this implies  that 
 $M$ contains an element of degree $1$ with leading term $xz$.
 This contradicts
 Lemma~\ref{eg-comp}(2).

(2) Both $w=x$ and $w=z$,
are normal elements in $T$. Thus, for any 
right ideal $M'$ of $T$ there exists a short exact sequence: 
$0\to M'w\to M'\to M'/M'w\to 0.$
The map $\psi: P\mapsto P[g^{-1}]_0$ is exact and 
so $\psi(M')=\psi(M'w)$. If $w=z$, then $z$ is central and so
$M'z=zM'\cong M'[-1]$. By induction, $\psi(M)=\psi(M[i])$ for all $i\in \mathbb Z$. 
On the other hand, $M'x=x\tau^{-1}(M')=\tau^{-1}(M')[-1]$.
Thus a second induction proves (2).

(3) This is immediate from (1) and (2).
\end{proof}

One  reason why Theorem~\ref{three-to-one} fails for $T$  (and other examples
where $E$ is not integral) is that $g$ is a
product of normal elements. This introduces extra units into $A(S)$ and ``too
many'' isomorphisms between projective modules. However, if one is willing to change the
ring $A(S)$ then  one can still get an analogue of
Theorem~\ref{three-to-one}.   We end this section by  outlining the result but
since it depends upon a case by case analysis  the details will be left to the
 reader.

Let  $S\in \ASthree'$ and suppose first that
$E$ has an integral component $X$ of
degree $r\leq 2$ that is fixed by $\sigma$. 
Using the arguments of \cite[Section~4]{AS}
one can show that $B'=B(X,\cL|X,\sigma|_X)=S/hS$ for some normal
element $h\in S_r$ that divides $g$.
This happens, for example, in Types $\mathbf
S_1$ (where $r=1$)  and $\mathbf S'_1$ (where $r$ can be $1$ or $2$). 
In Type 
$\mathbf S_2$ it only happens when $X$ is the line 
fixed by $\sigma$.  Now set
$A'=S[h^{-1}]_0$; thus $A(S)$ is a localization of $A'$.
 Since $X$ is integral, $B'$ is a domain
and the argument of Corollary~\ref{firstcor}(2) shows that
$\gr_{\Lambda'}A'$ is also a domain, where $\Lambda'_i = S_{ir}h^{-i}$.
The proof of Theorem~\ref{three-to-one} then goes through to give an
$r$-to-$1$
 correspondence between isomorphism classes of rank one torsion-free
$A'$-modules and elements of $\cV$ with $c_1=0$. When $S=U$ is the homogenized
Weyl algebra, $g=h^3$ and so $A'=A(S)$ and one recovers the result from
\cite{BW1}.

There are two further cases to be considered.  First  consider $S$ of type
$\mathbf S_2$ and  let $Y$ be the union of the lines interchanged by $\sigma$.
Then $B''=B(Y,\cL|_Y,\sigma|_Y)=S/h'S$  where $\deg h'=2$. In this case $B''$
is a prime ring that is not a domain, so the analysis of the last paragraph
will  fail. However, if one replaces $S$ by the Veronese ring $S^{(2)}$
then $h'$ will become a product of two  normal elements (this is because
$\sigma^2$ fixes the components of $Y$) and the above analysis can be
pushed through. The final case is when $E$ is a product of $3$ lines cyclicly
permuted by $\sigma$ (this is a degenerate   Type {\bf A} example,
 when $b=0$ in \eqref{skly-def2}).
 Here one has to use the $3$-Veronese ring~$S^{(3)}$.

\section{Cohomology and Base Change in Noncommutative Geometry}
\label{section-base}

In this section we describe two cohomological results that will be needed in
the sequel.  Both are minor variants of results from the literature. 
The first describes the general  
``Cohomology and Base Change'' machinery from \cite{EGA} in  the
appropriate generality to apply to noncommutative projective schemes.
These theorems give detailed information 
about the variation of cohomology
 in a family parametrized by a
noetherian base scheme $U$. The second result describes the noncommutative 
\v{C}ech cohomology of \cite{VW2} in a form appropriate to our needs. 

The results of this section need rather strong hypotheses, but these are
probably necessary. However, as will be shown in Lemma~\ref{beilinson
hypotheses} they do hold for $S_A=S\otimes_kA$, where $A$ is a commutative
noetherian $k$-algebra and $S\in \ASthree$.

\subsection{Cohomology and Base Change}\label{cbc-subsec}

The results of this subsection will require  the following hypotheses
on our algebras:

\begin{hypotheses}\label{CBC-hypotheses}
Fix a commutative noetherian ring $A$ and a finitely generated, 
connected graded
$A$-algebra $S$.  Assume  that:
\begin{enumerate}
\item $S$ is \emph{strongly  noetherian}\label{strong-defn}
 in the sense
that $S_B=S\otimes_A B$ is  noetherian
for every commutative noetherian $A$-algebra $B$.
\item $S$  satisfies
 $\chi$ and  $\Qgr S$ has  finite cohomological dimension.
 \end{enumerate}
 \end{hypotheses}

If $B$ is a commutative $A$-algebra and $S$ is a cg $A$-algebra,
we regard $S_B=S\otimes_AB$\label{S-sub-B}
 as a cg $B$-algebra. By \cite[Proposition~B8.1]{AZ2},
$\Qgr S_B$ is equivalent to the category of $B$-objects in $\Qgr S$. 
 If $\M\in \qgr S_B$, with $S_B$ noetherian, then 
\cite[Lemma~C6.6]{AZ2} implies that there is a canonical identification 
of cohomology
groups $\HH^i(\Qgr S_B, \M) = \HH^i(\Qgr S,\M)$ via the natural functor 
$\Qgr S_B\rightarrow \Qgr S$ and so we may write this group as 
$\HH^i(\M)$ without confusion.
 Generalizing  earlier notation, set $\theo_{S_B} =\pi(S_B)\in \qgr S_B$.
If the context is clear, we will typically write $\theo$ for 
  $\theo_{S_B}$.
  An object $\M\in \Qgr S$ is {\em $A$-flat}\label{flat-defn} if the functor
 $
-\otimes_A \M : \Module A\longrightarrow \Qgr S
$
is exact. When $S_A$ is strongly noetherian,
\cite[Lemma~E5.3]{AZ2} implies that $\M\in \qgr S$ is $A$-flat if and only if
$\HH^0(\qgr S, \M(n))$ is flat for $n\gg 0$.

\begin{thm}[Theorem on Formal Functions]\label{FF}
Assume that $A$ and $S$ satisfy Hypotheses~\ref{CBC-hypotheses} and let
 $\cF,\cG\in \qgr S$.  
Then for every $i$ and every ideal ${\mathfrak m}$ of $A$, 
the canonical homomorphism
\bd
\Ext^i_{\qgr S}(\cF,\cG)\otimes \widehat{A} \longrightarrow 
\underset{\longleftarrow}
{\lim}\; \Ext^i_{\qgr S}(\cF,\cG\otimes_A A/{\mathfrak m}^k)
\ed
 is an isomorphism.
\end{thm}

\begin{proof} 
 \cite[Proposition C6.10(i)]{AZ2} implies that 
$S$ satisfies the strong $\chi$ condition
of \cite[(C6.8)]{AZ2}, while
\cite[Proposition C6.9]{AZ2}  implies that $\Qgr S$ is
Ext-finite. Thus the hypotheses of \cite[Theorem D5.1]{AZ2}
(with $k=A$ and $R=\widehat{A}$) are satisfied.
\end{proof}

We remark that, for a commutative ring  $R$, a ``point $y$ of $\spec R$''
means any point (not necessarily a closed point).

\begin{thm}[Cohomology and Base Change]\label{CBC}
Assume that $A$ and $S$ satisfy Hypotheses~\ref{CBC-hypotheses} and let $y\in
\spec A$.
Pick $\cF,\cG\in\qgr S$ such that $\cG$ is $A$-flat.  
Then:
\begin{itemize}
\item[(1)] If the natural map 
\bd
\phi_i: \Ext^i(\cF,\cG)\otimes_A k(y)\longrightarrow
 \Ext^i(\cF,\cG\otimes_A k(y))
\ed
is surjective, then it is an isomorphism.

\item[(2)] If  $\phi_{i-1}$ and $\phi_i$ are  surjective,
then $\Ext^i(\cF,\cG)$ is a vector bundle in a neighborhood of $y$ in $\spec A$.
\item[(3)] If $\Ext^{i+1}(\cF,\cG\otimes k(y))
=\Ext^{i-1}(\cF,\cG\otimes k(y))=0$,
then $\Ext^i(\cF,\cG)$ is a vector bundle in a neighborhood of $y$ in $\spec A$
and $\phi_i$ is an isomorphism.
\item[(4)] If $\cF$ and $\theo_S$ are also $A$-flat, then 
\bd
\Ext^i_{\qgr S}(\cF,\cG\otimes k(y)) = 
\Ext^i_{\qgr S\otimes k(y)}(\cF\otimes k(y),\cG\otimes k(y)).
\ed
\end{itemize}
\end{thm}

\begin{proof}
Since $\cG$ is $A$-flat, the collection of functors  $M\mapsto T^i(M) =
\Ext^i(\cF,\cG\otimes_A M)$ forms a cohomological $\delta$-functor 
in  the sense of
\cite[Section~III.7]{EGA} or \cite[Section~III.1]{H}. 
 By our assumptions, $T^i(M)\in \text{mod-}A$ for every $M\in \text{mod-}A$.
   Moreover, $T^i$ commutes with
colimits and, by   Theorem~\ref{FF}, the canonical homomorphism
$T^i(M)\,\widehat{\hspace{.4em}}\longrightarrow 
\underset{\longleftarrow}{\lim}\; T^i(M\otimes A/I^k)$ is an isomorphism
for every $M\in \text{mod-}A$ and ideal
$I\subset A$. The
proof of \cite[Proposition~12.10]{H} now shows that, if $\phi_i$ is
surjective, then $T^i$ is a right exact functor of $A_P$-modules, where
$P$ is the ideal of $A$ corresponding to $y$. Thus $\phi_i$ is an
isomorphism by \cite[Proposition~III.7.3.1]{EGA}, proving (1).

If $\phi_{i-1}$ is also surjective
then \cite[Proposition~III.7.5.4]{EGA} implies that $T^i(A_P)$ is a free
$A_P$-module. Since $T^i$ commutes with colimits we have $T^i(A)_P =T^i(A_P)$,
proving (2).  Finally, part (3) follows from \cite[Corollaire~III.7.5.5]{EGA} 
and part (4)
from \cite[Proposition~C3.4(i)]{AZ2}.    \end{proof}

\begin{remark}\label{well-defn-h}
One consequence of Theorem~\ref{CBC}(4) is that, under the hypotheses of that
 result,  $\HH^i(\Qgr S, \cG\otimes k(y))=\HH^i(\Qgr S\otimes_Ak(p), 
 \cG\otimes k(y))$,  and so  we can write this group as $\HH^i(\cG\otimes k(y))$
 without ambiguity.
\end{remark}

\subsection{Schematic Algebras and the \v{C}ech Complex}\label{schematic-defn}
Let $A$ be a commutative ring and $S$ a noetherian connected
graded $A$-algebra.
Following \cite{VW2}, 
we say $S$ is {\em $A$-schematic}\label{schematic-defn2}
 if there is a finite set
$C_1, \dots, C_N$ of two-sided homogeneous Ore sets 
satisfying
\begin{enumerate}
\item $C_i\cap S_+ \neq\emptyset$ for $i=1,\dots, N$.
\item For all 
$\displaystyle (c_1,\dots,c_N) \in \prod_{i=1}^N C_i$
there exists $m\in{\mathbf N}$ with 
$\ds (S_+)^m\subseteq \sum_{i=1}^N c_i\cdot S$.
\end{enumerate}

Given an $A$-schematic algebra $S$, fix Ore sets
$C_1,\dots, C_N$ as in the definition and set
$I=\{1,2,\dots, N\}$.  Given $w=(i_0,\dots,i_{p-1})\in I^p$,
set $Q_w = S_{C_{i_0}} \otimes_S \dots \otimes_S S_{C_{i_{p-1}}}$,
the tensor product of localizations.
The reader should be warned that the various localizations will not commute
in general and so $Q_w$ need not be a ring and it 
does depend on the ordering of $w$. Of course $Q_w$  is an $S$-bimodule. 
  The {\em noncommutative \v{C}ech complex} $\boldc^{\bullet}$\label{cech-defn}
is then given by setting
\bd
\boldc^p = \prod_{w\in I^{p+1}} Q_w, \;\; p\geq 0,
\ed
and taking as differentials $\boldc^p\rightarrow \boldc^{p+1}$
the usual alternating sum of maps just as in the commutative
setting.  

Define a functor from   $\Ggr S$ to complexes of
graded $S$-modules by $M\mapsto M\otimes_S \boldc^{\bullet}$.
Note that, since $C_j\cap S_+ \neq \emptyset$ for every $j$,
if $T$ is a torsion module then $T\otimes_S\boldc^{\bullet}=0$.
Hence $-\otimes_S\boldc^{\bullet}$ descends to a functor from
$\Qgr S$ to the category of complexes of graded $S$-modules.  
We write 
${\check{\mathbf H}}^i(\M) =  \HH^i(\M\otimes_S\boldc^{\bullet})$
for $\M\in \Qgr S$.
This is a graded $S$-module and we define
${ \check{H}}^i(\M)
=\HH^i(\M\otimes_S\boldc^{\bullet})_0$; 
equivalently ${ \check{H}}^i(\M)$ is the cohomology of 
$\M\otimes \mathbf{C}^\bullet_0$.

Given $\M\in\Qgr S$, write
$
{\mathbf H}^i(\M) = \bigoplus_{n\in{\mathbf Z}} \HH^i({\Qgr S}, \M(n)).
$
This is naturally a graded $S$-module with 
${\mathbf H}^0(\M)$ being the module $\Gamma(\cM)$ defined after
\eqref{gamma-fn}.

\begin{thm}[\cite{VW2}]\label{cech cohomology}
Suppose that $S$ is a noetherian $A$-schematic algebra.
Then there is an isomorphism of $\delta$-functors
${\mathbf H}^i\rightarrow \HH^i(-\otimes_S\boldc^{\bullet})$.

Moreover, this isomorphism is graded in the sense that 
${\rm H}^i(\Qgr S_A, \mathcal M(j)) \cong \check{H}^i(\mathcal M(j))$
for all $ i$ and $j.$
\end{thm}

\begin{proof}
The first assertion of the 
theorem is proved in \cite[Theorem~4]{VW2} when $A$ is a field, 
but the proof works equally well in the more general context.

   The final claim of the theorem is only implicit in \cite{VW2} but
since the proof is short we will give it. It suffices to take $j=0$.
Recall that the $Q_w$ are flat $S$-modules. Thus, for any short exact sequence
$0\to M\to N \to P\to 0$ in $\Ggr S$ one obtains an exact sequence  $$0\to
{\mathbf C}^p(M)_0 \to {\mathbf C}^p(N)_0 \to {\mathbf C}^p(P)_0 \to 0.$$ It
follows that the homology functors $\check{\mathbf H}^i_0=\check{H}^i$  form a
cohomological $\delta$-functor.
If $E\in\Ggr S$ is injective, then \cite[Theorem~3]{VW2} implies that
$\check{\mathbf H}^i(E)=0=\check{\rm H}^i(E)$ for all $i>0$.
Thus the $\delta$-functor is effaceable. Finally, as is observed in 
\cite[p.79]{VW2}, $\check{\mathbf H}^0(\mathcal M)  =
\Gamma^*(\mathcal M) $ and so $\check{\rm H}^0(\mathcal M)= \Hom_{\Qgr
S}(\pi(S),\mathcal M) = {\rm H}^0(\mathcal M)$ for any $\mathcal M\in \Qgr S$. 
Therefore  the hypotheses of \cite[Corollary III.1.4]{H} are satisfied and  ${\rm
H}^i_{\Qgr S}(\mathcal M) \cong \check{H}^i(\mathcal M)$ for all $i.$
\end{proof}

\begin{remark}\label{rog-egs}
 We have already noted that the results in
this section apply to $S_A$, where $A$ is a commutative noetherian $k$-algebra
and $S\in \ASthree$. However, they do not apply to all cg noetherian algebras.
Specifically the noetherian cg algebras studied in \cite{rog} 
are not strongly noetherian, do not satisfy $\chi$ and are not schematic
(see \cite[Theorems~1.2, 1.3 and Proposition~11.8]{rog}, respectively). 

One can use the noncommutative \v{C}ech complex to prove
analogues of Theorem~\ref{CBC} under formally weaker hypotheses than
those used here, but we know of no applications of those generalizations. 
\end{remark}

\section{Monads and the Beilinson Spectral Sequence}\label{monads and beilinson}

A classical approach to the study of moduli of sheaves on ${\boldp}^2$ uses the
Beilinson spectral sequence to describe sheaves as the homology of monads. 
This can then be used to reduce the moduli problem to a question in linear
algebra. In this section we show that  an analogue of the  Beilinson spectral
sequence also works for vector bundles on a noncommutative $\boldp^2$; 
in other words for vector bundles in $\qgr S$ when $S\in \ASthree$.
 As will be seen in the  later sections, this will enable
us to  construct a projective moduli space as a GIT quotient of a subvariety of
a product of Grassmannians. 
  
   The
following properties of $S$ will be needed.
 We use the standard definition of a Koszul algebra as given, for example, in 
\cite[Definition~4.4]{KKO}.
 \begin{lemma} \label{beilinson hypotheses}
If $S\in \ASthree$,  then $S$ is strongly noetherian, schematic and Koszul.
\end{lemma}

\begin{proof} That   $S$ is Koszul follows from 
\cite[Theorem~5.11]{Sm}. By \cite[Theorem~2]{ATV1}, either $S$ or a factor ring
$S/gS$ is a twisted homogeneous coordinate  ring. In either case,  by
\cite[Propositions~4.9(1) and  4.13]{ASZ}, $S$ is strongly noetherian.

The fact that the Sklyanin algebra is schematic is proved in detail in 
\cite[Theorem~4]{VW3}. The idea of the proof is that one can reduce to the case
when the base field $k$ is a finite field. The algebra is then a PI algebra,
for which the result is easy. As in the proof of \cite[Theorem~7.1]{ATV2}, this
argument works in general. \end{proof}

 \begin{defn}\label{monad-defn}  Let $S\in \ASthree$.
  A  complex $\bK$  in $\qgr S_R$ 
 is called a \emph{monad} if it has the form
\begin{equation}\label{monad0}
\bK: \;\;\;
0 \rightarrow V_{-1}\otimes_R \theo(-1)
 \xrightarrow{A_{\bK}} V_0 \otimes_R 
\theo \xrightarrow{B_{\bK}} V_1 \otimes_R \theo(1)
\rightarrow 0
\end{equation}
with the following properties:
\begin{enumerate}
\item[(1)] The $V_i$ are finitely generated projective $R$-modules.

\item[(2)] For every 
 $p\in\spec R$, the complex $\bK\otimes_R k(p)$ 
is exact at $V_{-1}\otimes k(p)\otimes\theo(-1)$ and 
$V_1\otimes k(p)\otimes\theo(1)$.  
\end{enumerate}
Let $\on{Monad}(S_R)$ denote the category of monads, with 
 morphisms being homomorphisms of complexes.
 A monad  $\bK$ is called \emph{torsion-free} if
the cohomology of the complex $\bK\otimes_R k(p)$ at 
$V_0\otimes k(p)\otimes \theo$ is  torsion-free for all 
$p\in\spec R$.

The terms in the monad $\bK$ are indexed so that $V_i\otimes\theo(i)$ lies in 
cohomological degree $i$. Thus one has a functor
$\HH^0: \on{Monad}(S_R) \rightarrow \qgr S_R$ given 
by $\bK\mapsto \HH^0(\bK)$.
\end{defn}   

It is convenient through much of the paper to have a general definition
for a property $P$ of objects of $\qgr S$ to apply to a family:
\begin{defn}\label{families def} Assume that $\bK$ is an object or a complex
in $\qgr S_R$.
Let $P$ be a property of objects or complexes of $\qgr S$ (for example, $P$ 
could be ``torsion-free,'' or ``(semi)stable'' (page \pageref{stability 
definition for modules})).  We say that  $\bK$ 
{\em has $P$} (or {\em is a family of $P$-objects}) if 
 $\bK\otimes k(p)$ has property $P$ for every
point $p\in\spec R$.  Similarly, we say that $\bK$ is {\em geometrically $P$} 
if $\bK\otimes F$ has 
$P$ for every
geometric point $\spec F\rightarrow \spec R$.
\end{defn}

Although most of the earlier results of this paper were concerned with 
rank one torsion-free modules,
the arguments of this section work for any module in $\qgr S_R$
that satisfies the following hypotheses. As will be seen later, these 
conditions are  satisfied by suitably normalized, 
stable vector bundles and so the results of
this section will also form part of the proof of 
Theorem~\ref{firstthm}
from the introduction.

\begin{VC}\label{thesheaf}  Let $R$ be a commutative ring and pick
$S\in \ASthree$. 
Write $(\qgr S_R)_{\on{VC}}$ for the full
 subcategory of $\qgr S_R$
consisting of $R$-flat objects $\cM$ satisfying:
$$\HH^0(\cM(i)\otimes k(p)) = \HH^2(\cM(i)\otimes k(p)) 
 = 0\qquad{\text{for}}\ i=-1,-2\quad \text{and}\quad p\in\spec R.$$
\end{VC}

\begin{remark}\label{implied by VC} Combined with Theorem~\ref{CBC}(3), 
this vanishing condition implies
 that $\HH^1(\cM(i))$ is a projective $R$-module for $i=-1,-2$.
\end{remark}

The next theorem gives our version of the Beilinson spectral sequence.
Later in the section we will extend this result to 
produce an equivalence of categories between 
$\on{Monad}(S_R)$ 
and $(\qgr S_R)_{\on{VC}}$.

\begin{thm}\label{beilinson5'}  Assume that $S\in \ASthree$
and let $R$ be a commutative noetherian $k$-algebra.
Fix $\cM\in (\qgr S_R)_{\on{VC}}$. 
 Then ${\cM}$  is the cohomology of the monad:
$$\bK(\cM) :\ 
0\to V_{-1}\otimes_R \theo(-1) 
\to   V  \otimes_R \theo 
\to V_1\otimes_R \theo(1)\to 0,
$$
where $V_{-1}= \HH^1(\qgr S_R, \cM(-2))$,\ \ $V_1=\HH^1(\qgr S_R, \cM(-1))$
and each $V_j$ is a finitely generated projective $R$-module.
\end{thm}
The module $V$ is defined by \eqref{v0-defn}, but the definition is not
particularly helpful.

\begin{proof} The first part of the proof follows the argument of
\cite[Theorem~6.6]{KKO}.
By Lemma~\ref{beilinson hypotheses}, $S$ is a Koszul algebra.
In particular,  $S=T(S_1)/(\mathcal R)$ is a quadratic algebra
 with Koszul 
dual    $S^! = T(S_1^*)/(\mathcal R^\perp).$  
By \cite[Theorem~5.9]{Sm} the  
{\it augmented left Koszul resolution} for $S$ is 
$$0\to S(-3)\otimes(S^!_3)^* \to S(-2)\otimes(S^!_2)^*  \to 
S(-1)\otimes(S^!_1)^*  \to 
S\otimes(S^!_0)^* \to k\to 0.$$
Define $\Omega^1$\label{omega-defn} to be the cohomology of this Koszul complex 
truncated at  $S(-1)\otimes(S^!_1)^* $; equivalently, 
$\Omega^1$ is defined by the exact sequence 
\begin{equation}\label{omega1}
0\longrightarrow \Omega^1  \longrightarrow S(-1)\otimes S_1
\longrightarrow
S\longrightarrow  k \longrightarrow  0.
\end{equation}
Let $\widetilde{\Omega}^1$
 denote the image of
$\Omega^1$ in $S$-qgr. Similarly, given an $S$-bimodule $M$,  regard $M$ 
as a left $S$-module  and write $\widetilde{M}$
for the image of $M$ in $S$-qgr.

Define the {\it diagonal bigraded algebra of $S$}\label{delta-defn}
to be  $\Delta=\bigoplus_{i,j} \Delta_{ij}$, where $\Delta_{ij}=S_{i+j}$. 
Then, \cite[Equation~11, p.402]{KKO} combined with the above observations 
 proves that the following  complex of bigraded $S$-bimodules is exact:
\begin{equation}\label{diag-resn1} 
0\to S(-1)\boxtimes S(-2) \to \Omega^1(1)\boxtimes S(-1)\to S\otimes S
 \to \Delta \to 0. \end{equation}
Here, $\boxtimes$ stands for external tensor product. 
To save space we write  \eqref{diag-resn1} as
$$0\to K^{-2}\to K^{-1}\to K^0\to{\Delta}\to 0.$$

Let $\cM\in (\qgr S_R)_{\on{VC}}$. In the
commutative setting, the Beilinson spectral sequence has $E_1^{pq}$ term
$\HH^q\left({\boldp}^2, \cM\otimes  \widetilde{\Omega}^{-p}(-p)\right)\otimes_k
\cO(p)$ and this converges to $\cM$. Unfortunately, this does not make sense in our 
situation: $\Omega^1$ is only a left $S$-module and so  $\cM\otimes 
\widetilde{\Omega}^{1}(1)$ is no longer a module. 
We circumvent this problem by using  \v{C}ech cohomology
as   formulated in Subsection~\ref{schematic-defn}. 
Recall from Remark~\ref{well-defn-h} that 
$\HH^i(\Qgr S_R, \M) = \HH^i(\Qgr S,\M)$, which we write 
as $\HH^i(\M)$. By Theorem~\ref{cech cohomology}, $\HH^i(\M)$
can---and will---be computed as the $i$th
cohomology group  $\HH^i(\M\otimes_S\boldc^{\bullet})_0$.

 Given a module $\cM\in (\qgr S_R)_{\on{VC}}$, write $M=\Gamma^*(\cM)$ and
 consider the augmented  double complex  
\begin{equation}\label{beilinson11}
 M(-1)\otimes_S {\mathbf C}^{\bullet}
\otimes_S K^{\bullet}\longrightarrow
 M(-1)\otimes_S {\mathbf C}^{\bullet}\otimes_S \Delta\longrightarrow 0
\end{equation} 
As in \cite{OSS} and \cite{KKO} we have shifted $M$ by $(-1)$ to facilitate the
cohomological computations. We also note that, as $-\otimes {\mathbf C}^j$ 
kills finite dimensional modules,
$ \M(-1)\otimes_S {\mathbf C}^{\bullet}$ is well-defined (and equals 
 $M(-1)\otimes_S {\mathbf C}^{\bullet}$).
 
Since the sequence  $K^\bullet\to \Delta$ is a sequence of
graded $S$-bimodules,  each term in \eqref{beilinson11}
is also bigraded and so
we may take the degree zero summand ${}_0\{-\}$  under the left gradation and 
then take the image in $\Qgr S_R$ as right modules.  These last two
operations are exact functors and we will 
need to compute the cohomology of the 
spectral sequence obtained from the resulting double complex in $\Qgr S_R$:
\begin{equation}\label{beilinson44}
{\mathbf
C}^{\bullet,\bullet} = \pi (_0\{ M(-1)\otimes_S {\mathbf C}^{\bullet} \otimes_S
K^{\bullet}\}).
\end{equation}

We begin by considering the first filtration of this complex.
Thus, for fixed $i$, consider the $i$th
row of    \eqref{beilinson11}:
\begin{equation}\label{beilinson45}
M(-1)\otimes_S {\mathbf C}^{i} \otimes_S
K^{\bullet}\to M(-1)\otimes_S {\mathbf C}^{i}\otimes_S \Delta\to 0.
\end{equation}
As a left $S$-module,  $\Delta \cong \bigoplus_{j\geq 0}
S\Delta_{0j} \cong \bigoplus_{j\geq 0} S(j)_{\geq 0}$. Although $\Delta$ 
 is not a free left $S$-module, each summand  has
 finite codimension inside  $S(j)$ and so the factor $S(j)/S(j)_{\geq 0}$
is annihilated by tensoring with 
${\mathbf C}^i$. Hence, each ${\mathbf C}^i\otimes_S \Delta$ is a flat left
$S$-module. Using \eqref{omega1}  the same is true of ${\mathbf C}^i\otimes_S
\Omega^1$.  Therefore \eqref{beilinson45} is also  exact.
This shows that the spectral sequence associated to  \eqref{beilinson44}
has $E_1$ term:
\begin{equation}\label{beilinson89}
E_1^{pq} = \HH^q\big(\pi({}_0\{M(-1)\otimes_S {\mathbf C}^{\bullet} \otimes_S
K^p\})\big)  \Rightarrow 
\HH^{p+q}\big(\pi(_0\{M(-1)\otimes_S {\mathbf C}^{\bullet}\otimes_S
\Delta\})\big).
\end{equation}
We emphasize that this is a spectral sequence of objects in $\Qgr S_R$.

Consider the  right hand side of \eqref{beilinson89}.
For  any $r\in \mathbb Z$, note that 
 $_r\{M(-1)\otimes_S {\mathbf C}^{i}\} = M(-1)_t {\mathbf C}^{i}_{r-t}$
for  $t\gg 0$ and $_{-r}\Delta = \bigoplus_{j\geq 0} \Delta_{-r,j} =
\bigoplus_{j\geq 0} S_{j-r}$.  Thus, as
right $S$-modules
$$_0\{M(-1)\otimes_S {\mathbf C}^{i}\otimes_S \Delta\}
 \cong \bigoplus_{j\geq 0} \sum_{t,r} M(-1)_t{\mathbf C}^i_{r-t}S_{j-r} =
 M(-1){\mathbf C}^i_{\geq 0}. $$ 
By Theorem~\ref{cech cohomology},   
$\pi(_0\{M(-1)\otimes_S {\mathbf C}^{\bullet}\otimes_S \Delta\})$
therefore  has cohomology groups
  \begin{equation}\label{beilinson90}
\HH^i( \pi(_0\{M(-1)\otimes_S {\mathbf C}^{\bullet}\otimes_S \Delta\})) =
 {\mathbf H}^i\left(\M(-1)\right)_{\geq 0} =
 \bigoplus_{j\geq 0} \HH^i\left(\M(j-1)\right).
 \end{equation}
Here,  $\HH^i\left( \Qgr S_R,\M(j-1)\right)=0$
for $i>0$ and $j\gg 0$. 
Consequently, in $\Qgr S_R$, the 
only nonzero cohomology group in \eqref{beilinson90}
is 
$$\pi \HH^i( \pi(_0\{M(-1)\otimes_S {\mathbf C}^{\bullet}\otimes_S \Delta\}))=
\pi(M(-1))=\cM(-1).$$
Thus we have proved: 
\begin{equation}\label{beilinson70} 
\text{\em The spectral sequence  \eqref{beilinson89} 
converges to  $\cM(-1)$
concentrated
 in degree $0$.} 
 \end{equation}

We now compute all the terms $E^{pq}_1$ in \eqref{beilinson89}. 
By Condition~\ref{thesheaf},
$\HH^i\left(\M(-j)\right) = 0$ 
for  $i=0,2$ and $j=1,2$. Since $K^{-2}=S(-1)\otimes_k S(-2)$  and 
$K^0 = S\otimes_kS$, this means that 
$E^{pq}_1=0$ when both $q\in\{0,2\}$ and $p\in \{0,-2\}$.
By Lemma~\ref{chi-lemma}(6),
the cohomology is 
also zero for any $q\geq 3$ and by \eqref{diag-resn1},
it is  zero for $p\not\in \{0,-1,-2\}$.
Thus, it remains to consider the terms $E^{pq}$ when either $p$ or $q$ 
equals $1$.

\begin{sublemma}\label{beilinson middle term}
 Let  $M=\Gamma^*(\cM)$, as in the  proof of Theorem~\ref{beilinson5'}.
 Then:
\begin{enumerate}
\item 
$\HH^i(M(-1)\otimes_S \boldc^{\bullet}\otimes_S\Omega^1(1))=0$
for $i=0,2$.
\item $\HH^1(M(-1)\otimes_S \boldc^{\bullet}\otimes_S\Omega^1(1))$
is a finitely generated projective $R$-module.
\end{enumerate}
\end{sublemma}

\noindent
We postpone the proof of the sublemma until we have completed 
the proof of the theorem.
By Sublemma~\ref{beilinson middle term}(1)
and the observations beforehand,
$E^{pq}_1=0$ unless $q=1$. In other words, the spectral sequence
\eqref{beilinson89}
is  simply the  complex 
\begin{equation}\label{beilinson76}
0\to \HH^1\big( \M(-2)\big)\otimes_k{\mathcal O}_S(-2) 
\to V\otimes_k{\mathcal O}_S
\to \HH^1\big( \M(-1)\big)\otimes_k{\mathcal O}_S \to 0,
\end{equation}
where 
\begin{equation}\label{v0-defn}
V= \HH^1\big(\pi(_0\{M(-1)\otimes {\mathbf C}^\bullet\otimes
 \Omega^1(1)\})\big)\otimes _k{\mathcal O}_S(-1).  
 \end{equation}
 By the sublemma, $V$ is a finitely generated
 projective $R$-module, while $\HH^1(\M(-i))$ is a projective $R$-module 
 for $i=-1,-2$  by Remark~\ref{implied by VC}. 
Thus \eqref{beilinson70} implies that \eqref{beilinson76}
is a monad whose
cohomology is precisely $\cM(-1)$.
Shifting the  degree by $1$ therefore gives the desired complex
 $\bK(\M)$. 
   This completes the proof of  Theorem~\ref{beilinson5'}, 
   modulo the proof of the sublemma.\end{proof}

\begin{proof}[Proof of Sublemma~\ref{beilinson middle term}]
We first  prove the following assertion: 
\begin{equation}\label{beilinson-middle2}
\text{If } \cF\in\on{Gr-} S_R, \text{ then }
{\mathbf H}^i(\cF\otimes_S \boldc^{\bullet}\otimes_S\Omega^1(1)) =
\on{\underline{Ext}}^i_{\qgr S_R}\big((\widetilde{\Omega}^1(1))^{\vee},
 \pi \cF\big),
\end{equation}
where $\cN^\vee =\pi\uHom_{S\on{-qgr}}(N,\cO)$\label{M-vee-defn}
 for $\cN\in S\on{-qgr}$,
in the sense of \eqref{m-vee}.

By \eqref{omega1},
$\widetilde{\Omega}^1$ is a vector bundle and so
 $\boldc^{\bullet}\otimes\Omega^1(1)$ is a
complex of flat $S$-modules. Hence the functors 
\begin{equation}\label{cech delta functors}
M\mapsto {\mathbf H}^i\big(M\otimes_S\boldc^{\bullet}\otimes\Omega^1(1)\big)
\end{equation}
form a cohomological $\delta$-functor from $\Qgr S_R$ to $\Ggr R$. 
 If $I\in \Qgr S_R$ is injective then
 Theorem~\ref{cech cohomology} implies that 
$I\otimes_S \boldc^{\bullet}$ is exact in cohomological degrees greater than
zero.  Since $-\otimes_S\widetilde{\Omega}^1(1)$ is an exact functor,
 it follows that \eqref{cech delta functors}
is effaceable for all $i>0$. Thus
 \cite[Theorem~III.1.3A]{H} implies that 
the functors \eqref{cech delta functors} are the right derived functors
of  
$M\mapsto {\mathbf H}^0\big(M\otimes_S\boldc^{\bullet}\otimes\Omega^1(1)\big)$.

Recall the  localizations $Q_w$ defined in  Subsection~\ref{schematic-defn}.
If $\cF\in \gr S_R$  and  $N\in S_R\on{-qgr}$, 
then there is a natural map
\bd
\Psi_{\cF,N}: \cF\otimes Q_w\otimes N\rightarrow
\on{\underline{Hom}}_{\Qgr S_R}(N^{\vee},\pi \cF\otimes Q_w)
\ed
of graded $R$-modules
given by
$
f\otimes q\otimes n \mapsto \big[ \phi\mapsto f\otimes q\cdot \phi(n)\big].
$
The map $\Psi_{\cF,N}$ is easily seen to be compatible with homomorphisms
 $N\rightarrow  N'$ and
with the maps $Q_w\rightarrow Q_w'$ in the \v{C}ech complex.
Observe that, if $N=\theo_{S_R}(k)$ for some $k$, then $\Psi_{\cF,N}$ is an
isomorphism: indeed, 
$N^{\vee} = \theo_{S_R}(-k)$ and so $\Psi_{\cF,N}$ reduces to the isomorphism
$\cF\otimes Q_w (k) = \uHom_{\Qgr S_R}(\theo_{S_R}(-k),\cF\otimes Q_w)$.

Dualizing  \eqref{omega1} gives an exact
sequence
$
0\rightarrow \theo \rightarrow \theo(1)\otimes S_1^*\rightarrow
(\widetilde{\Omega}^1)^{\vee}\rightarrow 0
$
in $S\on{-qgr}$.  After tensoring with $R$, which
we suppress in our notation, this gives the exact sequence
\begin{equation}\label{Homs for Omega}
0\rightarrow \uHom\big((\widetilde{\Omega}^1)^{\vee},\pi \cF\otimes Q_w\big)
\rightarrow \uHom\big(\theo(1)\otimes S_1^*,\pi \cF\otimes
Q_w\big)\rightarrow \uHom\big(\theo,\pi \cF\otimes Q_w\big),
\end{equation}
where all the Homs are in $\Qgr S_R$.
On the other hand, by \eqref{omega1}, the sequence
\begin{equation}\label{tensors for Omega}
0\rightarrow \cF\otimes Q_w\otimes \widetilde{\Omega}^1
\rightarrow \cF\otimes Q_w\otimes\theo(-1)\otimes S_1
\rightarrow \cF\otimes Q_w\otimes \theo
\end{equation}
is exact since $\theo_{S_R}$ is flat in $\qgr S_R$.  The maps 
$\Psi_{\cF,\bullet}$ give
a homomorphism of complexes from 
\eqref{tensors for Omega} to \eqref{Homs for Omega} 
that is an isomorphism in the middle
and right-hand columns. Thus $\Psi_{\cF,\Omega^1}$ is
also an isomorphism.  Since this is true for every $Q_w$ compatibly
with the maps in the \v{C}ech complex, Theorem~\ref{cech
cohomology}
implies that 
$\uHom_{\Qgr S_R}\big((\widetilde{\Omega}^1)^{\vee}, \pi\cF\big) = 
{\mathbf H}^0\big(\cF\otimes \boldc^{\bullet}\otimes \Omega^1\big)$ 
for every graded 
$S_R$-module $\cF$.  Thus 
$\cF \mapsto \uHom_{\Qgr S_R}\big((\widetilde{\Omega}^1)^{\vee}, 
\pi\cF\big)$ and 
$\cF\mapsto {\mathbf H}^0\big(\cF\otimes \boldc^{\bullet}\otimes 
\Omega^1\big)$ have the
same right-derived functors and \eqref{beilinson-middle2} follows.

 We now turn to the proof of the sublemma.
The significance of \eqref{beilinson-middle2}
is that we can apply Theorem~\ref{CBC}. Indeed, by parts (1) and (3) of
that theorem, the sublemma will be immediate if we can  prove
that 
\begin{equation}\label{beilinson-middle3}
\HH^i(M(-1)\otimes k(p)\otimes_S \boldc^{\bullet}\otimes_S \Omega^1(1)) = 
0\qquad \text{for}\ i=0,2\ \text{and  } p\in\spec R.
\end{equation}

In order to prove \eqref{beilinson-middle3},
set $M'=M\otimes k(p)$, for some $p\in \spec R$. 
As before,  ${\mathbf C}^\ell$ kills finite dimensional
left $S$-modules and
 ${\mathbf C}^\ell\otimes S(j)$ is
a flat left $S$-module.
Thus,  tensoring the shift of
\eqref{omega1} with $M'(-1)\otimes_S{\mathbf C}^\bullet $
gives the exact sequence
\begin{equation}\label{beilinson73}
0\rightarrow M'(-1)\otimes\boldc^{\bullet}\otimes\Omega^1(1)
\rightarrow M'(-1)\otimes\boldc^{\bullet}\otimes S\otimes_k S_1
\rightarrow M'(-1)\otimes\boldc^{\bullet}\otimes S(1) \rightarrow 0.
\end{equation}
Taking homology  gives the exact sequence
 \begin{equation}\label{beilinson71}
0\to \HH^0\big(M'(-1)\otimes_S{\mathbf C}^{\bullet}\otimes \Omega^1(1)\big)
\to
\HH^0\big(M'(-1)\otimes_S{\mathbf C}^{\bullet}\otimes_S S\otimes_k S_1 
\big).
\end{equation}
By Theorem~\ref{cech cohomology}, the final term of \eqref{beilinson71} 
equals
 $ \HH^0\big(\M(-1)\otimes k(p)\big)\otimes_k S_1 $
which, by  Condition~\ref{thesheaf}, is zero. Thus \eqref{beilinson-middle3}
holds for $i=0$.

By construction, $\Omega^1(1)$
can also be included in the exact sequence
\begin{equation}\label{beilinson72}
0\to S(-2) \to S(-1)\otimes (S_2^!)^*
\to {\Omega}^1(1) \to 0.\end{equation}
Now, $\HH^3(\Qgr S,-)=0$. Thus,
if we tensor  \eqref{beilinson72} with  $M'(-1)\otimes_S{\mathbf C}^\bullet 
$
and take homology,
we obtain the exact sequence
$$\HH^2\big(M'(-1)\otimes_S {\mathbf C}^{\bullet}
  \otimes_S S(-1)\otimes_k (S_2^!)^* \big)
  \to \HH^2\big(M'(-1)\otimes_S {\mathbf C}^{\bullet}\otimes_S 
\Omega^1(1)\big)
  \to 0.$$
But now the first term equals $ \HH^2\big(\M\otimes k(p)(-2)\big)$ which, 
by
 Condition~\ref{thesheaf}, is zero.  Thus \eqref{beilinson-middle3}
holds for $i=2$ and so  \eqref{beilinson-middle3} is true.
This therefore 
completes the proof of both the sublemma and Theorem~\ref{beilinson5'}. 
\end{proof}

We next show that the functor  $\cM\mapsto \bK(\cM)$
 induces an equivalence of categories.

\begin{thm}\label{beilinson5}
Suppose that $S\in \ASthree$  and that
$R$ is a commutative noetherian $k$-algebra. Then:
\begin{enumerate}
\item[(1)] The functor $\HH^0$ induces an equivalence of categories from 
$\on{Monad}(S_R)$ 
to $(\qgr S_R)_{\on{VC}}$, with inverse $\cM\mapsto \bK(\cM)$.
\item[(2)] Let $R\rightarrow R'$ be a homomorphism of commutative noetherian
$k$-algebras. Then, as functors from $\on{Monad}(S_R)$ to $(\qgr
S_R)_{\on{VC}}$, the composite functors 
$(-\otimes_RR')\circ \HH^0$ and $\HH^0\circ(-\otimes_RR')$
are naturally equivalent. 
\end{enumerate}
\end{thm}

Before proving the theorem, we need several lemmas.

\begin{lemma}\label{beilinson5-1}
Let $R$ be a commutative noetherian $k$-algebra and $S\in \ASthree$.
  If $\bK$ is a monad for $S_R$, then
$\M=\HH^0(\bK)$ is $R$-flat. Moreover,
$\HH^i(\bK)\otimes_R N = \HH^i(\bK\otimes_R N)$
for every $R$-module $N$ and every $i$.
\end{lemma}

\begin{proof}
Keep the notation of \eqref{monad0} and set $\cB =\ker(B_{\bK})$.
From the exact sequence
$
0\rightarrow \cB\rightarrow \theo\otimes
V_0\rightarrow \theo(1)\otimes V_1\rightarrow 0
$
it is clear that $\cB$ is $R$-flat.  It follows from
\cite[Proposition~C1.4]{AZ2} that
$\on{ker}(B_{\bK})\otimes N = \on{ker}(B_{\bK}\otimes N)$.

  Now consider
the exact sequence
$
0\rightarrow \theo(-1)\otimes V_{-1}\rightarrow \cB
\rightarrow \M \rightarrow 0.
$
By Lemma~\ref{chi-lemma}(3), this induces the  exact sequence of
$R$-modules
$$
0\longrightarrow \HH^0(\theo(-1+n)\otimes V_{-1})
\ {\buildrel{\theta}\over{\longrightarrow}}\
\HH^0(\cB(n))\longrightarrow \HH^0(\M(n))\longrightarrow 0
$$
for all $n\gg 0$.
In order to prove that $\cM$ is $R$-flat, it suffices to prove
that $\HH^0(\M(n))$
is $R$-flat for all $n$ sufficiently large \cite[Lemma~E5.3]{AZ2}. Thus,
by  \cite[Theorem~6.8]{Ei}, it suffices to prove that
$\theta\otimes_Rk(p)$ is injective
for all $n\gg 0$ and all $p\in\spec R$.
Theorem~\ref{CBC}(3) implies that $\theta\otimes_Rk(p)$
 is  the natural map
\bd
\HH^0\left(\theo(-1+n)\otimes V_{-1}\otimes k(p)\right)\rightarrow
\HH^0\left(\cB(n)\otimes k(p)\right) =
\HH^0\left(\ker(B_{\bK\otimes k(p)})(n)\right).
\ed
This map is injective since $A_{\bK\otimes k(p)}$
is injective.  Thus $\M=\HH^0(\bK)$ is  $R$-flat.

By \cite[Proposition~C1.4]{AZ2}, again, it follows that
\bd
0\rightarrow \theo(-1)\otimes V_{-1}\otimes
 N\rightarrow \on{ker}(B_{\bK}\otimes N)
\rightarrow \M\otimes N\rightarrow 0
\ed
is exact for any $R$-module $N$.
  In particular, $\M\otimes N = \HH^0(\bK\otimes N)$ and
$A_{\bK\otimes N}$ is injective.  Since $B_{\bK}$ is surjective,
$B_{\bK\otimes N}$ is also surjective.
 \end{proof}

\begin{lemma}\label{monads give VC modules}
Keep the hypotheses of Theorem~\ref{beilinson5}.
Let $\bK$ be a monad and set  $\M=\HH^0(\bK)$.  
Then for every $p\in\spec R$, $\HH^0(\M(-i)\otimes k(p))=0$
for $i\geq 1$ and $\HH^2(\M(j)\otimes k(p))=0$ for $j\geq -2$.
\end{lemma}

\begin{proof} 
By Lemma~\ref{beilinson5-1},
$\M\otimes k(p) = \HH^0(\bK)\otimes k(p) = \HH^0(\bK\otimes k(p))$
for any point $p$ of $\spec R$.
Thus it suffices to prove the claim for a monad $\bK$ in $\qgr S$.
Keep the notation of \eqref{monad0}.

Write $\cK =\ker(B_{\bK})$, in the notation of \eqref{monad0}.
Then the   exact sequence
\bd
0\rightarrow \theo(-1-i)\otimes V_{-1}\rightarrow 
\cK(-i) \rightarrow \M(-i)\rightarrow 0
\ed
induces an exact sequence
$
\HH^0\left(\cK(-i)\right)\rightarrow \HH^0(\M(-i))
\rightarrow \HH^1(\theo(-1-i)).
$
Since $\cK(-i)\subseteq V_0\otimes \theo(-i)$, 
one has $\HH^0\left(\cK(-i)\right)=0$ for all $i>0$.
On the other hand, $\HH^1(\theo(-1-i))=0$ for all $i$, by 
Lemma~\ref{chi-lemma}(5). Thus, $\HH^0(\M(-i))=0$ for $i\geq 1$. 
By Lemma~\ref{chi-lemma}(5),
$\HH^2(\theo(j)\otimes V_0)=0$ for $j\geq -2$ while 
$\HH^3(\cN)=0$ for any $\cN\in \qgr S$.
Thus the long exact sequence in homology associated to the  surjection 
$\theo(j)\otimes V_0\twoheadrightarrow \on{coker}(A_{\bK})(j)$
 shows that $\HH^2\left(\on{coker}(A_{\bK})(j)\right)=0$
for $j\geq -2$.  

Finally, the  exact sequence
$$
0\rightarrow \M(j)\rightarrow \on{coker}(A_{\bK})(j)\rightarrow 
\theo(1+j)\otimes V_1\rightarrow 0
$$
induces the exact sequence
\bd
\HH^1(\theo(1+j)\otimes V_1)\rightarrow \HH^2(\M(j))\rightarrow 
\HH^2\left(\on{coker}(A_{\bK})(j)\right).
\ed
By Lemma~\ref{chi-lemma}(5), respectively 
the conclusion of the last paragraph, the outside terms vanish for $j\geq -2$. 
Thus, the middle term is also zero.
\end{proof}

The last two lemmas   imply that $\HH^0$ maps $\on{Monad}(S_R)$ 
to $(\qgr S_R)_{\on{VC}}$.
We next want to prove that $\HH^0$ is a full and faithful functor on monads,
which will complete the proof of  Theorem~\ref{beilinson5}(1). 
This follows from the following more general fact.

\begin{lemma}\label{exts agree} 
Keep the hypotheses of Theorem~\ref{beilinson5} and 
suppose that $\bK$ is a monad and $\bL$ is a complex of the form
\bd
\bL:\;\;\; W_{-1}\otimes\theo(-1)\rightarrow W_0\otimes\theo \rightarrow
W_1\otimes\theo(1)
\ed
with each $W_i$ a finitely generated $R$-module and $\HH^n(\bL)=0$
for $n\neq 0$.
 Let $\cE = \HH^0(\bK)$ and $\cF = \HH^0(\bL)$.  
Then for any $i$, $\Ext^i(\cE, \cF) = \HH^i(\Hom^{\bullet}(\bK,\bL))$.
\end{lemma}

\begin{proof}
Let $L^n$ denote the $n$th term in $\bL$ and choose an injective 
resolution $L^n\rightarrow I^n_{\bullet}$ for each
$n$; these combine to give a double complex $(I^{\bullet}_{\bullet})$
 the total complex of which is an injective resolution    
of $\bL$.  Let $d_1$ and $d_2$ denote the two differentials 
(the first ``in the $\bL$-direction'' and the second 
``in the $I$-direction''), with signs adjusted as usual so that
$d_{1}d_{2} +d_{2}d_{1} = 0$.  Let $d_1$ denote the given differential on 
$\bK$ and define a second differential $d_{2}$ 
on $\bK$ to be zero everywhere. 

Define a double complex by setting 
$\cC^{p,q} =\bigoplus_n \Hom(K^n, I^{n+p}_q)$ and defining 
differentials by $\delta_i(f) = d_i\circ f - (-1)^{\deg (f)}f\circ d_i$. 
 The associated total complex then satisfies
$\on{Tot}(\cC^{\bullet,\bullet}) = 
\Hom^{\bullet}\left(\bK, \on{Tot}(I^{\bullet}_{\bullet})\right)$, 
which yields
$\HH^i(\on{Tot}(\cC^{\bullet, \bullet})) = \Ext^i(\cE,\cF)$ 
since $\bK$ is quasi-isomorphic to $\cE$.  On the other hand,
the double complex has associated spectral sequence with 
$E_2^{\bullet,\bullet} = H_I(H_{I\! I}(\cC^{\bullet,\bullet}))$.  
Now 
\bd
H_{I\! I}^{p,q}(\cC^{\bullet,\bullet}) = \bigoplus_n\Ext^q(K^n,L^{n+p})
=\bigoplus_{n=-1}^1 \HH^q(\theo(p))\otimes V_n^*\otimes W_{n+p} ,
\ed 
which vanishes by hypothesis when $q\neq 0$
and satisfies $H_{I\! I}^{p,0} = \bigoplus_n\Hom(K^n, L^{n+p})$. 
 The spectral sequence thus degenerates at $E_2$, completing the proof.
\end{proof}

\begin{proof}[Proof of Theorem~\ref{beilinson5}]
   Lemmas~\ref{beilinson5-1} and \ref{monads give VC modules}
show that  $\HH^0$  maps $\on{Monad}(S_R)$ 
to $(\qgr S_R)_{\on{VC}}$, while  Theorem~\ref{beilinson5'}
shows that $\HH^0$ is surjective on objects.  Lemma~\ref{exts agree}
implies that $\HH^0$ is fully faithful, which by \cite[Theorem~I.5.3]{Po} 
is sufficient to prove part (1).
Part (2) of the theorem is immediate from Lemma~\ref{beilinson5-1}.
\end{proof}

\section{Semistable Modules and Kronecker Complexes}\label{section3}

We  want to use the results of the previous section to construct a  projective
coarse moduli space for modules over AS regular rings. As   will be seen in
Lemma~\ref{vanishing cohomology for ss modules}, the construction of
Theorem~\ref{beilinson5} applies to  a large class of flat families of modules
including  all torsion-free modules of rank one, up to a shift. To construct
the moduli space as  a projective scheme, it will then suffice to give a
convenient realization of a parameter space for  monads.  The condition that
the pair of maps  $A_K$ and  $B_K$  in \eqref{monad0} actually define a monad
rather than just  a complex is awkward to describe directly in terms of linear
algebra. The standard way round this difficulty \cite{Drezet-Le Potier} is to
characterize monads as so-called Kronecker complexes satisfying an appropriate
stability condition.  This approach also works in our  noncommutative setting
and the details are given in  this section.

\subsection{Modules
 in $\qgr S$}\label{families of modules and Kronecker complexes}
There are several invariants that can be attached to a  module $\cM\in \qgr
S$ for $S\in \ASthree$ and we begin by describing them. 

As before, we write $\theo_S=\pi(S)\in \qgr S$.
The {\it Euler characteristic}\label{euler-defn} of $\cM$
is defined  to be $\chi(\cM)=\sum (-1)^{i} h^i (\cM)$, where $h^i(\cM)=\dim
\HH^i(\cM)$, and the {\it Hilbert polynomial}\label{hilbert-defn} of $\cM$ is
$p_{\cM}(t) = \chi(\cM(t))$. 
 If $\M$ is nonzero
and torsion-free, the {\em normalized Hilbert polynomial}  of $\M$ is
$p_{\M}(t)/\on{rk}(\M)$. The other invariant that will be frequently used is the
first Chern class, as defined on page~\pageref{chern-defn}.
When $\qgr S\simeq \coh \boldp^2$,   these definitions 
coincide with the usual commutative ones.

The next lemma gives some standard properties of these invariants. 

\begin{lemma}\label{hilbert polynomial}
{\rm (1)}  If $\M\in \qgr S$ then
$p_\M$ is a polynomial.  The Hilbert polynomial is additive in exact 
sequences and has positive leading coefficient.

{\rm (2)} If
$\M\in \qgr S$ then
\begin{equation}\label{RR}
p_{\M}(t) = \frac{1}{2}\on{rk}(\M)t(t+1) + \big(c_1(\M)+ \on{rk}(\M)\big)t + 
\chi(\M).
\end{equation}
\end{lemma}

\begin{proof}  To begin with, assume that $\cM=\theo_S(m)$ for  
some $m\in \mathbb Z$. Since $S\in \ASthree$, $S$ has the Hilbert series of a
polynomial ring in $3$ variables and so  
$\HH^0(\qgr S, \theo_S(m))=\HH^0(\mathbf P^2, \theo_{\mathbf P^2}(m))$ 
for all $m$. By  Lemma~\ref{chi-lemma}(5) and 
Serre duality (Proposition~\ref{serre-duality}),
$\HH^i(\qgr S, \theo_S(m))=\HH^i(\mathbf P^2, \theo_{\mathbf P^2}(m))$, for 
$i\geq 1$ and $m\in \mathbb Z$.  Thus, by
Riemann-Roch, as stated in \cite[p.154]{LP2},  
the lemma   holds for $\cM=\theo(m)$.

Now let $\cM$ be arbitrary.  As usual, additivity of the Hilbert polynomial on
exact sequences follows from the long exact cohomology sequence. Thus, both
sides of \eqref{RR} are additive on exact sequences. 
Since $S$
has finite global dimension, every finitely generated graded $S$-module admits
a finite free resolution. By additivity it therefore 
suffices to prove the lemma for
$\cM = \theo_S(m)$, as we have done.  \end{proof}

\begin{corollary}\label{two term}
Let $
\chi(\cE,\cF) = \sum_i (-1)^i \dim \Ext^i(\cE,\cF)$
for $\cE,\cF\in \qgr S$.
Then
\begin{equation}\label{two term euler char}
\chi(\cE,\cF) = \on{rk}(\cE)[\chi(\cF)-\on{rk}(\cF)]
-c_1(\cE)[3\on{rk}(\cF)+c_1(\cF)] + \chi(\cE)\on{rk}(\cF).
\end{equation}
\end{corollary}
\begin{proof} The formula reduces to \eqref{RR} when $\cE=\theo(k)$ and
 $\cF = \theo(\ell)$.  Since both sides of
\eqref{two term euler char} are separately additive in $\cE$ and $\cF$,
the formula follows by taking resolutions of both $\cE$ and $\cF$ by
direct sums of $\theo(n)$'s.
\end{proof}

A torsion-free module  $\cM\in \qgr S$ is 
{\it normalized}\label{normalized-defn}
 if $-\on{rk}(\cM)< c_1(\cM)\leq 0$.
If $\cM$ is torsion-free then Lemma~\ref{hilbert polynomial'} implies that
$c_1(\cM(1))=c_1(\cM)+\on{rk}(\cM)$ and so 
there is a unique normalized  shift of $\cM$. 
In particular, a torsion-free,  rank
1 module $\cM$ is normalized if and only if $c_1(\cM)=0$.

Order polynomials lexicographically and write {\em (semi)stable} to mean
``stable, respectively semistable''. 
A torsion-free module $\cM\in \qgr S$ is {\em
(semi)stable}, \label{stability definition for modules} if, for every proper
  submodule $0\not= \cF\subset \cM$, one has $\on{rk}(\cM)p_{\cF} -
\on{rk}(\cF)p_{\cM} < 0$ (respectively $\leq 0$). $\M$ is called {\em
geometrically (semi)stable} if $\M\otimes \overline{k}$ is (semi)stable where
$\overline{k}$ is an algebraic closure of $k$. $\cM$ is called  {\em
$\mu$-(semi)stable} if $\on{rk}(\cM)c_1(\cF) - \on{rk}(\cF)c_1(\cM) <0$
(respectively $\leq 0$).

If $\cF\subset\cE$ is de-(semi)stabilizing, then
$\cF\otimes\overline{k}\subset \cE\otimes\overline{k}$ is 
de-(semi)stabilizing and so a geometrically (semi)stable module
is (semi)stable.
One also has the standard implications
\bd
\mu\text{-stable}\; \Rightarrow\; \text{stable}\; \Rightarrow\; 
\text{semistable}\; \Rightarrow\; \mu\text{-semistable.}
\ed

As Lemma~\ref{vanishing cohomology for ss modules} will
 show, semistable modules have tightly controlled cohomology.
The argument uses the following standard consequence of
 \cite[Theorem 3]{Rudakov}.

\begin{lemma}\label{JH}
Let $\cM\in \qgr S$ be  torsion-free and semistable. Then $\cM$ admits
a \emph{Jordan-Holder filtration} in the   sense that 
there exists a filtration
\bd
\{0\} = \cM_0\subsetneq \cM_1 \subsetneq\dots \subsetneq \cM_k = \cM
\ed
such that
\begin{enumerate}
\item each $\on{gr}_i(\cM)=\cM_i/\cM_{i-1}$ is a torsion-free  
stable  module in $\qgr S$, and
\item $p_{\cM}/\on{rk}(\cM) = p_{\on{gr}_i(\cM)}/\on{rk}(\on{gr}_i(\cM))$
for every $i$. \qed
\end{enumerate}
\end{lemma}

\begin{lemma}\label{vanishing cohomology for ss modules} Suppose that 
$\cM\in \qgr S$ is  torsion-free. Then:
\begin{enumerate} 
\item[(1)] If $\cM$ is
$\mu$-semistable and normalized,  then $\HH^0(\cM(-i))=0$ for $i\geq 1$
and $\HH^2(\cM(i))= 0$ for
$i\geq -2$. 
\item[(2)] If $\cM$ is semistable and normalized, then either 
$\cM\cong \theo^{r}$
for $r=\on{rk}(\cM)$ or 
$\HH^0(\cM)= 0$. 
\end{enumerate} 
\end{lemma} 
 
\begin{proof} Although the proof is very similar to that of 
\cite[Lemme~2.1]{Drezet-Le Potier} we will give it since it is 
fundamental to our approach.

(1) If $\HH^0(\cM(-i))\not= 0$ for some $i>0$,
then $\theo_S(i)\hookrightarrow \cM$.
As  $c_1(\cM)\leq 0$, 
 $\mu$-semistability forces $$\on{rk}(\cM)c_1(\theo(i))\ \leq \
\on{rk}(\cM)c_1(\theo(i))-i\cdot c_1(\cM) \leq 0,$$ a contradiction. 
If $\HH^2(\cM(-3+i))\not= 0$ for some $i\geq 1$ then Serre duality
(Proposition~\ref{serre-duality}) implies that there exists
$0\not=\theta\in \on{Hom}(\cM,\theo(-i))$. It is then a simple exercise to see 
that $\cF=\ker(\theta)$ contradicts the $\mu$-semistability of $\cM$.

(2)  As $\cM$ is $\mu$-semistable,
 part (1) implies that  $p_{\cM}(-1)=\chi(\cM(-1))\leq 0$.  
By \eqref{RR},
this forces $\chi(\cM) \leq c_1(\cM)+\on{rk}(\cM)$. 
 Now assume that  $h^0(\cM)\neq 0$; thus
 $\theo_S\hookrightarrow \cM$.
 By $\mu$-semistability and the fact that $c_1(\theo_S)=0$ 
 we obtain $c_1(\cM)=0$. Hence
$\chi(\cM)\leq \on{rk}(\cM)$.  Conversely, semistability 
implies that $\on{rk}(\cM) p_{\theo_S}-p_{\cM}\leq 0$.
Substituting this into \eqref{RR} shows that $\chi(\cM)\geq \on{rk}(\cM)$.
Thus, $\chi(\cM)= \on{rk}(\cM)$ and so $\cM$ has normalized Hilbert
polynomial equal to that of $\theo_S$.

By Lemma~\ref{JH},  choose a Jordan-H\"older filtration $\{\cM_i\}$
 of $\cM$.
By the last paragraph, each $\on{gr}_i\cM$
 is stable and  has normalized Hilbert
polynomial equal to that of $\theo$. Also,  
$\chi(\on{gr}_i\cM)>0$ by Lemma~\ref{JH}(2) but $\HH^2(\on{gr}_i\cM)=0$ 
by part (1) of this proof. Thus, $\on{gr}_i\cM$ has a nonzero global
section. The corresponding inclusion 
$\theo\hookrightarrow \on{gr}_i\cM$
contradicts the stability of $\on{gr}_i\cM$
unless $\on{gr}_i\cM\cong \theo$ for each $i$. 
Finally,  
Lemma~\ref{chi-lemma}(5) implies that 
$\Ext^1_{\qgr S}(\theo,\theo)=0$
and hence that  $\cM\cong \theo^r$.
\end{proof}

\begin{remark}\label{chi ineq}
The proof of the lemma shows that, if 
$\M$ is torsion-free, semistable and normalized then 
$\chi(\M)\leq c_1(\M)+\on{rk}(\M)$.
\end{remark}

Combined with  the Beilinson spectral sequence, the observations of this
section have strong consequences for the  cohomology of rank one torsion-free
modules.

\begin{corollary}\label{monad6}
Let $S=S(E,\cL,\sigma)\in \ASthree'$.
 Assume that $\cM\in \qgr S$ is a torsion-free, rank one module with 
$c_1(\M)=0$. Then:
\begin{enumerate}
\item  $\cM$ is stable and $\cM \in (\qgr S)_{\on{VC}}$, as defined
 in Condition~\ref{thesheaf}.
\item  $\dim_k\HH^1(\cM(-1)) = \dim_k\HH^1(\cM(-2))=n$, where $n=1-\chi(\cM)$.
\item Suppose that $\cM\in \cV_S$; that is, suppose that $\cM|_E$ is a vector bundle.  Then 
 $\cM|_E \cong \left(\cL\otimes 
(\overline{\cL})^{-1}\right)^{\otimes n}$ where 
 $\overline{\cL}=\cL^{\sigma^{-1}}$ and $n=1-\chi(\cM)$.
\end{enumerate}
\end{corollary}

\begin{proof} (1) If $0\not= \cF\subsetneq \cM$
then $\on{rk}(\cF)=\on{rk}(\cM)$ but 
$p_\cF<p_\cM$, simply because $p(\cM/\cF)>0$. Thus, $\cM$ is stable. 
By Lemma~\ref{vanishing cohomology for ss modules},
$\cM \in (\qgr S)_{\on{VC}}$. 

(2) By part (1) and 
Theorem~\ref{beilinson5'}, 
$\cM$ is the cohomology of a monad, which we assume has the form \eqref{monad0}.
By the additivity of $c_1$ on exact sequences, 
$$0=c_1(\M)= c_1(\theo(-1)\otimes V_{-1})+c_1(\theo(1)\otimes V_{1})
=\dim V_{1}-\dim V_{-1}.$$  Equivalently,
$\dim_k\HH^1(\cM(-1)) = \dim_k\HH^1(\cM(-2))$.
  Lemma~\ref{chi-lemma}(5)
and the additivity of $\chi$ on exact sequences 
then imply that $\chi(\cM)=1-\dim V_1$.

(3)  By Lemma~\ref{shifts}, the restriction of \eqref{monad0} to $E$
is a complex of the form 
\begin{equation}\label{diagg1}
0\rightarrow V_{-1}\otimes \cL^* \rightarrow V_0\otimes\theo_E 
\rightarrow V_1\otimes\baL \rightarrow 0.
\end{equation}
 Equation~\ref{monad0}
can be split into two short exact sequences of torsion-free modules and
by  Lemma~\ref{restriction of
torsion-free} the restrictions of those exact sequences to $E$ are 
again exact. Thus \eqref{diagg1} is  a complex whose only cohomology
 is the vector bundle  $\cM|_E$ in degree zero.
 Part (3) therefore  follows by taking
determinants of  \eqref{diagg1}.  \end{proof}

\subsection{Kronecker Complexes}\label{kronecker-defn}
We continue to assume that $S\in \ASthree$ in this subsection. 
To treat moduli of monads, it is more convenient
 to work with (a priori) more general complexes that satisfy suitable 
 stability properties. In this subsection we define the 
 relevant complexes, called Kronecker complexes,  and 
show that semistability forces them to  be monads.
Our treatment closely follows \cite{Drezet-Le Potier}.

\begin{defn} 
A {\em Kronecker complex} in $\qgr S$ is a complex of the form
\begin{equation}\label{typical Kronecker complex}
{\mathbf K}: \;\;\; \theo(-1)\otimes V_{-1}\xrightarrow{A} \theo\otimes V_0 
\xrightarrow{B} \theo(1)\otimes V_1
\end{equation}
where $V_{-1}$, $V_0$ and $V_1$ are finite-dimensional vector spaces. 
Clearly a monad is a special case of a Kronecker complex.
We index the complex so that $\theo(i)$ occurs in cohomological degree $i$;
 thus  $\HH^1(\mathbf K)$ denotes the homology at $\theo(1)\otimes V_1$.
\end{defn}

Morphisms of Kronecker complexes are just   morphisms of complexes, 
and so  are defined by maps of the 
defining vector spaces $V_i$.   It is then 
easy to see that the category of Kronecker complexes in $\qgr S$ is 
 an abelian category.

The invariants we   just defined for modules have their natural
 counterparts for Kronecker complexes.
Thus, if $\bK$ is a Kronecker complex  as in \eqref{typical Kronecker
complex},   the  {\em rank}\label{kronecker-rank}
 of  $\bK$ is defined to be 
  $\on{rk}(\bK) = \on{rk}(V_0) - \on{rk}(V_{-1})-\on{rk}(V_1)$.
This can be negative. The {\em first Chern
  class}\label{kronecker-chern} $c_1(\bK)$
is $c_1(\bK) = \on{rk}(V_{-1})-\on{rk}(V_1)$ and the 
{\em Euler characteristic}\label{euler-kronecker}
$\chi(\bK)$ is $\chi(\bK) = \on{rk}(V_0)-3\on{rk}(V_1)$.  The {\em Hilbert
polynomial}\label{hilbert-kronecker} $p_\bK$ of $\bK$ is given by
the formula
\bd
p_\bK = \on{rk}(V_0)\cdot p_{\theo} - \on{rk}(V_{-1})
\cdot p_{\theo(-1)} - \on{rk}(V_1)\cdot p_{\theo(1)}.
\ed
The {\em normalized Hilbert polynomial} of a Kronecker complex $\bK$ of
positive rank is $p_{\bK}/\on{rk}(\bK)$. If $\bK$ is a monad,
  all these invariants coincide with the corresponding 
invariants of the cohomology  $\HH^0(\bK)$.

The correct notion of (semi)stability for Kronecker complexes is the following.

\begin{defn}\label{semistability 
definition for complexes}  A Kronecker complex
$\bK$ is  {\em (semi)stable} if, for every proper 
subcomplex $\bK'$ of $\bK$, one has
$\on{rk}(\bK)p_{\bK'} - \on{rk}(\bK')p_{\bK} <0$ (respectively $\leq 0$)
 under the lexicographic order on polynomials.  
Equivalently:
\begin{enumerate}
\item[(1)]
$\on{rk}(\bK)c_1(\bK') - \on{rk}(\bK')c_1(\bK) \leq 0$, and
\item[(2)] if   equality holds in (1) then 
$\on{rk}(\bK)\chi(\bK') - \on{rk}(\bK')\chi(\bK) <0$ (respectively $\leq 0$).
\end{enumerate}
A Kronecker complex $\bK$ is {\em geometrically (semi)stable}\label{geom-Kr}
if $\bK\otimes \overline{k}$ 
is (semi)stable where $\overline{k}$ is an algebraic closure of $k$.
A Kronecker complex $\bK$ of rank $r>0$ is 
{\em normalized}\label{norm-complex} if $-r<c_1(\bK)\leq 0$.
\end{defn}

The following proposition is the main result of this section. 
It obviously  allows us to replace monads by semistable Kronecker
complexes in describing $S$-modules and this will be  important since it is much
easier to describe the moduli spaces of Kronecker complexes than those of
monads.

\begin{prop}\label{semistable Kronecker complex is monad}
Let $S\in \ASthree$.
Suppose that $\bK$ is a semistable normalized Kronecker complex.  Then
$\bK$ is a torsion-free 
monad in the sense of Definition~\ref{monad-defn}.
\end{prop}

 Before proving the proposition, we need several preliminaries.
 Since we are proving a result about $\qgr S$, it suffices to prove 
 the result when $S=S(E,\cL,\sigma)\in \ASthree'$ and we assume 
 this throughout the proof.

The first three results closely follow the strategy
from \cite{Drezet-Le Potier}
where similar results are proved in the commutative case. 

\begin{lemma}\label{enumeration of Kronecker complexes}  
Let $S\in \ASthree'$ and suppose that  $\bK$ is a Kronecker complex 
of the form \eqref{typical Kronecker complex}.  
\begin{enumerate}
\item[(a)]\label{subcomplexes first part} If $B|_E$ is not surjective then
 $\bK$ has a quotient complex of one of the forms:
\begin{alignat*}{6}
(3)\quad\quad\quad\quad  &\theo(-1) &\ \rightarrow\
 & \theo^2 & \ \rightarrow\ &\theo(1) & 
\quad\quad & \text{which is exact at } \theo^2\\
(4)\quad\quad\quad\quad &\qquad 0 & \ \rightarrow\ & \theo &\ \rightarrow\
 &\theo(1) & \\
(5)\quad\quad\quad\quad  &\qquad  0 & \ \rightarrow\
 & 0 &\ \rightarrow\  &\theo(1) &  \\
(6)\quad\quad\quad\quad  & \qquad 0 & \ \rightarrow\ & \theo^2 &\ \rightarrow\
 &\theo(1) & 
\end{alignat*}
\item[(b)]\label{subcomplexes second part} Suppose either that 
$\theo(-1)|_E\otimes V_{-1}\xrightarrow{A|_E} \theo_E\otimes V_0$ is not
injective or that $\on{coker}(A|_E)$ has a simple subobject.
Then $\bK$ has a subcomplex  either of type (3) or of one of the forms
\begin{alignat*}{6}
\qquad\qquad\ (1) \quad\quad\quad\quad&\qquad  \theo(-1) &\ \rightarrow\ & 
 0 &\ \rightarrow\ & 
0 \quad\quad\quad\quad\qquad & \phantom{\text{which is exact at } \theo^2} \\
(2) \quad\quad\quad\quad&\qquad  \theo(-1) &\ \rightarrow\ & 
 \theo &\ \rightarrow\ & 0 & \\
(7) \quad\quad\quad\quad&\qquad   \theo(-1) &\ \rightarrow\ &
  \theo^2 & \ \rightarrow\ & 0 & 
\end{alignat*}
\end{enumerate}
\end{lemma}

\begin{remark}\label{kronecker-remark}
Our numbering of these complexes is chosen to be consistent with that in 
\cite{Drezet-Le Potier}.  We will refer to a complex
of one of the forms above as a {\em standard 
complex of type (3)},\label{standard-defn}
 etc. Note that the exactness hypothesis of a complex of type (3) implies that 
the first map in that complex is automatically an injection.
\end{remark}

\begin{proof} (a) Assume that 
    $B|_E$ is not surjective. Then there exists a 
projection $\pi:\theo(1)\otimes 
V_1\rightarrow\theo(1)$ such that   $(\pi\circ B)|_E$ 
is not surjective.  
Since $\theo(1)$ is globally generated and $h^0(\theo(1)) = 3$, 
the image $W= (\pi\circ B)( \HH^0(\theo\otimes V_0))$ 
of $\HH^0(\theo\otimes V_0) $ satisfies $\dim W\leq 2$. 
This induces the commutative diagram of complexes
\begin{equation}\label{minimal res quotient}
\xymatrix{\theo\otimes V_0 \ar[r]^{B}\ar[d]^{\rho_1} & \theo(1)\otimes 
V_1\ar[d]^{\rho}\ar[r]^{\tau} & P_p(1)\ar[d]^{=}\\
\theo\otimes W\ar[r]^{\overline{B}} & \theo(1)\ar[r]^{\tau'} & P_p(1)}
\end{equation}
where $P_p$ is the  module in $\qgr S/Sg\subset \qgr S$ 
corresponding to some closed point 
$p\in E$. Moreover,
the maps $\tau$ and $\tau'$ are surjections.

If
$\dim W$ is $0$ or $1$ we obtain a complex of type 
(5), respectively (4). 

So, assume that 
  $\dim W = 2$; 
thus $\overline{B}$ induces an 
  injection of global sections. We now want to use 
a result from \cite{ATV2} which requires an algebraically closed field. 
  Thus, let $k$ have algebraic closure $F$ and use ``superscript $F$'' 
  to denote 
$-\otimes_kF$. Let $P_q^F$ denote a simple factor module of 
$P_p^F$ in $\qgr S^F/S^Fg\subset \qgr S^F$.
(Although we do not need it, this module does correspond to 
a closed point $q\in E^F$.)
 So the final row of the commutative diagram gives 
a complex 
\begin{equation}\label{alg-closed}
\theo\otimes W^F \xrightarrow{\overline{B}^F}  \theo^F(1)
\xrightarrow{\tau''} P_q^F(1)\rightarrow 0,
\end{equation}
where $\overline{B}^F$ still gives an injection of global sections 
and $\tau''$ is  surjective.
 By \cite[Proposition~6.7(iii)]{ATV2} the minimal resolution of
  $P_q^F(1)$ has the form 
  \begin{equation}\label{alg-closed2}
0\rightarrow \theo^F(-\epsilon-1) \xrightarrow{\alpha}
    \theo^F\oplus \theo^F(-\epsilon) \xrightarrow{\beta} \theo^F(1)
   \xrightarrow{\tau''} P_q^F(1)\rightarrow 0,
   \end{equation}
   for some $\epsilon \geq 0$. Since $\overline{B}^F$ is 
   injective on global sections, $\epsilon=0$. 
   We may therefore assume that $\beta=\overline{B}^F$, in which case 
 \eqref{alg-closed2}
  is an 
   extension of \eqref{alg-closed}. Since $\tau''$ factors through
   $(\tau')^F$, this forces $P_p^F=P_q^F$ to be simple and \eqref{alg-closed}
   to be exact.
   
Returning to $k$-algebras, this implies that the second row of 
\eqref{minimal res quotient} is also exact. 
Applying $\Gamma^*$ to that sequence gives the complex
\begin{equation}\label{alg-closed4}
S\otimes_kW \xrightarrow{\Gamma^*(\overline{B})} S(1) 
\xrightarrow{\Gamma^*(\tau')}\Gamma(P_p(1)).
\end{equation}
By Lemma~\ref{chi-lemma}(1), 
   $\Gamma^*(\tau')$ is surjective in high degrees and it
    cannot be zero in any degree $\geq 1$. 
   Since $S$ has Hilbert series $(1-t)^{-3}$, computing the dimension of 
  \eqref{alg-closed4}  shows that 
$\Gamma^*(\overline{B})$ is not surjective in  degree one. 
Equivalently, the map $\alpha$ of \eqref{alg-closed2}  is defined
over $k$.

Summing up, this implies that
   \eqref{minimal res quotient} can be extended  to the following
commutative diagram for which the second row is exact
and the columns are surjections:
\begin{equation}\label{alg-closed3}
\xymatrix{\theo(-1)\otimes V_{-1}\ar[r]^A&\theo\otimes V_0 
\ar[r]^{B}\ar[d]^{\rho_1} & \theo(1)\otimes V_1\ar[d]^{\rho}\\
\theo(-1)\ar[r]^{\overline{A}}&\theo^2\ar[r]^{\overline{B}} & \theo(1).}
\end{equation}
We may extend the vertical maps to a map
 of Kronecker complexes, thereby showing that 
$\bK$ has a factor   of type (3) or type (6).

(b) This follows from part 
(a) applied to the dual complex. In more detail, first
suppose that  $\bK$ is a Kronecker complex for which
$A_{\bK}|_E$ is not injective.  Then $\bK^*$ is a Kronecker
complex (of left $S$-modules!) in which 
$B_{\bK^*}|_E = (A_{\bK}|_E)^*$ is not surjective.  
Thus, we may apply part
(a) to obtain a quotient complex of $\bK^*$ of one of the given types;
dualizing this gives the appropriate subcomplex of $\bK$.

On the other hand, suppose that $A_{\bK}|_E$ is injective but its cokernel
$C$ has a subobject $T$ of finite length. 
As $E$ is a Gorenstein curve, 
$\on{\underline{Ext}}^2_E(C/T,\theo_E)=0$
but $\on{\underline{Ext}}^1_E(T,\theo_E)\neq 0$.
Thus, by the long exact sequence in cohomology, 
$\on{\underline{Ext}}^1_E(C,\theo_E)\not=0$.
Consider the  exact sequence 
$0\to P\to Q\to C\to 0$, where 
$P=\theo_E(-1)\otimes V_{-1}$ and 
$Q=\theo_E\otimes V_0$. Dualizing 
gives the exact sequence 
$$0\longrightarrow C^*\longrightarrow Q^*
\ {\buildrel{\theta}\over{\longrightarrow}}\
P^*\longrightarrow \on{\underline{Ext}}^1_E(C,\theo_E)\longrightarrow 0.$$
Thus, $A_{\bK^*}|_E=\theta$ is not surjective. 
As in the last paragraph, we can now apply part (a)
 to the dual complex $\bK^*$ to 
find the  appropriate subcomplex of $\bK$.
\end{proof}

\begin{remark}\label{big table}
One has the following table of numerics for the various standard complexes.  
Let $r=\rank$, $c_1$, $\chi$ denote the invariants
of the complex $\bK$ and use the same letters 
with primes attached to  denote the corresponding invariants of the standard 
complex $\bK'$.
\bd
\begin{array}{|l| r|r|r| c| c| c|}\hline
\text{type} & r' & c_1' & \chi' & rc_1' - r'c_1 & r\chi' - r'\chi
& r p_{\bK'} - r'p_{\bK}\\ \hline
(3) & 0 & 0 & -1 & 0 & -r & -r\\
(4) & 0 & -1 & -2 & -r & -2r & -r(t+2)\\
(5) & -1 & -1 & -3 & -r+c_1 & -3r+\chi & (-r+c_1)t + (-3r+\chi)\\
(6) & 1 & -1 & -1 & -r-c_1 & -r-\chi & (-r-c_1)t + (-r-\chi)\\
(1) & -1 & 1 & 0 & r+c_1 & \chi & (r+c_1)t + \chi\\
(2) & 0 & 1 & 1 & r & r & rt+ r\\
(7) & 1 & 1 & 2 & r-c_1 & 2r-\chi & (r-c_1)t +(2r-\chi) \\ \hline
\end{array}
\ed
\end{remark}

\noindent
An immediate consequence of this table is:

\begin{corollary}\label{specials} Let $S\in \ASthree'$.
A semistable, normalized Kronecker complex $\bK$ has no 
quotient complex of types (3--6) and no subcomplex of types (1), (2) or (7).
\qed
\end{corollary}

\begin{lemma}\label{filtration by type (3)} Let $S\in \ASthree'$ and 
suppose that $\bK$ is a normalized semistable Kronecker complex of the
form \eqref{typical Kronecker complex}.  Then $\bK$ 
admits a filtration by subcomplexes 
\bd
0=F_0\subset F_1\subset F_2 \subset\dots\subset F_\ell\subseteq \bK
\ed
such that 
\begin{enumerate}
\item[(i)]\label{filtration by type (3) subcomplexes 1} 
each $F_i/F_{i-1}$ is of type (3), 
\item[(ii)]\label{filtration by type (3) subcomplexes 2} 
each $\bK/F_i$ contains no subcomplex of type (1), (2) or (7), and
\item[(iii)]\label{filtration by type (3) subcomplexes 3} 
the first map $A_{\bK/F_\ell}|_E$ in $(\bK/F_\ell)|_E$ is injective
and $\on{coker}(A_{\bK/F_\ell}|_E)$ has no simple subobjects in $\coh E$.
\end{enumerate}
\end{lemma}
\begin{proof}  If
 $A_{\bK}|_E$ is injective and $\on{coker}(A_\bK|_E)$ is torsion-free
then, by Corollary~\ref{specials}, 
 the lemma holds with $\ell =0$. If $A_{\bK}|_E$ does not have both these
properties, then
  Lemma~\ref{enumeration of Kronecker complexes} and 
 Corollary~\ref{specials} imply that $\bK$ has a subcomplex 
 $F_1$ of type (3).

Suppose that we have constructed $0=F_0\subset\cdots \subset F_k$
so that each  subquotient $F_i/F_{i-1}$  is of type (3). If $A_{\bK/F_k}$ 
satisfies condition (iii) above then we are finished, so suppose not.  Then
Lemma~\ref{enumeration of Kronecker complexes}
implies that $\bK/F_k$ has a subcomplex $\bK'/F_k$ of type (1), (2), (3) or
 (7). Additivity of rank and $c_1$ implies that
\bd
\on{rk}(\bK)c_1(\bK') - \on{rk}(\bK')c_1(\bK) = 
\on{rk}(\bK/F_k)c_1(\bK'/F_k) - \on{rk}(\bK'/F_k)c_1(\bK/F_k).
\ed
By Remark~\ref{big table}, this is strictly positive if $\bK'/F_k$ is
 of type (1), (2) or (7)
contradicting the semistability of $\bK$.  Thus
 $\bK/F_k$  must contain a complex $F_{k+1}/F_k$ of 
type (3). Now apply  induction.
\end{proof}

There is one possibility that does not occur 
in \cite{Drezet-Le Potier} and is more difficult to treat.  This is when $\bK$ 
is a complex of the form \eqref{typical Kronecker complex}
for which $\on{coker}(B)$ has a composition series of fat points. Here, 
a \emph{fat point in $\qgr S$}\label{fat-defn} is a simple 
object $\cM$ whose module of global sections $\Gamma^*(\cM)$
has Hilbert series $n(1-t)^{-1}$  for some $n>1$. 
Fat points do not belong to $\qgr S/gS$. The next few lemmas are
concerned with this case.

\begin{lemma}\label{general-lifting} Let $S\in \ASthree'$ and 
assume that we have an exact sequence 
\begin{equation}\label{general-lifting1}
\theo^a_{S}\ {\buildrel{\alpha}\over{\longrightarrow}}\
\theo^b_{S}(1)\ {\buildrel{\beta}\over{\longrightarrow}}\
{\mathcal K}\longrightarrow 0.\end{equation}
for which $\cK\not=0$ but $\beta|_{E} = 0$. 
\begin{enumerate}
\item If $K=\Gamma^*({\mathcal K})$ then  
GKdim$(K)= 1$ and $K_ng=K_{n+3}$, for all $n\gg 0$.

\item Let $\beta(t) : \theo^b(t+1)\to {\mathcal K}(t)$ be the
 natural map induced by $\beta$ for $t\in {\mathbb Z}$. Then 
$\dim_k {\mathrm{Im}}\, \HH^0(\beta(t)) \geq 
\dim_k {\mathrm{Im}}\, \HH^0(\beta(-1))$
 for all $t\geq 0.$
\end{enumerate}
\end{lemma}

\begin{proof} (1)  By right exactness, ${\mathcal K}|_E=0$. Thus if
$K=\Gamma^*({\mathcal K})$, then $K/Kg $ is finite dimensional 
and so $K_n=K_{n-3}g$ for all $n\gg 0$. 

(2)
 Break \eqref{general-lifting1} into two exact sequences:
\begin{equation}\label{general-lifting11}
0\longrightarrow {\mathcal A} \longrightarrow\theo^a_{S}\ 
{\buildrel{\alpha}\over{\longrightarrow}}\
{\mathcal B} \longrightarrow 0.\end{equation}
and 
\begin{equation}\label{general-lifting12}
0\longrightarrow{\mathcal B}\ {\buildrel{\iota}\over{\longrightarrow}}\
\theo^b_{S}(1)\ {\buildrel{\beta}\over{\longrightarrow}}\
{\mathcal K}\longrightarrow 0.\end{equation}
 Taking cohomology of \eqref{general-lifting12} gives the exact sequence
$$ \HH^0(\theo^b(1+u)) \xrightarrow{\HH^0(\beta(u))} 
\HH^0({\mathcal K}(u)) \to \HH^1({\mathcal B}(u)) \to 0
\qquad\mathrm{for}\quad u\in \mathbb Z.$$
 By
 \cite[Proposition~6.6(iv)]{ATV2} and Lemma~\ref{chi-lemma}(3),
 $h^0({\mathcal K}(u))=\dim \HH^0({\mathcal K}(u))$ 
is constant for $u\in \mathbb Z$. 
Thus, by  part (1), in order to prove the lemma  it suffices to show that 
$h^1({\mathcal B}(-1))\geq h^1({\mathcal B}(t))$, for $t\geq 0$.

By  Lemma~\ref{chi-lemma}(5),  $\HH^j(\theo^a(u))=0$ for $j=1,2$ and
$u\geq-1$.  Thus, taking cohomology of shifts of 
\eqref{general-lifting11} shows that
$\HH^1({\mathcal B}(u)) \cong  \HH^2({\mathcal A}(u))$ for
$u\geq-1$. So it suffices to prove that  $h^2({\mathcal A}(-1))\geq
h^2({\mathcal A}(t))$ for $t\geq 0$.
Finally,  Serre duality (Proposition~\ref{serre-duality}) gives 
\bd
h^2(\cA(t)) = \dim_k\Hom(\cA,\theo(-3-t))\,\,\ \leq \,\, 
\dim_k\Hom(\cA,\theo(-2)) =
 h^2(\cA(-1)),
\ed
for $t\geq 0$.
 \end{proof}

 An $S$-module $M$ is called {\it $1$-homogeneous}\label{homog-defn}
  if $\GKdim M = 1$ and 
 $M$ has no finite-dimensional submodules. 
 
\begin{lemma}\label{q-resolution}
Let $Q$ be a graded  $1$-homogeneous right  $S$-module that is
 generated by a $t$-dimensional 
subspace of $Q_1$. Then the minimal projective resolution of $Q$ has the form
\begin{equation}\label{q-resolution1}
0 \longrightarrow P_2\ {\buildrel{\theta}\over\longrightarrow}\
P_1 \ {\buildrel{\phi}\over\longrightarrow }\ S(1)^t \longrightarrow
  Q\longrightarrow 0,
\end{equation}
such that the following holds:
\begin{enumerate}
\item[(1)] $P_2\cong S^s(-1) \oplus P_2'$
where $P_2'\cong \bigoplus_m S(j_m)$ and each $j_m<-1$;
\item[(2)] $P_1\cong S^r \oplus P_1'$
where $P_1'\cong \bigoplus_m S(k_m)$ and each $k_m<0$;
\item[(3)] $r\leq 2t$ and $\theta(S^s(-1))\subseteq 
S^r$.
\end{enumerate}
\end{lemma}

\begin{proof}
 \cite[Proposition~2.46]{ATV2}
implies that $\on{hd}(Q)\leq 2$ and so $Q$
 has a minimal resolution of the form \eqref{q-resolution1}
satisfying (1) and (2).
By Lemma~\ref{general-lifting}, $\dim Q_2\geq t$ 
and so computing the dimension of \eqref{q-resolution1} in degree $2$ shows
that $r\leq 3t-\dim Q_2\leq 2t$. It follows easily from the minimality of 
\eqref{q-resolution1} that $\theta(S^s(-1))\subseteq  S^r$.
\end{proof}

\begin{prop}\label{kill-q}
 Suppose that
\begin{equation}\label{kill-q1}
\bK: \;\;
\theo(-1)^a
\xrightarrow{\theta}
\theo^b \xrightarrow{\phi}
\theo(1)^c
\end{equation}
is a Kronecker complex with $\theta$ injective and
 $\GKdim(\Gamma^*{\mathrm{coker}}(\phi))\leq 1$.  Then either 
$\on{coker}(\phi) = 0$ or there exists a nonzero quotient 
complex of $\bK$
of the form 
 \begin{equation}\label{kill-q2}
\bK':\;\; \theo(-1)^s \xrightarrow{\theta'} \theo^r
 \xrightarrow{\phi'} \theo(1)^t 
 \end{equation}
with $s+t\leq r\leq 2t$.
\end{prop}

\begin{proof} Assume that $\on{coker}(\phi)\neq 0$ and  pick a nonzero, 
$1$-homogeneous
factor  $Q$ of $\Gamma^*\left({\mathrm{coker}}(\phi)\right)$. Then
$Q$ is generated by a subspace of
$Q_1$ of dimension $t\leq c$. By Lemma~\ref{q-resolution} there 
is a resolution of the form \eqref{q-resolution1}, the 
image of which in $\qgr S$ has the form
\bd
{\mathbf R}: \;\; 0\rightarrow \theo(-1)^{s'} \oplus \pi P_2' 
\xrightarrow{\theta'}  \theo^{r'} \oplus \pi P_1' \xrightarrow{\phi'} 
\theo(1)^t\xrightarrow{\psi'} \pi Q\to 0
\ed
 where the $P_i'$ are generated in   degrees $\geq i$.
 By construction, the map 
$S(1)^c\twoheadrightarrow Q$  factors through $S(1)^t$.  
We next show that this induces a map $\bK\to {\mathbf R}$; more precisely, 
we construct a commutative diagram of the form
 \bd
  \begin{CD}
 0@>>> S(-1)^a @>\theta >> 
 S^b @>\phi>> S(1)^c
 @>\psi >> Q @>>> 0  \\
 @. @V\alpha_2 VV @V \alpha_1 VV @V \alpha_0 VV @V\alpha_{-1} VV \\
 0@>>>  S(-1)^{s'}\oplus P_2'@>\theta' >> 
S^{r'}\oplus P_1' @>\phi'>> S(1)^t
 @>\psi' >> Q @>>> 0
 \end{CD}
 \ed
whose image in $\qgr S$ is the desired  map $\bK\to {\mathbf R}$.
 The maps $\alpha_0$ and $\alpha_{-1}$ are the surjections already defined
 and so the final square is commutative.
Since the second row of this
 diagram is exact, the usual diagram chase constructs  $\alpha_1$ and
 $\alpha_2$.

By the definition of the $P_i'$, one has $\alpha_1(S^b)\subseteq S^{r'}$ and
$\alpha_2(S(-1)^a)\subseteq S(-1)^{s'}$. Thus, if  
  $\bK'$ is  the image of $\bK$ in ${\mathbf R}$, then $\bK'$ is a Kronecker
  complex of the the form 
  $$ 
  \theo(-1)^s \xrightarrow{\theta'}  \theo^r \xrightarrow{\phi'}  \theo(1)^t
 ,$$
  where $r = \on{rk}(\alpha_1(\theo^b))$ and 
  $s=\on{rk}(\alpha_1(\theo(-1)^a))$ .  The   map $\theta'$ is injective
  since it is induced from the injection $\theta$.  
Moreover, $r \leq r'\leq 2t$
by Lemma~\ref{q-resolution}(3).

It remains to prove that $s+t\leq r$.
Since this is just a question of ranks, we may pass to the (graded or ungraded)
division ring of fractions $D$ of $S$.  Then $\bK'\otimes_S D$ gives a complex
$D^s \hookrightarrow D^r \xrightarrow{\phi'\otimes D} D^t$.
By hypothesis, $\phi\otimes D$ is surjective and hence so is 
 $\phi'\otimes D$. This gives the desired inequality.
\end{proof}

We are now ready to prove Proposition~\ref{semistable Kronecker complex 
is monad}.

\begin{proof}[Proof of Proposition~\ref{semistable Kronecker complex 
is monad}]  
Recall that we are given a semistable, normalized Kronecker complex of the form 
$$\bK:\qquad \theo(-1)\otimes V_{-1} \xrightarrow{A} 
\theo\otimes V_0 \xrightarrow{B} \theo(1)\otimes V_1,$$
over an algebra $S\in \ASthree$ and we wish to prove that $\bK$ is a 
torsion-free monad.
As was remarked earlier, since the proposition is a result about $\qgr S$,
we may assume that $S=S(E,\cL,\sigma)\in \ASthree'$.

Choose a filtration $F_0\subset \cdots \subset 
F_\ell\subset \bK$ by Lemma~\ref{filtration by type (3)}.
By Remark~\ref{kronecker-remark},
$\HH^{-1}(F_i/F_{i-1}) = 0$ for each $i$ and so $\HH^{-1}(F_\ell) = 0$.
Next, consider $\ker(A_{\bK/F_\ell})$.  We have an exact sequence
\bd
0\rightarrow \ker(A_{\bK/F_\ell})\rightarrow \theo(-1)\otimes V'_{-1}
\xrightarrow{A_{\bK/F_\ell}} \im(A_{\bK/F_\ell})\rightarrow 0,
\ed
and $\im(A_{\bK/F_\ell})$ is torsion-free since it is a submodule of a vector
bundle.  By Lemmas~\ref{restriction of torsion-free} and 
\ref{filtration by type (3)}(iii), this implies
that $\ker(A_{\bK/F_\ell})|_E = \ker(A_{\bK/F_\ell}|_E) = 0$.
  Thus $\ker(A_{\bK/F_\ell})$ is a torsion
   submodule of $\theo(-1)\otimes V_{-1}'$ and so
it is zero. Therefore, $A_{\bK/F_\ell}$ and $A_{F_\ell}$ are both injective, 
and so $A_{\bK}$ is also injective.

We next show that $B$ is surjective. By Corollary~\ref{specials}
and Lemma~\ref{enumeration of Kronecker complexes}(a),
$B|_E$ is surjective. 
Suppose that $B_\bK$ is not surjective.  Then Lemma~\ref{general-lifting}(1)
implies that $\GKdim\Gamma^*(\mathrm{Coker}(B))=1$ and so,
 by Proposition~\ref{kill-q},
there is a quotient complex $\bK'$ of
$\bK$ of the form $\theo(-1)^{t-a-b} \rightarrow \theo^{2t-a} \rightarrow 
\theo(1)^t$ with $a, b\geq 0$ and $t>0$.  By definition,
$c_1(\bK') = t-a -b -t = -a-b$ and $\on{rk}(\bK') = b$.  This gives
$$\on{rk}(\bK)c_1(\bK') - \on{rk}(\bK')c_1(\bK)
= -a\on{rk}(\bK) - b\left\{\on{rk}(\bK)+c_1(\bK)\right\}=x,
$$
say. Since $\bK$ is normalized, $x\leq 0$
 with equality if and only if $a=b=0$. So $\bK'$ de-semistabilizes $\bK$ 
unless $a=b=0$, in which case we find that
$$\on{rk}(\bK)\chi(\bK') - \on{rk}(\bK')\chi(\bK) = \on{rk}(\bK')(t-a-3t)<0$$
and so  $\bK'$ still de-semistabilizes $\bK$.
This contradiction shows that $B$ is surjective.
 
 Finally, we must show that $\HH^0(\bK)$ is torsion-free.
 Write  $\bK/F_\ell$ as the complex 
 \begin{equation}\label{cokernel1}
 0 \to \theo(-1)\otimes V_{-1}' \xrightarrow{A'} \theo\otimes V_0'
 \xrightarrow{B'} \theo(1)\otimes V_1'\to 0.
 \end{equation}
 We first want to show that $\on{coker}(A')$ is torsion-free, so assume that
its torsion submodule $\cT$ is nonzero.
 By the earlier part of the proof,  $A'$ is injective. 
Thus, by Lemma~\ref{chi-lemma}(5), applying  $\Gamma(-)$ to
 \eqref{cokernel1} gives the exact sequence 
\begin{equation}\label{cokernel2}
0\to S(-1)^a \xrightarrow{A''}  S^b \xrightarrow{B''}  C \to 0,
\end{equation}
where $A'' = \Gamma^*(A')$ and $C = \on{coker}(A'')$.
 Let $T$ denote the torsion submodule of $C$.
 By Lemma~\ref{chi-lemma}(2), $\pi(C)=\on{coker}(A')$
 and $\pi(T)=\cT$, so $T\not=0$. 
 The sequence \eqref{cokernel2} shows that
 $C$ has projective dimension $\on{hd}\,C\leq 1$.
 Since $C/T$ has no finite dimensional submodules, 
 \cite[Proposition~2.46(i)]{ATV2}
 implies that $\on{hd}(C/T)\leq 2$ and hence that $\on{hd}(T)\leq 1$.
 By \cite[Theorem~4.1(iii)]{ATV2}, $\GKdim T = 2$.
 
 By Lemma~\ref{filtration by type (3)}(iii), $A'|_E$ is injective and hence so
 is $A'' \otimes_S(S/gS)$. Thus $\underline{\on{Tor}}_1^S(C,S/gS)=0$. 
 But this Tor group is also $\on{Ker}\left[C\otimes Sg\to C\right]$. 
 Therefore multiplication by $g$ is injective on $C$ and its submodule $T$.
By \cite[Lemma~2.2(3)]{Zh2}, this implies that 
 $\GKdim T/Tg=\GKdim T-1=1$. On the other hand, as 
 $C/T$ is torsion-free, the argument in Lemma~\ref{restriction of torsion-free}
 shows that $\underline{\on{Tor}}_1^S(C/T,S/gS)=0$. 
 Therefore,  $T/Tg\hookrightarrow C/Cg$ and so  $\cT/\cT g=\pi(T/Tg)$ is 
 submodule of $\pi(C/Cg)=
 \on{coker}(A')|_E$ of finite length. 
 This contradicts Lemma~\ref{filtration by type (3)}(iii).
 
   Thus $\on{coker}(A')$ is torsion-free.
By the definition of a type (3) complex in Lemma~\ref{enumeration of Kronecker
complexes}(a), 
 $\on{coker}(A_{F_i/F_{i-1}})$ is torsion-free for each $i$.
Thus
$\on{coker}(A_{\bK})$ is an   iterated extension of torsion-free modules 
and hence is  torsion-free.  Therefore
its submodule $\HH^0(\bK)$ is torsion-free.
  This completes the proof of
  Proposition~\ref{semistable Kronecker complex is monad}.
\end{proof}

It is now easy  to generalize 
 Proposition~\ref{semistable Kronecker complex is monad} to families of
 Kronecker complexes, defined as follows.

\begin{defn}\label{typical Kronecker-defn} Let $R$ be a commutative 
noetherian $k$-algebra and $S\in \ASthree$. Recall that $S_R=S\otimes_kR$,
regarded as a graded $R$-algebra.
A {\em family of Kronecker complexes} parametrized
by a scheme $\spec R$, also called a Kronecker complex in 
$\qgr S_R$,  is a complex of the form
\begin{equation}\label{typical Kronecker family}
{\mathbf K}: \;\;\; \theo(-1)\otimes_R V_{-1}\xrightarrow{A} 
\theo\otimes_R V_0 \xrightarrow{B} \theo(1)\otimes_R V_1
\end{equation}
in $\qgr S_R$
where $V_{-1}$, $V_0$ and $V_1$ are finite-rank vector bundles on $\spec R$. 
The definition of families of (geometrically, semi)stable Kronecker complexes
follows from Definition~\ref{families def}. 
As before,  a family of geometrically (semi)stable Kronecker complexes
is automatically (semi)stable.

As usual, morphisms of Kronecker complexes in $\qgr S_R$ are morphisms of
complexes.  Two such  complexes (or modules) $\bK$ and $\bK'$
     are {\em equivalent}\label{equiv-def}
      if there is a line bundle $M$ on $\spec R$ such 
     that $\bK\otimes_RM$ and $\bK'$ are isomorphic in $\qgr S_R$.
\end{defn}

\begin{corollary}\label{monad families}
Suppose that $\bK$ is a  semistable Kronecker complex in $\qgr S_R$.
Then $\bK$ is a torsion-free monad in $\qgr S_R$ 
 in the sense of Definition~\ref{monad-defn}.
\end{corollary}
\begin{proof}
By definition, if $p\in \spec R$, then $\bK\otimes k(p)$ is semistable and so 
Proposition~\ref{semistable Kronecker complex is monad} implies that 
$\bK\otimes k(p)$ is a torsion-free monad.
\end{proof}

We end the section by checking that the notions of
semistability for monads and 
for their cohomology do correspond.

\begin{prop}\label{semistable module from complex} Let $S\in \ASthree$
 and let $R$ be a commutative 
noetherian $k$-algebra. Then the map $\bK\mapsto \on{H}^0(\bK)$ 
induces an equivalence of categories between:
\begin{enumerate}
\item[(i)] the category of  geometrically (semi)stable
 normalized Kronecker complexes $\bK$  in $\qgr S_R$ and
\item[(ii)] the category of 
$R$-flat families $\cM$ of  geometrically (semi)stable, normalized
torsion-free modules in $\qgr S_R$.
\end{enumerate}
\end{prop}

\begin{proof}  By  Corollary~\ref{monad families}, the Kronecker complex $\bK$
in   (i) is a torsion-free monad in $\on{Monad}(S_R)$, while
Lemma~\ref{vanishing cohomology for ss modules}  implies that $\cM\in (\qgr
S_R)_{\text{VC}}$ in part (ii). Now consider the equivalence of categories
$\on{Monad}(S_R)\simeq (\qgr S_R)_{\text{VC}}$ from
Theorem~\ref{beilinson5}(1). By definition 
 the torsion-free objects correspond
under this equivalence and it is elementary that the normalized objects also
correspond.  Thus in order to prove the proposition, it remains to
check that the  geometrically (semi)stable objects correspond. If $\cM\in (\qgr
S_R)_{\text{VC}}$, then Theorem~\ref{beilinson5}(2) implies that
$\bK(\cM)\otimes F = \bK(\cM\otimes F)$ for every geometric point $\on{Spec}
F\to \on{Spec} R$ and so we need only prove that the (semi)stability of
$\cM\otimes F$ is equivalent to that of $\bK(\cM\otimes F)$. In other words, 
after replacing $k$ by $F$,  it remains to prove
that  $\cM\in \qgr S_k$ is (semi)stable if and only if $\bK=\bK(\cM)$ 
is (semi)stable.

If $\bK$ is (semi)stable then  the proof of \cite[Proposition~2.3(3)]{Drezet-Le
Potier} can be used, essentially without change, to prove that
$\cM=\on{H}^0(\bK)$  is (semi)stable.

The other direction does not quite follow from the corresponding result in 
 \cite{Drezet-Le Potier} and before proving it we need a definition and a lemma.
Let $\bK$ denote a Kronecker complex over $S\in \ASthree$.
  A {\em maximal subcomplex}\label{max-subcomplex} of $\bK$ is 
a subcomplex $\bK'$ that realizes the maximum of 
$\on{rk}(\bK)p_{\bK'} - \on{rk}(\bK')p_{\bK}$ among all subcomplexes of
 $\bK$.

\begin{sublemma}\label{maximal subcomplex has surjective B}
{\rm (}\cite[Lemme~2.4]{Drezet-Le Potier}{\rm )}
Suppose that $\bK$ is a normalized monad over $S\in \ASthree'$.
 If $\bK'$ denotes a maximal subcomplex  of $\bK$, then $\HH^1(\bK'|_E) = 0$.
\end{sublemma}

\begin{proof} 
By Lemma~\ref{enumeration of Kronecker complexes},
it suffices to check that $\bK'$ has no quotient complex ${\mathbf L}$
of types  (3--6). Suppose that such an ${\mathbf L}$ exists and set 
  ${\mathbf L}' = \mathrm{Ker}(\bK'\rightarrow {\mathbf L})$. Then
\bd
\on{rk}(\bK)p_{{\mathbf L}'} - \on{rk}({\mathbf L}')p_{\bK} = 
\left\{\on{rk}(\bK)p_{\bK'} - \on{rk}(\bK')p_{\bK}\right\} 
-\left\{\on{rk}(\bK)p_{{\mathbf L}} - \on{rk}({\mathbf L})p_{\bK}\right\}
\ed
by the additivity of ranks and Hilbert polynomials.  
Remark~\ref{big table} shows that 
$\on{rk}(\bK)p_{{\mathbf L}} - \on{rk}({\mathbf L})p_{\bK}<0$.
This  contradicts the maximality of $\bK'$.
\end{proof}

We now return to the proof of 
Proposition~\ref{semistable module from
complex}. 
By the first paragraph of the proof it remains to show that,
 if $\cM\in \qgr S$ is (semi)stable,
then so is $\bK=\bK(\cM)$.

As usual, we may assume that $S=S(E,\cL,\sigma)\in \ASthree'$ and we
write $\bK$ as
$$\theo(-1)\otimes V_{-1}\xrightarrow{A} \theo\otimes V_0 
\xrightarrow{B} \theo(1)\otimes V_1.$$
If $\bK'$ is
a maximal subcomplex of $\bK$ then
Sublemma~\ref{maximal subcomplex has surjective B} implies that 
 $\HH^1(\bK'|_E)= 0$.  
Let $\bK'' =\bK/\bK'$.  We first prove that $\HH^{-1}(\bK'')= 0$.
There is a filtration
\bd
\bK' = F_0\subset F_1 \subset \dots\subset F_l\subseteq \bK
\ed
of $\bK$ with each $F_i/F_{i-1}$ of type (3), such that $\bK/F_l$ 
contains no subcomplex of type (3).  In particular, 
$c_1(F_l) = c_1(\bK')$, by Remark~\ref{big table}.  

We claim that 
 $\bK/F_l$ contains no subcomplex 
of type (1), (2), or (7). Suppose there were such a 
subcomplex, say ${\mathbf L}/F_l$.  Then Remark~\ref{big table} implies that
 $c_1({\mathbf L}) = c_1(F_l)+1=c_1(\bK')+1$. Thus
\bd
\on{rk}(\bK)c_1({\mathbf L}) - \on{rk}({\mathbf L})c_1(\bK) 
= \on{rk}(\bK)c_1(\bK') -\on{rk}(\bK')c_1(\bK) + 
\left[\on{rk}(\bK)- c_1(\bK)\on{rk}({\mathbf L}/F_l)\right].
\ed
In the three cases, Remark~\ref{big table} implies that 
 $\on{rk}({\mathbf L}/F_l) = -1, 0, 1$, respectively. This in turn implies that 
 $\left[\on{rk}(\bK)- c_1(\bK)\on{rk}({\mathbf L}/F_l)\right]>0$ (when 
$\on{rk}({\mathbf L}/F_l) = -1$ this needs  the fact 
that $\bK$ is normalized). By the displayed equation, this  
contradicts the maximality of $\bK'$ and  proves the claim.

 By Lemma~\ref{enumeration of Kronecker complexes},
the last paragraph implies that $\HH^{-1}((\bK/F_l)|_E)=0$.
Applying  Lemma~\ref{restriction of torsion-free} to the short exact sequence
\bd
0\rightarrow \HH^{-1}(\bK/F_l) \rightarrow (\bK/F_l)_{-1}  
\xrightarrow{A_{\bK/F_l}} \on{Im}(A_{\bK/F_l})\rightarrow 0
\ed
shows
that $\HH^{-1}(\bK/F_l) =0$.  But $\HH^{-1}(F_l/\bK')=0$ by construction, so 
$\HH^{-1}(\bK'') = \HH^{-1}(\bK/\bK') =0$, as desired.

It follows that the map $\HH^0(\bK')\rightarrow \HH^0(\bK)=M$ is injective.  
Let  $C=\on{coker}(B_{\bK'})$. 
As $A_{\bK'}$ is injective, 
$p_{\HH^0(\bK')} =  p_{\bK'}+p_C\geq 
p_{\bK'}$. By Sublemma~\ref{maximal subcomplex has surjective
B}, $C$ must be torsion and hence have rank zero.
In particular, $\on{rk}(\HH^0(\bK')) = \on{rk}(\bK')$.
 Since the rank of a monad equals the rank of
its cohomology, these observations imply  that 
$$
\on{rk}(\bK)p_{\bK'} - \on{rk}(\bK')p_{\bK} 
\ \leq\ \on{rk}(M)p_{\HH^0(\bK')} - \on{rk}(\HH^0(\bK'))p_M.
$$
Therefore, the (semi)stability of $M$ implies (semi)stability 
 of $\bK$, completing the proof of Proposition~\ref{semistable module from
 complex}.
\end{proof}

\section{Moduli Spaces: Construction}\label{section5}

In this section we use a Grassmannian embedding to construct
a projective moduli space of semistable modules in $\qgr S$ for $S\in \ASthree$
 (see Theorem~\ref{projective moduli spaces}) and determine cases when the
 moduli space is fine (see Proposition~\ref{finemodulispace} and its corollaries).
 This will complete the proof of Theorems~\ref{firstthm} and
 \ref{three-to-one-intro} from the introduction.
 The more subtle properties of these moduli spaces (notably,
 determining cases
 where they are smooth or connected) will be examined in
 Section~\ref{moduli-smooth}.
 
 The idea behind the proofs of this section is similar to the classical case: 
in Subsection~\ref{semistability via Grassmannians}
we prove that   the moduli functor of
 framed semistable Kronecker
complexes (semistable Kronecker complexes that have been appropriately
rigidified) can be embedded  in a product of
Grassmannians $\GR$, and that the image is  the semistable locus
for the natural group action on a closed subscheme of $\GR$.
This allows us
to prove in Subsection~\ref{moduli-subsec} that the GIT quotient of
this subscheme is exactly the moduli space we wish to construct.
The description in terms of Kronecker complexes then gives a convenient
tool for proving smoothness and existence of universal modules.

 We remind the reader of the relevant definitions.
Let  $F$ be a contravariant functor from the category of (typically,
 noetherian
affine) $k$-schemes to the category of sets. Then a scheme $Y$
together with a morphism of functors $F\rightarrow h_Y =
\Hom_{k\on{-Sch}}(- ,Y)$  {\em corepresents}\label{corep-defn} $F$ if
this morphism is universal for morphisms $F\rightarrow h_X$ from $F$ to schemes
$X$.  In this case $Y$ is also called a  \emph{coarse moduli space for
$F$}\label{moduli-defn}. The scheme $Y$ \emph{represents} $F$, equivalently
is a \emph{fine moduli space for $F$}, if $F=h_Y$ under the usual embedding
of the category of schemes into the category of contravariant functors.

 We are interested in classifying torsion-free modules over $S\in \ASthree$ 
 and by shifting there is no harm in assuming that they are normalized (see
 Remark~\ref{other-c1} for the formal statement).
  By Proposition~\ref{semistable
 module from complex} this means that we
 can work with normalized 
Kronecker complexes over $S\in \ASthree$.
Thus, through the end of Subsection~\ref{semistability via Grassmannians} 
we can and will fix the following data.

\begin{notation}\label{invariants-notation} Fix 
integers $r\geq 1$, $c_1$ and $\chi$ with $-r<c_1\leq 0$.
A normalized
Kronecker complex with rank $r$, first Chern class $c_1$ and Euler 
characteristic $\chi$ will be called 
 \emph{a Kronecker complex with invariants $\{r,c_1,\chi\}$.}
In order to avoid excessive subscripts it will typically be written  
as 
\begin{equation}\label{Gr Kr cplx}
K\otimes\theo(-1)\xrightarrow{A} H\otimes\theo\xrightarrow{B} 
L\otimes\theo(1),
\end{equation}
where $\dim K = d_{-1}$, $\dim H=n$ and $\dim L=d_{1}$. 
These numbers are determined by the other invariants and so 
they too will be fixed throughout the section. Specifically:
\begin{equation}\label{lin-alg}
d_{-1} = 2c_1+r-\chi,\qquad 
n=3r+3c_1-2\chi \quad\text{and}\quad 
 d_{1}=c_1+r-\chi.
\end{equation}
 Finally, the vector space $H$ of dimension $n$ will also be
  fixed throughout.
 \end{notation}

Let $R$ be a commutative $k$-algebra. 
A   complex in $\qgr S_R$ is a {\em framed Kronecker
complex}\label{framed-complex}
$(\bK,\phi)$ if it is a Kronecker complex of the form
 \eqref{typical Kronecker family}
with a specific choice of isomorphism $\phi: V_0\cong H_R=R\otimes H$.
Equivalently, it has the form  
\begin{equation}\label{framed cplx}
K\otimes_R \theo_{S_R}(-1) \xrightarrow{A} H_R\otimes_R \theo_{S_R}
\xrightarrow{B} L\otimes_R \theo_{S_R}(1),
\end{equation}
for some bundles $K$ and $L$.  It is implicit in this description
that isomorphisms of framed Kronecker complexes restrict to the 
identity on $H_R\otimes\theo_{S_R}$.

In the notation of Definitions~\ref{typical
Kronecker-defn} and \ref{families def},
 let $\cK^{ss}$ denote the moduli functor for equivalence
classes of families  of geometrically
semistable Kronecker complexes in $\qgr S$ with
invariants $\{r,c_1,\chi\}$. Write $\cK^s$ for its
subfunctor of geometrically stable Kronecker complexes.
Let $\widehat{\cK}^{ss}$ denote the moduli functor of families of framed
geometrically semistable Kronecker 
complexes with invariants $\{r,c_1,\chi\}$ with 
 subfunctor $\widehat{\cK}^s$
of geometrically stable framed Kronecker complexes.
\label{Kronecker functors}
The group $\on{GL}(H)$ acts on $\widehat{\cK}^{ss}$: an element $g\in\on{GL}(H)$
takes a complex $\bK$ as in
\eqref{framed cplx} to a complex $g\cdot \bK$ 
with the same objects but with the maps $A$ and $B$ 
replaced by $(g\otimes 1)\circ A$ and $ B\circ (g^{-1}\otimes 1)$.

\subsection{Semistability via Grassmannians}
\label{semistability via Grassmannians}

In this subsection, we convert  the problem of classifying  Kronecker complexes
into a problem about Grassmannians.  We begin with a general framework,
since this will also enable us to describe our moduli spaces for families of AS
regular algebras as well as for individual algebras.

Let $\cB$ be a $k$-scheme. We define  $S= S_{\cB}$ to be  \emph{ a $\cB$-flat 
family of algebras in $\ASthree$}\label{family-as-defn}
 if $S$ is a flat, connected graded 
$\cB$-algebra such that $S_{\cB}\otimes k(p)\in \ASthree(k(p))$
for all points $p\in \cB$.  If $U=\spec R\xrightarrow{f}\cB$
is an affine $\cB$-scheme, then we write 
$S_U=S_\cB \otimes_{\theo_{\cB}}\theo_U$ for the corresponding 
family of $\theo_U$-algebras.

We will be interested in subschemes of the
following product of Grassmannians:
\begin{equation}\label{Grass-defn}
\GR_\cB = \Gr_{d_{-1}}(H_{\cB}\otimes_{\theo_\cB} S_1)\times 
\Gr^{d_{1}}( H_{\cB}\otimes_{\theo_\cB} S_1^*),
\end{equation}
where $\Gr_{d_{-1}}$
 denotes the relative 
Grassmannian of rank $d_{-1}$ $\theo_\cB$-subbundles of 
$ H_{\cB}\otimes S_1$ 
and $\Gr^{d_{1}}$ denotes the relative Grassmannian
of rank $d_{1}$ quotient bundles of $ H_{\cB}\otimes S_1^*$. As in
\cite[Example 2.2.3]{HLbook}, $\GR_\cB$ represents the
 functor that, to
a $\cB$-scheme $U\xrightarrow{f}\cB$, associates the set of pairs (called
\emph{KL-pairs})
\begin{equation}\label{KL-pair}
 \big(K\ \buildrel{i}\over{\hookrightarrow}\
 H_U\otimes_{\theo_U}(f^*S_1),\ H_U \otimes_{\theo_U}(f^*S_1)^*\
\buildrel{\pi}\over{\twoheadrightarrow}\  L\big)
\end{equation}
 of a subbundle $K$ of $H_U\otimes S_1 $
of rank $d_{-1}$ 
and a quotient bundle $L$ of $H_U\otimes S_1^*$ of rank $d_{1}$
on $U$.  In particular, by the Yoneda Lemma,
 the identity map $\GR_{\cB}\rightarrow\GR_{\cB}$
determines a $KL$-pair 
$$\cP = (\overline{K} \xrightarrow{i_{\cP}} S_1\otimes H_{\GR}, \,\,
S_1^*\otimes H_{\GR}
\xrightarrow{\pi_{\cP}}\overline{L})$$
on $\GR_{\cB}$ that is {\em universal}\label{universal-k-l} 
in the sense that a $KL$-pair on $U$ is the pullback 
$\tilde{f}^*\cP$ along a unique 
map $U\xrightarrow{\tilde{f}} \GR_{\cB}$.
We will typically suppress the pull-back maps $f^*$ or $\tilde{f}^*$
when it causes no confusion.

\begin{remark}\label{inv construction}
Let $\cH$ be a vector bundle of rank $n$ on a $k$-scheme $\cB$,
respectively $\cH = H_{\cB}$ and let
$$P = (K\subset \cH_U\otimes S_1, \cH_U\otimes S_1^*\rightarrow L)$$
be a vector subbundle and quotient bundle over $U=\spec R$. Then
one obtains a diagram
of the form \eqref{typical Kronecker family}, respectively 
\eqref{framed cplx},
by taking for $A$ the composite map
\bd 
A:\; K\otimes \theo_{S_U}(-1) \rightarrow 
 \cH_U \otimes S_1\otimes \theo_{S_U}(-1)
\rightarrow \cH_U\otimes\theo_{S_U},
\ed
(where the final map comes from multiplication in $S_U$ and all tensor products
are over $R$) and for
$B$ the composite map
\bd 
B:\; \cH_U\otimes \theo_{S_U}\rightarrow
 \cH_U \otimes S_1^* \otimes S_1\otimes\theo_{S_U}
\rightarrow  L\otimes S_1\otimes\theo_{S_U}
\rightarrow L\otimes \theo_{S_U}(1).
\ed
Let $\bK(P)$ denote the diagram associated to $P$. Since
$\bK(P)$ need not be a complex, we define
 $N_{\cB,\cH}$\label{straight-N-def}
  to be the subfunctor of the relative Grassmannian
functor $\GR_{\cB,\cH}$ consisting
of those $P$ for which the associated diagram  
$\bK(P)$ is a complex.  When $\cH = H_{\cB}$,   write 
 $N_{\cB,H_{\cB}}=N_{\cB}$;
then the  above construction defines a map from $N_{\cB}$
to the functor of isomorphism classes 
 of framed Kronecker complexes.
\end{remark}

\begin{remark}\label{equivariance} Given a vector space $W$, 
   $\on{GL}(H)$ acts on the Grassmannian 
$\on{Gr}= \on{Gr}_{d}(W\otimes H)$ as follows:
$g\in\on{GL}(H)$  takes a subspace
$V\subset W\otimes H$ to $gV = (1\otimes g)\cdot V$ and the quotient
$H/V$ to $H/(1\otimes g)V$.   
Since this action is defined via 
the subspaces themselves it gives a $\on{GL}(H)$-equivariant structure on
 the universal subbundle and universal quotient.
Note that if
$g$ is in the centre of $\on{GL}(H)$ then $gV=V$ but $g$ acts by
 scalar multiplication $m_g$ on $V$ and $H/V$.
  
More globally, let 
$\cV\xrightarrow{\iota}
\theo_{\on{Gr}}\otimes W\otimes H$ denote the 
universal subbundle on the Grassmannian. The composition
$$\widetilde{\iota}: \cV\xrightarrow{\iota} 
\theo_{\on{Gr}}\otimes W\otimes H
\xrightarrow{1\otimes g}\theo_{\on{GL}(H)\times\on{Gr}}\otimes W\otimes H$$
has image $g\cV_x$ in the fibre over $x\in\on{Gr}$ and so  
the map $m_g: \on{Gr}\rightarrow \on{Gr}$
is the one determined by the subbundle $\on{Im}(\widetilde{\iota})$.  
Conversely, if 
$\theo_{\on{Gr}}\otimes W\otimes H\xrightarrow{\pi} \cQ$ is the universal quotient
bundle (that is, the cokernel of $\iota$), then the composite map
$$\theo_{\on{Gr}}\otimes W\otimes H
\xrightarrow{1\otimes g^{-1}}\theo_{\on{Gr}}\otimes W\otimes H
\xrightarrow{\pi} \cQ$$ has the same kernel as $\widetilde{\iota}$, so it
also corresponds to $m_g:\on{Gr}\rightarrow \on{Gr}$.

Now let $\on{GL}(H)$ act diagonally on $\GR_{\cB}$.  Then the last
paragraph implies that $m_g^*\cP$ is isomorphic to the pair given
by the maps 
$\big((1\otimes m_g)\circ i_{\cP}, \pi_{\cP}\circ (1\otimes m_{g^{-1}})\big)$.
Thus $\bK(m_g^*\cP)\cong g\cdot \bK(\cP)$.  In particular, the subfunctor $N_{\cP}$ 
is $\on{GL}(H)$-invariant and the map from $N_{\cB}$ to the functor
of isomorphism classes  of framed Kronecker complexes is 
$\on{GL}(H)$-equivariant.
\end{remark}
\begin{lemma}\label{subscheme N}
Let $S=S_\cB$ be a flat
family of algebras in $\ASthree$
parametrized by a $k$-scheme $\cB$. Then:  
\begin{enumerate}
\item[(1)] There is a closed $\on{GL}(H)$-invariant subscheme 
$\cN_\cB\subset\GR_\cB$ that 
represents the subfunctor $N_{\cB}$.
\item[(2)]
If $\cB'\rightarrow \cB$ is a morphism of $k$-schemes,
 then $\cN_{\cB'} = \cN_{\cB}\times_{\cB}\cB'$.
\end{enumerate}
\end{lemma}

\begin{proof}  
Use Remark~\ref{inv construction} to 
construct the family of diagrams \eqref{framed cplx} corresponding to the
universal KL-pair $\cP$ on $\GR$.
 The composite map $B\circ A$
is zero if and only if
\bd
\HH^0\big((B\circ A)(1)\big): \overline{K} = \overline{K}\otimes
\HH^0(\theo)
\rightarrow \HH^0(\overline{L}\otimes \theo(2))=
\overline{L}\otimes S_2
\ed
is zero.  But this is the map of vector bundles
\bd
\phi: \overline{K} \xrightarrow{i}  H_{\GR}\otimes S_1\rightarrow 
H_{\GR}\otimes S_1^*\otimes S_1\otimes S_1
\rightarrow \overline{L}\otimes S_1\otimes S_1
\xrightarrow{\overline{L}\otimes \on{mult}} \overline{L}\otimes S_2.
\ed
The subfunctor of part (1) is therefore represented by
the vanishing locus $\cN_\cB$ of 
this map  of vector bundles on $\GR_\cB$ and so it is closed.
By Remark~\ref{equivariance}, $\cN_{\cB}$ is
 $\on{GL}(H)$-invariant.

Part (2) follows since the functor represented by $\cN_{\cB'}$ is the fibre
product of the functors represented by $\cN_\cB$ and $\cB'$.
\end{proof}

We now follow \cite{Le Potier} to show that framed semistable 
Kronecker complexes correspond to points of an open subscheme
of $\cN$. For this we can specialize to $\cB=\spec k$ and hence
work with a fixed 
AS regular algebra $S\in \ASthree$.
The conventions concerning $\{r,c_1,\chi\}$ and $H$
  from Notation~\ref{invariants-notation} will be maintained.
Write $\GR=\GR_S$   for the corresponding 
product $\GR_{\on{Spec} k}$
of Grassmannians associated to $S$, as in \eqref{Grass-defn},
and $\cN=\cN_S\subset\GR$\label{en-defn} for the closed subscheme 
$\cN_{\on{Spec} k}$ 
defined by  Lemma~\ref{subscheme N}.

The next theorem identifies $\widehat{\cK}^{ss}$ as a subscheme
of $\cN$, for which we need the following construction. 

\begin{construction}\label{H' construction}
Consider a KL-pair $(K,L)$ 
as in \eqref{KL-pair}, where $U=\spec F$ for some field extension $F$ of $k$.
Given a subspace $H'\subset H$, set
$K' = K\cap (H'\otimes S_1)$ and $L'= \im(H'\otimes S_1^*\rightarrow L)$.
If $(K,L)$ lies in ${\mathcal N}$, then $K'\otimes\theo(-1)
\rightarrow H'\otimes\theo\rightarrow L'\otimes\theo(1)$ 
is also a Kronecker complex, and we let $r'$, $c_1'$, $\chi'$ 
denote the rank, first Chern class and Euler characteristic
of this complex. 
\end{construction}

\begin{lemma}\label{description of Nss} (1) Let $S\in \ASthree$. Then
$\widehat{\cK}^{ss}$  is $\on{GL}(H)$-equivariantly isomorphic to
 the subfunctor $\cN'$ 
 of $\cN$ defined as follows: $\cN'(\spec R)$ consists of those KL-pairs $(K,L)$ 
such that for each geometric point 
$\spec F\rightarrow \spec R$, the KL-pair
$(K\otimes F, L\otimes F)$  satisfies the following condition.
\begin{equation}\label{dagger}
\begin{array}{rl}
 &\text{For every  }\{0\}\subsetneq H'\subsetneq H_F, 
\text{ one has } r(c_1'm+\chi') - r'(c_1m+\chi) \leq 0\\
&  \text{ as polynomials in }m.
\end{array}
\end{equation}

(2) The same result holds for $\widehat{\cK}^s$ provided one replaces 
$\;\leq 0$ by $<0$ in \eqref{dagger}.
\end{lemma} 

\begin{proof}  (1) 
Suppose that  $K$ and $L$ are vector bundles on $U=\spec R$ and that 
\bd  \bK: \ 
K\otimes\theo_{S_R}(-1)\xrightarrow{A} H\otimes_R\theo_{S_R}
\xrightarrow{B} L\otimes\theo_{S_R}(1)
\ed
is a family of geometrically semistable Kronecker complexes
parametrized
by $U$. By
Corollary~\ref{monad families}, $\bK$ is a family of torsion-free monads.  
The map $\HH^0(A(1)):K\rightarrow H\otimes_R S_1$
satisfies $ \HH^0(A(1))\otimes F = \HH^0(A(1)\otimes F)$
for each geometric point $\spec F \rightarrow\spec R$.
Thus $K\otimes_R F \hookrightarrow H\otimes_R S_1\otimes 
_RF$ for each such $F$. Therefore, by \cite[Theorem~6.8]{Ei},
$H\otimes_R S_1/K$ is flat, hence a vector bundle on $U$. 
Thus $K$ is   a subbundle of $H\otimes_R S_1$.

Next, consider the map
\bd \phi(B) : 
H_R\otimes S_1^*\xrightarrow{\HH^0(B)\otimes S_1^*} 
L\otimes S_1\otimes S_1^*
\rightarrow L.
\ed
We claim that $\phi(B)$ is surjective. 
To see this, note that the composite 
map $$\psi:H\rightarrow H\otimes S_1^*\otimes S_1 
\xrightarrow{\phi(B)\otimes S_1} 
L\otimes S_1$$ is just $\HH^0(B)$. Since 
$B=\HH^0(B)\otimes \theo_{S_R}$, this implies that 
$B=\psi\otimes \theo_{S_R}$.
Thus, if  $\phi(B)$ has image  $L'\subsetneq L$,
then 
$\on{Im}(B)\subseteq L'\otimes  \theo_{S_R}(1).$
Since $(L/L')\otimes S(1)$
is not a bounded $S$-module it has a nonzero image in $\qgr S_R$,
contradicting the surjectivity 
of $B$. Thus $\phi(B)$ is surjective and $L$ is a quotient bundle 
of $H\otimes_RS_1$.

Thus $\bK\mapsto (K,L)$  gives a 
 map $\alpha: \widehat{\cK}^{ss}\rightarrow\GR$.
 It is easily checked
that the construction of Remark~\ref{inv construction} sends
this pair $(K,L)$ back to the original Kronecker complex $\bK$. Thus,
$\alpha$ is injective and (by Lemma~\ref{subscheme N})
 $\alpha$ factors through $\N$.
Since a geometrically semistable monad 
certainly satisfies $\eqref{dagger}$ at each geometric point,
the image of $\alpha$ is contained in $\cN'$.  
The equivariance of this embedding of functors follows from
the last sentence of 
Remark~\ref{equivariance}.

It remains to show that $\alpha$ surjects onto $\cN'$.  Suppose that
$(K,L)\in \cN(\spec R)$ is a KL-pair for which $(K\otimes F,L\otimes F)$
satisfies \eqref{dagger} for every geometric point 
$\spec F\rightarrow\spec R$. 
Let $\bK$ be the Kronecker complex determined by $(K,L)$.
Suppose that
\bd
\bK': \;\; K'\otimes \theo(-1)\rightarrow H'\otimes\theo 
\rightarrow L'\otimes\theo(1)
\ed
is a maximal subcomplex of $\bK_F$ for some geometric point $\spec F$, as
defined in the proof of 
 Proposition~\ref{semistable module from complex}.
 Then, as in \cite[Lemme~2.1]{Le Potier},
it is straightforward to check that
$\bK'$ is the subcomplex associated to $H'$ 
by Construction~\ref{H' construction}.
Since $\bK_F$
satisfies $\eqref{dagger}$, it follows that 
$\on{rk}(\bK_F)p_{\bK'} - \on{rk}(\bK')p_{\bK_F}\leq 0$,
implying that $\bK_F$ is semistable.  Consequently $\bK$
is geometrically semistable.

(2) The proof for stable complexes is essentially the same.
\end{proof}

We next prove that   $\cN'$ is exactly the semistable locus of $\N$ in 
the GIT sense for the action of $\on{SL}(H)$.  We proceed as follows.
For each fixed $r$, $c_1$ and $\chi$ there are only finitely
many possible values of $r'$, $c_1'$ and $\chi'$ that can occur
as invariants associated to a subspace $H'\subset H$.
Thus, for any integer $m=m(r,c_1,\chi)\gg 0$
and  for all possible values of $r'$, $c_1'$ and $\chi'$
associated to subspaces of $H$, the inequality in
$\eqref{dagger}$ is satisfied if and only if it is satisfied for this
fixed value $m(r,c_1,\chi)$ of $m$.

The variety $\GR$ has natural ample line bundles $\theo(k,\ell)$ obtained by
pulling back $\theo(k)$, $\theo(\ell)$ from projective spaces under the
Pl\"ucker embeddings of the two factors of $\GR$
 (see \cite[Example~2.2.2]{HLbook}) and these bundles are equivariant
for the diagonal $\on{GL}(H)$-action.
The open subset of ${\mathcal N}$ consisting
of (semi)stable
points (see \cite[Definition~2]{Seshadri}) 
for the action of $\on{SL}(H)$ under the 
linearization $\theo(k,\ell)$ of $\on{GR}$
will be written $\cN^{s}$, respectively $\cN^{ss}$.\label{ns-defn}

 One now  has a result 
analogous to  \cite[Th\'eor\`eme~3.1]{Le Potier}.   

\begin{prop}\label{semistability via GR}
Let $S\in \ASthree$, pick integers $\{r\geq 1,c_1,\chi\}$  and 
choose $m=m(r,c_1,\chi)\gg 0$  as above.
Set $k=(2m+3)(c_1+r) - n$ and $\ell=-(2m+3)(c_1-r) +n$.  
Then $\widehat{\cK}^{s}=\cN^s$ and 
$\widehat{\cK}^{ss}=\cN^{ss}$ under the linearization  $\theo(k,\ell)$ of
$\on{GR}$.
\end{prop}
\begin{proof} Let $F$ be an algebraically closed field. 
By  \cite[Lemme~3.3]{Le Potier}, a pair 
$(K,L)\in\on{GR}(\spec F)$ 
is (semi)stable for the $\on{SL}(H)$-action 
if and only if 
  \begin{equation}\label{GR semistability}
k[n\dim K'-n'\dim K] - \ell[n\dim L' - n'\dim L] < 0 \;
 (\text{respectively }\leq 0)
 \end{equation}
for every 
proper nonzero $H'\subset H_F$ of dimension $n'$ (where  $K'$ and
$L'$ are defined by Construction~\ref{H' construction}).

Now consider the point $(K,L)$ associated to a Kronecker complex
$\bK$. Use \eqref{lin-alg} to rewrite the left hand side of \eqref{GR
semistability} in terms of $r,\dots,d_{1}'$. After  simplification, 
this gives:
\begin{equation}\label{GR-semi2}
k\left[n\dim K' -n'\dim K\right] -\ell \left[n\dim L'-n'\dim L\right]
= Z(H'),
\end{equation}
 where 
$Z(H')=2n\left[r(c_1' m+\chi') - r'(c_1m+\chi)\right].
$
By Lemma~\ref{description of Nss} and the  choice of $m$, 
$Z(H') < 0$ (respectively $Z(H')\leq 0$) for every $H'\subset H$ if and only
if the Kronecker complex $\bK$ associated to $(K,L)$ is (semi)stable.
By the first paragraph of the proof, \eqref{GR-semi2}
 implies that the  (semi)stability  
of the Kronecker complex $\bK$ is equivalent to the 
(semi)stability  
of the associated point $(K,L)$ of the Grassmannian.
\end{proof}

\begin{corollary}\label{khat-fine}
Let $S\in \ASthree$ and pick integers $\{r\geq 1,c_1,\chi\}$. Then $\cN^s $ and 
$\cN^{ss}$ are fine moduli spaces for isomorphism classes of
  geometrically (semi)stable
framed Kronecker complexes with invariants $\{r,c_1,\chi\}$.\qed
\end{corollary}

\subsection{Existence of Moduli Spaces} \label{moduli-subsec}

Corollary~\ref{khat-fine} does not directly give fine moduli spaces for 
isomorphism  classes of 
$S$-modules, since Kronecker complexes cannot be canonically framed in families.
(In fact, as in the commutative setting \cite{HLbook}, 
one can only hope to get a fine  moduli space
for equivalence classes of modules, as defined in Definition~\ref{typical
Kronecker-defn}.) However, using the results of 
Sections~\ref{monads and beilinson} and  \ref{section3} it is 
not hard   to produce the relevant
  coarse and fine moduli spaces. 
 
Before stating the result, we show that, although the  last
subsection required the  invariants $\{r,c_1,\chi\}$ to satisfy
 $-r<c_1\leq 0$, the results proved there can be applied to $S$-modules
with any prespecified invariants. We emphasize, however, that the spaces
$\cN^s$, etc, have only been defined for invariants satisfying the hypotheses of 
Notation~\ref{invariants-notation} and we do not wish to define them more
generally.

\begin{remark}\label{other-c1} Let $S\in \ASthree$ and 
 pick a commutative noetherian $k$-algebra $R$. Write 
 $\cC(r,c_1,\chi)$ for   the category of  $R$-flat
families of geometrically semistable (or stable) 
torsion-free modules $\cM\in \qgr S_R$
with invariants $\{r,c_1,\chi\}$.
If $\cM \in \cC(r,c_1,\chi)$, it follows from 
Lemma~\ref{hilbert polynomial'}(2) and Corollary~\ref{two term} that 
$c_1(\cM(1)) = c_1+r $ and 
$\chi((\cM(1)) =\chi+2r+c_1$.
Thus, by induction, 
there exist unique invariants $\{r,c_1',\chi'\}$ such that $-r<c_1'\leq 0$ and 
$\cC(r,c_1,\chi)\simeq \cC(r,c'_1,\chi')$.
\end{remark}

The next result proves Theorem~\ref{firstthm} from the introduction.
\begin{thm}\label{projective moduli spaces}
Let $S\in \ASthree$ and fix integers 
 $r,c_1,\chi$, with $r\geq 1$.
 \begin{enumerate}
 \item There exists a projective coarse moduli space\label{moduli-1-defn}
 $\M^{ss}_S(r,c_1,\chi)$ for equivalence classes of 
  geometrically semistable torsion-free 
modules in $\qgr S$ of rank $r$, first Chern class $c_1$ 
and Euler characteristic $\chi$. 
 \item   $\M^{ss}_S(r,c_1,\chi)$ contains an open subscheme 
$\M^s_S(r,c_1,\chi)$ that is a coarse moduli space
for  the geometrically stable modules.
 \item If $-r <c_1\leq 0$ then 
 $\M^{ss}_S(r,c_1,\chi)=\N^{ss}/\!\!/\on{PGL}(H)$
 and $\M^s_S(r,c_1,\chi)=\N^s/\!\!/\on{PGL}(H)$.
 \end{enumerate}
\end{thm}

\begin{proof} By Remark~\ref{other-c1} we may assume in parts (1) and (2) 
that  $-r<c_1\leq 0$ and hence that our modules are normalized.  By
Proposition~\ref{semistable module from complex}, the moduli functor  for
equivalence classes  of geometrically semistable torsion-free modules in $\qgr
S$ with invariants $\{r,c_1,\chi\}$ is now 
isomorphic to $\cK^{ss}$, the functor
for the  corresponding Kronecker complexes. So, it suffices to prove the
theorem for $\cK^{ss}$.

The group functor  $\on{PGL}(H)$ acts on $\widehat{\cK}^{ss}$,
and there is a (forgetful) map of functors 
$\widehat{\cK}^{ss}/\on{PGL}(H)\to\cK^{ss}$.
Although this map is not an isomorphism of functors, we claim that
it is a {\em \'etale local isomorphism}
in the sense that it induces an isomorphism of the sheafifications
for the \'etale topology of affine $k$-schemes \cite[p.60]{Simpson}.  
To see this, observe
that the map of functors is injective, so it suffices to prove that
it is \'etale locally surjective. Let
$$\bK : \theo_{S_R}(-1)\otimes_R V_{-1} \to \theo_{S_R}\otimes
 V_0\to \theo_{S_R}(1)
\otimes
V_1$$  be an $R$-flat family
of geometrically semistable Kronecker complexes,
as in \eqref{typical Kronecker family}.  
 After an \'etale base change $R\rightarrow R'$,
the $R$-module $V_0$ becomes trivial and so $\bK\otimes_R R'$
is in the image of $\widehat{\cK}^{ss}(R')$, proving the claim.

By \cite[p.60]{Simpson}, this implies that any scheme $Z$  that
corepresents $\widehat{\cK}^{ss}/\on{PGL}(H)$ also
 corepresents $\cK^{ss}$.
But Proposition~\ref{semistability via GR}
implies that $\mathcal{N}^{ss}=\widehat{\cK}^{ss}$,
 while \cite[Theorem~4.2.10]{HLbook}
implies that the GIT quotient $Z={\mathcal N}^{ss}/\!\!/\on{PGL}(H)$
does  corepresent $\widehat{\cK}^{ss}/\on{PGL}(H)$. 
So  $Z$ corepresents
$\cK^{ss}$ and gives the desired moduli space. 
\end{proof}

\begin{remark}\label{stack-quotient} Pick integers $\{r,c_1,\chi\}$
with $-r<c_1\leq 0$.
 For future reference, note that, by
Corollary~\ref{khat-fine} and third paragraph of the proof of
Theorem~\ref{projective moduli spaces}, the stack-theoretic quotient
$[\cN^s/GL(H)]$ equals the moduli stack of geometrically stable Kronecker
complexes with invariants $\{r,c_1,\chi\}$.
\end{remark}

\begin{corollary}\label{loc free rmk} Let $S=S(E,\cL,\sigma)\in \ASthree'$
and pick integers  $\{r\geq 1,c_1,\chi\}$.
Then there is an open subscheme of 
$\M_S^s(r,c_1,\chi)$ parametrizing  modules whose
restriction to $E$ is locally free.
\end{corollary}
\begin{proof}
Suppose that one is given a 
flat family of coherent sheaves on $E$ parametrized by a scheme $X$.
 Then the vector bundles in that family are parametrized by an open
subset of $X$. Thus, the corollary follows from the next lemma.
\end{proof}

\begin{lemma}\label{loc free rmk-lemma} Let $S=S(E,\cL,\sigma)\in \ASthree'$.
Let $\cM$ be an $R$-flat family of torsion-free objects in $\qgr S_R$.
 Then $\cM|_{E}$ is 
$R$-flat.
\end{lemma}

\begin{proof} By \cite[Lemma~E5.3]{AZ2}, 
$M=\Gamma^*(\cM)_{\geq n}$ is a flat $R$-module for any $n\gg 0$ and so it
suffices to show that $M/Mg$ is $R$-flat.
By \cite[Lemma~C1.12]{AZ2} it is enough to show that  $\Tor_1^R(M/Mg, R/P)=0$
for every prime ideal $P\subset R$. 
By \cite[Proposition~C1.9]{AZ2} we may assume that $R=(R,P)$ is local. 

Write $S$ for $S_R$. By hypothesis, $\cM/\cM P$ is a torsion-free module in $\qgr S_{R/P}$
 and so, possibly after increasing $n$, $M/MP$ is a torsion-free module in $\gr
 S_{R/P}$.
 Thus $g$ acts without
torsion. Equivalently, $M g\cap MP = M Pg$ and so the natural homomorphism 
$MP\otimes_{S} S/gS \to M\otimes_S S/gS$ is injective.
As $M_R$  is flat, $MP\cong M\otimes_RP$ and so the natural map
$$\theta:(M\otimes_RP)\otimes_S S/gS \to M\otimes_SS/gS\cong M/Mg$$ is injective.
Since the actions of $S$ and $R$ on $M$ commute, 
 $$(M\otimes_RP)\otimes_S S/gS\cong (M\otimes_SS/gS)\otimes_R P
 \cong M/Mg\otimes_RP.$$
Thus $\theta$ is just the natural map $M/Mg\otimes_RP \to M/Mg\otimes_RR$. 
In particular,  $\Tor_1^R(\cM|_{E}, R/P)=\pi\on{Ker}(\theta)=0$.
\end{proof}

In many cases the stable locus of the moduli space of Theorem~\ref{projective
moduli spaces} is a fine moduli space and the next proposition gives one such 
example. The general technique for producing fine moduli spaces is nicely
explained in \cite[Section~4.6]{HLbook}, but we will give a detailed proof of
this result for the benefit of the reader who is unfamiliar with the main
ideas.
We first need a variant of the 
 standard fact that stable objects have trivial endomorphism rings when
the ground field $k$ is algebraically closed.

\begin{lemma}\label{stable is simple}
Suppose that $\cF$ is a geometrically stable torsion-free
object in $\qgr S_F$ for a field $F$.  Then the natural map 
$F\rightarrow\Hom_{\qgr S_F}(\cF,\cF)$ is an isomorphism.
\end{lemma}
\begin{proof}
Suppose that $F\rightarrow F'$ is a field extension and that 
$\phi: \cF\otimes F' \rightarrow \cF\otimes F'$
is a nonzero endomorphism.  By the observation before Lemma~\ref{JH},
$\cF\otimes F'$ is stable. Since $\on{Im}(\phi)$
is both a subobject and a quotient object of $\cF\otimes F'$, this implies that
$\on{rk}(\cF\otimes F')p_{\on{Im}(\phi)} - \on{rk}(\on{Im}(\phi))
p_{\cF\otimes F'}$ is both positive and negative. This contradicts stability 
 unless  
$\on{Im}(\phi) = \cF\otimes F'$.  Thus, $\phi$ is an automorphism
and $\End(\cF\otimes F')$ is a division ring containing $F'$.  

By \cite[Theorem~7.4]{AZ}, 
 $\End(\cF\otimes F')$ is a finite-dimensional  $F'$-vector space.
 Thus   $F'=\End(\cF\otimes F')$ if $F'$ is algebraically closed.
In particular, if $\overline{F}$ is the  algebraic closure of $F$, then
$\End(\cF)\otimes \overline{F} 
\subseteq \End(\cF\otimes \overline{F}) =\overline{F}$ and 
so the natural map $F\rightarrow \End(\cF)$ must
be an isomorphism.  
\end{proof}

\begin{prop}\label{finemodulispace}
Let $S\in \ASthree(k)$, where $\on{char}\,k=0$.  Fix integers 
$\{r,c_1,\dots, d_1\}$  by Notation~\eqref{invariants-notation}.
If $\on{GCD}(d_{-1},n,d_1)=1$, then
 $\M_S^s(r,c_1,\chi)$ is a fine moduli space whenever it is nonempty.
\end{prop}

\begin{proof}  
Theorem~\ref{projective moduli spaces} defines a map of functors 
$a:\cK^s\rightarrow \M^s$.  We think of both
functors as presheaves on the category of (noetherian affine)
$k$-schemes in the \'etale topology. Then $\M^s$ is a sheaf
since it is represented by a scheme.  We will construct a map
$b:\M^s\rightarrow \cK^s$ of functors and show that
\begin{enumerate}
\item[(i)] $b\circ a = 1_{\cK^s}$
\item[(ii)] $a$ is \'etale-locally surjective, in the sense that
 for any affine noetherian $k$-scheme $U$ and  map $x\in \Hom(U,\M^s)$
there is an \'etale cover $V\rightarrow U$ and object $\tilde{x}$
in $\cK^s(V)$ such that the image of $\tilde{x}$ in $\Hom(V,\M^s)$
is $x|_V$. 
\end{enumerate}
The isomorphism of $\cK^s$ and $\M^s$ is then a special case
of the following general fact (proved by a simple diagram chase):
suppose that $F$ is a presheaf of sets and $G$ is a sheaf of sets.  
Suppose further that $a:F\rightarrow G$ and $b:G\rightarrow F$
are natural transformations of these presheaves such that
$b\circ a = 1_F$ and such that $a$ is  \'etale-locally surjective.
Then $F=G$. 

We first prove (ii). Let $x\in\Hom(U,\M^s)$.
As in the proof of \cite[Corollary~4.2.13]{HLbook}, 
Luna's \'etale slice theorem together with Lemma~\ref{stable is simple}
 implies that the projection
$\N^s\rightarrow \M^s$
is a principal $\on{PGL}(H)$-bundle.  It follows that 
there exists an \'etale cover $V\rightarrow U$
such that the fibre product $\N^s\times_{\M^s} V$ is a 
$\on{PGL}(H)$-bundle with a section.  Hence $x|_V$ is in the
image of the map $\N^s(V)\rightarrow \M^s(V)$.  But this 
map factors through $\cK^s(V)$, so there is an element
 $\tilde{x}\in \cK^s(V)$ lifting $x|_V$.

We next construct the map $b$.
By Yoneda's Lemma, in order to produce $b$ it is enough to produce
a family $\cU$ of geometrically stable Kronecker complexes in 
$\qgr S\otimes\theo_{\M^s}$.  
The universal sub and quotient modules $K$ and $L$ on 
$\GR$ are $\on{GL}(H)$-equivariant and the centre ${\mathbf G}_m$
acts with weight $1$ in their fibres.
Pick $a,b,c\in{\mathbb Z}$ satisfying $ad_{-1}+ bd_1 + cn=1$.
Then $$M = (\det K)^{\otimes a} \otimes (\det L)^{\otimes b}
\otimes (\det H)^{\otimes c}$$ is a $\on{GL}(H)$-equivariant
line bundle on $\GR$ and ${\mathbf G}_m$ acts with weight
$1$ in its fibres (see Remark~\ref{equivariance}).  Hence 
$$\cP'=(K\otimes M^*\rightarrow S_1\otimes H_{\N^s}\otimes M^*,\
S_1^*\otimes H_{\N^s}\otimes M^* \rightarrow L\otimes M^*)$$
forms a pair of $\on{GL}(H)$-equivariant maps on $\N^s$ 
with  ${\mathbf G}_m$ acting trivially on the sheaves.

These sheaves and maps $\cP'$ are therefore   $\on{PGL}(H)$-equivariant.
Since $\N^s\rightarrow\M^s$ is a principal $\on{PGL}(H)$-bundle,
$\cP'$ descends to $\M^s$.  Now Remark~\ref{inv construction}
associates to $\cP'$ a diagram $\cU$ in $\qgr S\otimes\theo_{\M^s}$.   
The pullback of $\cU$ 
to $\N^s$ is the Kronecker complex associated to the universal
KL-pair on $\N^s$ and    $\N^s\rightarrow\M^s$ is faithfully flat. Thus
\cite[Lemma C1.1]{AZ2} implies that $\cU$ is a complex.  Furthermore, any map 
$\spec F\rightarrow \M^s$ for a field $F$ lifts to $\spec F\rightarrow\N^s$. 
Thus $\cU$ is a geometrically stable Kronecker complex since its pullback
to $\N^s$ is.  This defines the map $b$.

Finally, we will prove that 
 $\cU$ is {\em weakly universal} in the following sense:
any family of geometrically stable Kronecker complexes 
in $\qgr S\otimes\theo_U$ for a noetherian
affine $k$-scheme $U$ is equivalent to the pullback of 
$\cU$ under the induced composite
$U\rightarrow \cK^s \rightarrow \M^s$.
Yoneda's Lemma will then imply that $b\circ a=1_{\cK^s}$.

Suppose that $U$ is an affine noetherian $k$-scheme with
a family $\bK$ of geometrically stable Kronecker complexes 
\eqref{typical Kronecker family}.  Let $\on{Fr}(V_0)\xrightarrow{p}U$
denote the bundle whose fibre over $u\in U$ consists of 
vector space isomorphisms of the fibre $(V_0)_u$ with $H$.
This is a principal $\on{GL}(H)$-bundle over $U$ and the pullback
of $V_0$ to $\on{Fr}({\mathcal H})$ comes equipped with a
canonical
isomorphism with $H\otimes\theo_{\on{Fr}(V_0)}$.  Thus,
$\bK$ pulls back to a framed Kronecker complex $p^*\bK$
on $\on{Fr}(V_0)$.  This pullback determines a $\on{GL}(H)$-equivariant map
$\on{Fr}(V_0)\rightarrow \N^s$ by Corollary~\ref{khat-fine} and so  we 
obtain a commutative diagram
\bd
\xymatrix{ \on{Fr}(V_0)\ar[d]^{p}\ar[r]^{f_1} 
& \N^s\ar[d]^{q}\\
U\ar[r]^{f_2} & \M^s}
\ed
where $f_2$ is the canonical composite $U\rightarrow \cK^s\rightarrow \M^s$.
By construction, then, we have $\on{GL}(H)$-equivariant isomorphisms
\bd
p^*\bK = f_1^*(q^*\cU\otimes M) = p^*f_2^*\cU \otimes f_1^*M.
\ed
Since $f_1^*M$ is a $\on{GL}(H)$-equivariant line bundle on 
$\on{Fr}(V_0)$, it descends to a line bundle $\overline{M}$
on $U$.  So  there is a $\on{GL}(H)$-equivariant isomorphism
$p^*(\bK \otimes \overline{M}^*) \cong p^*f_2^*\cU$.
  This isomorphism then descends
to show that $\bK\otimes\overline{M} \cong f_2^*\cU$. In other words,
 $\bK$ is equivalent
to the pullback of $\cU$ from $\M^s$ by the natural map.  
This shows that  $b\circ a = 1_{\cK^s}$ and completes the proof of the 
proposition.
\end{proof}

An analogue of Proposition~\ref{finemodulispace}   holds for all values of the
invariants $\{r\geq 1,c_1,\chi\}$ but is awkward to state precisely since 
the condition on GCD's must be assumed for the 
normalization of the given modules. The next result gives two special cases where 
it is easy to unravel this condition.

\begin{corollary}\label{finemodulispace'}
Let $S\in \ASthree(k)$ with $\on{char}\,k=0$ and fix integers 
$\{r\geq 1,c_1,\chi\}$. 
Assume that \emph{either}  
{\rm (1)} $c_1=0$ and $\on{GCD}(r,\chi)=1$ \emph{or} 
{\rm (2)} $\on{GCD}(r,c_1)=1$.

Then
 $\M_S^s(r,c_1,\chi)$ is a fine moduli space whenever it is nonempty.
\end{corollary}

\begin{proof} (1)  
Use \eqref{lin-alg} to see that this case is already covered by 
Proposition~\ref{finemodulispace}.

(2) By Remark~\ref{other-c1}, the hypothesis $\on{GCD}(r,c_1)=1$
is unchanged when one shifts a module and so we can assume that 
$-r<c_1\leq 0$. Now, \eqref{lin-alg} implies that $\on{GCD}(d_{-1},n,d_1)=1$
and so the result follows from Proposition~\ref{finemodulispace}.
\end{proof}

The next result gives a simple application of this observation.

\begin{corollary}\label{stable in rank 1}
Let $S\in \ASthree(k)$, where $\on{char}k=0$.
For   $n\geq 1$, $\M^{ss}(1, 0, 1-n)=\M^{s}(1, 0, 1-n) $ is a  nonempty,
projective, fine moduli space
for equivalence  classes of rank one torsion-free modules $\cM$ 
in $\qgr S$ with $c_1(\cM)=0$ and $\chi(\cM)=1-n$.

If $\qgr S\simeq \coh(\mathbf{P}^2)$, then   $\M^{ss}(1, 0, 1-n)$
is isomorphic to the Hilbert scheme $(\mathbf P^2)^{[n]}$ 
of $n$ points on $\mathbf{P}^2$. 
\end{corollary}

\begin{proof} By Corollary~\ref{monad6}(1),   $\M^{ss}(1, 0, 1-n)
=\M^s(1, 0, 1-n)$.
We next check that
there exist torsion-free modules with the specified invariants. Given a point
module $P\in \gr S$,  \cite[Proposition~6.7(i)]{ATV2} implies that  it has
a resolution $$0\to \cO_S(-2)\to \cO_S(-1)^2\to \cO_S\to P\to 0,$$ from which
it follows that $c_1(P)=0$ but $\chi(P)=1$. Thus, if one takes $n$
nonisomorphic point modules $P_i\cong S/I_i$, then  $M=\bigcap I_i$ is a
torsion-free module with $c_1(M)=0$ and  $\chi(M)=1-n$, as required. Now,
Theorem~\ref{projective moduli spaces} and Corollary~\ref{finemodulispace} 
 combine to prove that $\M^{ss}(1, 0, 1-n)$ is a
projective fine moduli space.
 The final assertion is standard---see for example 
\cite[Example 4.3.6]{HLbook}.
\end{proof}

Suppose one has a family $S_\cB(E,\cL,\sigma)$ of algebras in $\ASthree$  that
give a deformation of the polynomial ring $k[x,y,z]$.  Then, as we will show in
the next section,  one can view Corollary~\ref{stable in rank 1} as showing
that $\M^{ss}(1, 0, 1-n)$ is a deformation of  $(\mathbf P^2)^{[n]}$. 
Since we proved that  $\M^{ss}(1, 0, 1-n)$ is nonempty by
finding the module corresponding to $n$ points on $E$, this may not seem so
surprising. However this analogy does not work for 
the subset $(\mathbf P_S\smallsetminus E)^{[n]}\subset \M^{ss}(1, 0, 1-n)$ 
that deforms $(\mathbf P^2\smallsetminus E)^{[n]}$. The reason is that, when
$|\sigma|=\infty$,   the
only ``points'' in $\qgr S$  are those that lie on $E$ and so, necessarily, 
the modules parametrized by $(\mathbf P_S\smallsetminus E)^{[n]}$ have to be
more subtle. In fact they are line bundles---see
 Theorem~\ref{connectedness1}(3). These modules   correspond in turn to 
the projective right ideals of the ring $A(S)$ that  we considered in
Section~\ref{sectionthereandback}.

\section{Moduli Spaces: Properties}\label{moduli-smooth}
In this section we study the more detailed properties of the moduli spaces
constructed in the last section. We prove that, if nonempty, they are smooth
and behave well in families (see Theorem~\ref{smooth family}).
Moreover, in many cases they are irreducible and hence connected
(see Proposition~\ref{connectedness}). In  rank one, we are able to
give a much more detailed picture of $\M^{ss}(1, 0, 1-n)$
and its open subvariety $(\mathbf P_S\smallsetminus E)^{[n]}$, thereby 
 justifying the comments from the end of the
last section.  

\subsection{Basic Properties}

We first consider a family of AS regular algebras $S_\cB$ parametrized by a
$k$-scheme $\cB$, where   $\on{char}\,k=0$, since we wish to study the 
global structure of our moduli spaces.
 Pick invariants
$\{r,c_1,\chi\}$ as in Notation~\ref{invariants-notation}.  
Then Lemma~\ref{subscheme N} and general properties of GIT 
quotients in the relative setting \cite[Remark 8]{Seshadri} 
yield a subscheme  
$\N^{ss}_{\cB}$\label{smooth-subdefn}
of $\on{GR}_\cB$ such that $\N^{ss}_\cB/\!\!/\on{PGL}(H)$\label{M_S-defn}
  is a universal categorical quotient in the sense that
for any $\cB$-scheme $\cB'\rightarrow \cB$, the fibre product
$(\N^{ss}_{\cB}/\!\!/\on{PGL}(H))\times_{\cB}\cB'$ corepresents
the quotient functor $\N^{ss}_{\cB'}/\on{PGL}(H)$.
 In particular, the fibre of this quotient over $b\in\cB$
is $\M^{ss}_{S_b}(r,c_1,\chi)$ and the same holds for the stable locus. 
 We write
$\M_{\cB}^s(r,c_1,\chi) = \N^s_{\cB}/\!\!/\on{PGL}(H)$.
Applying Luna's \'etale slice theorem as in \cite[Corollary~4.2.13]{HLbook},
and using Lemma~\ref{stable is simple},
shows that the projection $\N^s_{\cB}\rightarrow \M_{\cB}^s(r,c_1,\chi)$ 
is a principal $\on{PGL}(H)$-bundle.
When $\{r\geq 1,c_1,\chi\}$ are arbitrary integers, use Remark~\ref{other-c1}
to choose the unique integers $\{r,c_1',\chi'\}$ with 
$-r<c_1'\leq 0$ and $\cC(r,c_1,\chi)\simeq\cC(r,c_1',\chi')$. 
Then set $\M_{\cB}^{ss}(r,c_1,\chi)=\M_{\cB}^{ss}(r,c_1',\chi')$ and 
 $\M_{\cB}^s(r,c_1,\chi)=\M_{\cB}^s(r,c_1',\chi')$.
 
 We are now ready to prove Theorem~\ref{smooth family-intro} from the
 introduction.
\begin{thm}\label{smooth family}  Suppose that $\on{char}\,k=0$
and that $S$ is a flat family of algebras in $\ASthree$ parametrized by a
$k$-scheme $\cB$.
If $\{r\geq 1,c_1,\chi\}\subset \mathbb Z$, then 
 \begin{enumerate}
  \item $\M^s_{\cB}(r,c_1,\chi)$ is smooth over $\cB$.
 \item  If $p\in \cB$ then
  $\M^s_{\cB}(r,c_1,\chi)\otimes_{\cB} k(p)=\M^s_{S_b}(r,c_1,\chi)$, and
  similarly for $\M^{ss}_{\cB}(r,c_1,\chi).$
\end{enumerate}\end{thm}

\begin{proof} Part (2) is a general property of universal quotients
\cite{Seshadri},
so only part (1) needs proof.
We may assume that the invariants satisfy $-r<c_1\leq 0$. 

Since the projection   $\N^s_{\cB}\to \M_{\cB}^s(r,c_1,\chi)$ is a
principal bundle, \cite[Proposition~IV.17.7.7]{EGA}
 implies that 
it is enough to prove that
$\N_{\cB}^s$ is smooth over $\cB$. To prove this we will use the local criterion
for smoothness. Let
  $R'$ be a local commutative $k$-algebra with a noetherian 
factor ring $R=R'/I$, where $I^2=0$.
  Suppose that
\bd
\bK_R:\;\;\; \theo_{S_R}(-1)\otimes V_{-1} \xrightarrow{d^1_R}
 \theo_{S_R}\otimes V_0 \xrightarrow{d^2_R} \theo_{S_R}(1)\otimes V_1
\ed
is an $R$-flat family of stable monads.
Because $R=(R,\mathfrak{m})$ is local, we may assume that the  $V_{i}$ are 
$k$-vector spaces and the tensor products are over $k$. 

Suppose that  we can lift $\bK$  to a complex
\begin{equation}\label{smooth-family2}
\bK_{R'}:\;\;\; \theo_{S_{R'}}(-1)\otimes V_{-1}\xrightarrow{d^1_{R'}}
\theo_{S_{R'}}\otimes V_0 \xrightarrow{d^2_{R'}}
\theo_{S_{R'}}(1)\otimes V_1
\end{equation}
satisfying $\bK_{R'}\otimes_{R'} R \cong \bK_{R}$. 
Then
$\bK_{R'}\otimes_{R'} I \cong \bK_R\otimes_R I$ and so 
Lemma~\ref{beilinson5-1} implies that 
 $\HH^n(\bK_{R'}\otimes_{R'}I) = \HH^n(\bK_R)\otimes_R I$.  
From the long exact cohomology sequence
associated to the exact sequence of complexes 
$$
0\rightarrow \bK_{R'}\otimes I \rightarrow \bK_{R'}
 \rightarrow \bK_R\rightarrow 0,
$$
and the vanishing of $\HH^n(\bK_R)$ for $n\neq 0$, we find that $\bK_{R'}$
 is 
automatically a monad. Moreover,  from the local criterion for flatness 
\cite[Proposition~C7.1]{AZ2}, the cohomology of $\bK_{R'}$ is flat
 over $R'$.  Thus, $\cN^s$ satisfies the lifting property and so 
 $\cM^s_{\cB}(r,c_1,\chi)$ is smooth.

So it remains to show that we can always lift $\bK$ to $\bK'$.
Choose any lift of $\bK_R$ to a diagram $\bK_{R'}$ of the form 
\eqref{smooth-family2}. We need to adjust the differential  $d_{R'}$
to convert $\bK_{R'}$ into a complex. Note that 
$d_{R'} \in \Hom^1(\bK_{R'},\bK_{R'})$, where
\bd
\Hom^j(\bK_{R'},\bK_{R'}) = \bigoplus_i \Hom(\bK_{R'}^i, \bK_{R'}^{i+j}),
\ed
and  $d_{R'}\circ d_{R'}\in \Hom^2(\bK_{R'},\bK_{R'})$.
The only nonzero term in this product is 
$d^2_{R'}\circ d^1_{R'}$.  Since $d_{R'}$ reduces to $d_R$ mod $I$, we find 
that 
\bd d_{R'}\circ d_{R'}\in 
\Hom^2(\bK_{R'},\bK_{R'}\otimes I)\cong \Hom^2(\bK_R, \bK_R\otimes_R I).
\ed  
As in the proof of Lemma~\ref{exts agree},
 equip  $\Hom^{\bullet}(\bK_R,\bK_R\otimes I)$
with the   differential 
$\delta(f)=  d_{R'} \circ f - (-1)^{\deg(f)} f\circ d_{R'}.$  Then
$\delta(d_{R'}\circ d_{R'}) = 0$ and 
so $d_{R'}\circ d_{R'}$ defines a class in 
$\HH^2(\Hom^{\bullet}(\bK_R,\bK_R\otimes I))$.

Let $\cE=\HH^0(\bK_R)$ and set 
$\tilde{k}=R/\mathfrak{m}$.
Since $\cE$ is geometrically stable, Lemma~\ref{stable is simple} implies that 
$$
\Ext^2_{\Qgr S_{\tilde{k}}}(\cE\otimes \tilde{k},\cE\otimes \tilde{k}) =
 \Ext^0_{\Qgr S_{\tilde{k}}}(\cE\otimes \tilde{k},
 (\cE\otimes \tilde{k})(-3))^* = 0.
$$
By Theorem~\ref{CBC}(1,4) this implies that 
$\Ext^2_{\Qgr S_R}(\cE,\cE)\otimes \tilde{k}  = 0$ and then Nakayama's 
Lemma   gives $\Ext^2_{\Qgr S_R}(\cE,\cE)=0$.  
Since $\Qgr S$ has homological dimension two, $\Ext^3_{\Qgr S_R}(\cE, -)=0$. 
Writing  $I$ as a factor of a free $R$-module, this implies that 
$\Ext^2_{\Qgr S_R}(\cE,\cE\otimes I)=0$.
Consequently, by Lemma~\ref{exts agree},
 there 
exists $\phi\in\Hom^1(\bK_R,\bK_R\otimes I)$ such that 
$\delta(\phi) = d_{R'}\circ d_{R'}$.
Finally, if $D_{R'} = (d_{R'}-\phi) \in\Hom^1(\bK_{R'},\bK_{R'})$,
then $D_{R'} = d_R$ mod $I$ and
$$D_{R'}\circ D_{R'} = (d_{R'} -\phi)\circ (d_{R'}-\phi) =
 d_{R'}\circ d_{R'} - \delta(\phi) =0.$$  Thus 
$D_{R'}$ determines the desired complex $\bK_{R'}$.
\end{proof}

\begin{remark}\label{smooth family-remark} 
The proof of Theorem~\ref{smooth family} also proves its stack-theoretic
analogue. In other words, if one keeps 
  the hypotheses of the proposition and assumes that $-r<c_1\leq 0$, then  the
stack-theoretic quotient 
 $[\cN^s_{\cB}/GL(H)]$ is smooth over $\cB$.
 \end{remark}

Taking $\cB=k$ in the theorem gives:

\begin{corollary}\label{smooth}
If $S$ is in $\ASthree(k)$ where $\on{char} k=0$, then $\M_S^s(r,c_1,\chi)$
is a smooth quasi-projective $k$-scheme for any integers $\{r\geq 1,c_1,\chi\}$.
\qed
\end{corollary}

\begin{corollary}\label{dimension}
Let  $S\in \ASthree(k)$ where $\on{char} k=0$ and pick integers 
 $\{r\geq 1,c_1,\chi\}$. If
 $\M_S^s(r,c_1,\chi)$ is nonempty, then 
 $\dim \M_S^s(r,c_1,\chi) =  r^2 +3rc_1 +c_1^2 -2r\chi+1.$
  \end{corollary}
 
 \begin{proof}
 Pick a geometrically stable module $\cE$ with 
 invariants $\{r,c_1,\chi\}$.  
As in \cite[Theorem 4.5.1 and Corollary 4.5.2]{HLbook},  
 the tangent space to $\M_S^s(r,c_1,\chi)$ at the (necessarily smooth) 
 point $[\cE]$  is 
$T_{[\cE]}\M_S^s = \Ext_{\qgr S}^1(\cE,\cE)$.
By Lemma~\ref{stable is simple}, $\Ext^0(\cE,\cE)=k$ and 
 $\Ext^2(\cE,\cE)=0$.
 Corollary~\ref{two term} now implies that 
$T_{[\cE]}\M_S^s$ has the required dimension 
$1-\chi(\cE,\cE) = r^2 +3rc_1 +c_1^2 -2r\chi+1$.
\end{proof}

We want to examine  
when our moduli spaces are connected. Although we are mainly concerned with   case of 
modules of rank one, the proof actually works for
invariants $\{r,c_1,\chi\}$ for which  
 $\M^{ss}(r,c_1,\chi) = \M^{s}(r,c_1,\chi)$
 and the next lemma provides various cases where this is automatic.
 
 \begin{lemma}\label{semi-equals-stable}
Let $S\in \ASthree$ and consider invariants $\{r,c_1,\chi\}$ with $-r<c_1\leq
0$. Then  $\M^{ss}(r,c_1,\chi)=\M^{s}(r,c_1,\chi)$,
 provided that \emph{either} 
{\rm (1)} $c_1=0$ and $(r,\chi)=1$ \emph{or}
{\rm (2)} $c_1\not=0$ and $(r,c_1)=1$. 
\end{lemma}

\begin{proof} In each case, use \eqref{RR} to show 
that there cannot exist   torsion-free modules 
$0\not=\cF\subsetneq\cM\in \qgr S$  for which
$\on{rk}(\cM)p_\cF =\on{rk}(\cF)p_\cM$.
\end{proof}

\begin{prop}\label{connectedness}
Let $\cB$ be an irreducible $k$-scheme of finite type, 
where $\on{char}(k)=0$,
and  fix  a $\cB$-flat
family of algebras $S_{\cB}\in \ASthree$. Pick integers $\{r\geq 1,c_1,\chi\}.$

\begin{enumerate}
\item
 Suppose that $\M^{s}_{S_b}(r,c_1,\chi)$
is  nonempty for some $b\in \cB$. Then 
$\M^{ss}_{S_{b'}}(r,c_1,\chi)$
is  nonempty for every $b'\in\cB$.
\item Suppose that $\cB$ is a smooth $k$-curve.  
Pick invariants $\{r,c_1,\chi\}$ 
for which $\M^{ss}_{S_\cB}(r,c_1,\chi)=\M^{s}_{S_\cB}(r,c_1,\chi) $.
If $\M^s_{S_b}(r,c_1,\chi)$ is irreducible for some $b\in\cB$,
then $\M^s_{S_{b'}}(r,c_1,\chi)$ is irreducible for every $b'\in\cB$.
\end{enumerate}
\end{prop}

\begin{proof}
Set $\M^{ss} = \M^{ss}_{S_\cB}(r,c_1,\chi)$
and $\M^{ss}_p = \M^{ss}_{S_p}(r,c_1,\chi)\cong 
\M^{ss}\otimes_{\cB} k(p)$ for $p\in \cB$ (and similarly for $\M^s$).

(1)
The morphism $\M^{ss}\xrightarrow{\rho}\cB$ is proper \cite[Theorem~4]{Seshadri} and thus its image is
closed in $\cB$.  On the other hand, Theorem~\ref{smooth family} implies that
$\rho|_{\M^s}: \M^s\rightarrow \cB$
is flat; consequently
 its image is open  in $\cB$ \cite[Th\'eor\`eme IV.2.4.6]{EGA}
 and, by assumption, nonempty.  
Thus $\on{Im}(\rho)$ is closed and contains a nonempty open subset of $\cB$.
Since $\cB$ is irreducible,
this implies that $\on{Im}(\rho)=\cB$.

(2)
By Theorem~\ref{smooth family} and  the fact that
 $\cM^{ss}=\cM^s$, the morphism 
$\M^{ss}\xrightarrow{\rho}\cB$ is flat.  The morphism $\rho$ is
also proper; 
consequently, \cite[Proposition~III.8.7]{H} implies that 
$p_*\theo_{\M^{ss}}$ is a coherent subsheaf of a 
torsion-free $\theo_{\cB}$-module. Thus it is a vector bundle.
 We have $\HH^0(\theo_{\M^{ss}}\otimes k(b)) =
\HH^0(\theo_{\M^{ss}_{b}}) =k(b)$ since the fibre of
$\M^{ss}$ over $b$ is irreducible and projective. Since
$k\subseteq p_*\theo_{\M^{ss}}$, it follows that the map
$p_*\theo_{\M^{ss}}\otimes
k(b)\rightarrow\HH^0(\theo_{\M^{ss}}\otimes k(b))$ is surjective
and so, by Theorem~\ref{CBC}, the fibre of $p_*\theo_{\M^{ss}}$
over $b$ is 1-dimensional.  Thus, $p_*\theo_{\M^{ss}}$ is a line
bundle. By
\cite[Corollary~III.11.3]{H} this implies that $p$ has connected fibres.
Since those fibres are also smooth they must be irreducible.  
\end{proof}

By Corollary~\ref{stable in rank 1}, $\M^{ss}(1,0,\chi)\not=\emptyset$ for
$\chi\leq 1$.   In contrast, even over $\mathbf P^2$  it is a subtle question
to determine when $\M^{ss}(r,c_1,\chi)$ is nonempty for $r>1$ and that question
has been studied in detail in  \cite{Drezet-Le Potier}. In the noncommutative
case,  the question is likely to be similarly subtle, although a number of
positive results can be obtained by combining Proposition~\ref{connectedness}
with   \cite[Th\'eor\`eme~B]{Drezet-Le Potier}.  Since that theorem is rather
technical we merely note:

\begin{corollary}\label{big-nonempty}
 Assume that $\cB$ is an irreducible $k$-scheme of finite type, where
$\on{char} k = 0$,  and  let $S_{\cB}$ be a $\cB$-flat family of algebras in
$\ASthree$ such that $S_b\cong k[x,y,z]$, for some $b\in \cB$. 
If $(r,c_1)=1$, then $\M^{s}_{S_p}(r,c_1,\chi)\not=\emptyset$ for all $p\in \cB$ and all
 $\chi\ll 0$. \end{corollary}

\begin{proof} By \cite[Th\'eor\`eme~B]{Drezet-Le Potier},
 $\M^{s}_{S_b}(r,c_1,\chi)\not=\emptyset$
for $\chi\ll 0$. Now apply Proposition~\ref{connectedness}.
\end{proof}

Using \cite{Drezet-Le Potier}, one can give precise bounds  in this corollary.
As an illustration,  by combining Corollary~\ref{big-nonempty} with  the 
computation in \cite[p.196]{Drezet-Le Potier} shows:

\begin{example}\label{rank20} In
Corollary~\ref{big-nonempty},
$\M_{S_p}(20,9,\chi)\not=\emptyset$ if $\chi\leq 24$.
\end{example}

\subsection{Rank One Modules}
Let $S=S(E,\cL,\sigma)\in \ASthree'$.
In this subsection we examine   $\M_S^s(1,c_1,\chi)$ in more detail
and justify the assertions made at the end of Section~\ref{section5}. 
By Remark~\ref{other-c1}, there is no harm in assuming that $c_1=0$
and we always do so. 
As in the introduction we define 
 $(\boldp_S\smallsetminus E)^{[n]}$\label{hilb-defn} to be the 
subscheme of $\M_S^{ss}(1,0, 1-n)=\M_S^{s}(1,0, 1-n)$ 
parametrizing modules whose restrictions
to $E$ are line bundles. In other words,
 we are concerned with modules in the set $\cV_S$
from Section~\ref{sectionthereandback}. By Corollary~\ref{loc free rmk}, 
$(\boldp_S\smallsetminus E)^{[n]}$ is an open
subscheme of $\M_S^{ss}(1,0, 1-n)$.
As we show next it is always nonempty and generically  consists of
line bundles. 

There are two ways to prove the nonemptiness of
 $(\boldp_S\smallsetminus E)^{[n]}$.  We will give a
constructive proof. The other method is to  prove inductively that 
the complement of $(\boldp_S\smallsetminus E)^{[n]}$ has
dimension less than $2n=\dim\M_S^{ss}(1,0, 1-n).$

\begin{prop}\label{existence} Let $S=S(E,\cL,\sigma)\in \ASthree'$
and let $\epsilon=\epsilon(S)$ \label{epsilon-defn}
denote the minimal integer $m>0$ 
(possibly $m=\infty$) for which  $\cL^{\sigma^m}\cong \cL$ in $\on{Pic}\, E$.

Then for every $0\leq n < \epsilon$,
 there exists a line bundle $\cM_n$ in $\qgr S$ with first Chern class $c_1=0$
and Euler characteristic $\chi = 1-n$.
\end{prop}

\begin{proof} Clearly we may take $\M_0=\theo$, so suppose that $\epsilon>n>0$. 
We will define $\cM_n$ by constructing a short exact sequence 
\begin{equation}\label{existence1}
0\to \cM_n\to \cO(1)\oplus\cO(n)\to \cO(n+1)\to 0.
\end{equation}
Note that $\on{H}^0(\cO(m))=\frac{1}{2}m^2+\frac{3}{2}m+1$ for $m\geq 0$
while $\on{H}^i(\cO(m))=0$, for $m,i\geq 1$,
by Lemma~\ref{chi-lemma}(5).  An
elementary computation therefore shows that $\chi(\cM_n)=1-n$ and  
  $c_1(\cM_n)=0$, as required.

By Lemma~\ref{shifts}, $\cO(m)|_E \cong \cL_m^{\sigma^{-m}}$ for $m\geq 0$.
Suppose first that $n>1$. Since $\cL$ is  very ample, there exists an
injection  $\theta_1 : \cL_{n-1}^{\sigma^{-(n-1)}}\hookrightarrow
\cL_{n}^{\sigma^{-(n)}}.$ As $E$ has dimension one, $\cK=\on{Coker}(\theta_1)$
is a sheaf  of finite length. Since $\cK$ is  a homomorphic image of a very
ample invertible sheaf, it is also a homomorphic image of $\cO$. Let
$\theta_2: \cO\twoheadrightarrow \cK$ be the corresponding map. Since 
$\on{H}^1(E,\cL_{n-1}^{\sigma^{-(n-1)}})=0$,  the map  $\theta_2$ lifts to a
map $\cO_E\to \cL_{n}^{\sigma^{-(n)}}$. Combined with $\theta_1$ we have
therefore constructed a surjection $\theta':  \cO|_E\oplus \cO(n-1)|_E
\twoheadrightarrow \cO(n)|_E$. When $n=1$ the existence of $\theta'$ is just
the standard assertion that the  very ample invertible sheaf
$\cL^{\sigma^{-1}}$ is a homomorphic image of $\cO_E\oplus \cO_E$. Shifting
$\theta'$ by $1$  gives (for all $n\geq 1$)
 the surjection $$\theta: \cO(1)\oplus \cO(n)
\twoheadrightarrow \cO(1)|_E\oplus \cO(n)|_E \twoheadrightarrow \cO(n+1)|_E.$$

 By Lemma~\ref{chi-lemma}(5), the map $\on{H}^0(\cO(m))\to
\on{H}^0(\cO(m)|_E)$ is surjective for all $m\in \mathbb Z$ and so $\theta$
lifts to a  morphism $\phi: \cO(1)\oplus \cO(n)\to \cO(n+1)$.  Let
$\cF=\on{Coker}(\phi)$ and suppose that $\cF\not=0$. By construction,
$\cF|_E=0$ and so, by Lemma~\ref{general-lifting}, $\cF$ has finite length.
Thus $\epsilon<\infty$, by \cite[Proposition~7.8]{ATV2}. Since $\cF$ is a
homomorphic image of $\cO(n+1)/\phi(\cO(n))$, the shift $\cF(-n)$ is a
homomorphic image of a line module, in the notation of \cite{ATV2}. By
\cite[Proposition~7.8]{ATV2} and in the notation of that result,
$\eta=\epsilon(\cF(-n))\geq \epsilon$ and hence, by
\cite[Proposition~6.7(i)]{ATV2},  $\cF(-n)$ has a minimal resolution of the form 
$$\cO(-\eta-1)\to\cO(-1)\oplus \cO(-\eta)\to \cO\to \cF(-n)\to 0.$$
 By Serre duality and \cite[Lemma~4.8]{ATV2}, $\on{H}^i(\cF(m))=0$ 
 for all $i>0$ and
$m\in \mathbb Z$.  A simple computation using Hilbert series then ensures
that  $c_1(\cF)=0$ but $\chi(\cF)=\on{H}^0(\cF)=\eta$. By the
argument of the first paragraph of this proof, this implies  that
$\cM_n=\on{Ker}(\phi)$ satisfies $c_1(\cM_n)=0$ but 
$\chi(\cM_n)=1-n+\eta>1$. By Corollary~\ref{monad6}(1),
 this contradicts Remark~\ref{chi ineq}.

 Thus $\cF=0$ and we have  constructed the desired exact sequence 
 \eqref{existence1}.
 \end{proof}
 
 When $\epsilon\leq n$ we need to combine Proposition~\ref{existence}
  with the idea behind
 Corollary~\ref{stable in rank 1} in order to show that 
 $(\mathbf P_S\smallsetminus E)^{[n]}\not=\emptyset.$
 \begin{corollary}\label{existence2}
 Let $S=S(E,\cL,\sigma)\in \ASthree'$. 
 Then, for any $m\geq 0$, there exists a torsion-free, rank one module
  $\cM_m\in \qgr S$ 
 such that $c_1(\cM_m)=0$, $\chi(\cM_m)=1-m$ and 
 $\cM|_E$ is torsion-free. 
 In particular, $(\mathbf P_S\smallsetminus E)^{[m]}$
 is nonempty.
 \end{corollary}
 
 \begin{proof} Define $\epsilon$ as in Proposition~\ref{existence}. When 
 $\epsilon =1$, \cite[Theorem~7.3]{ATV2} implies that 
 $\qgr S\simeq \coh(\mathbf P^2)$,
  where the result is standard. So 
 we may assume that $\epsilon>1.$ If $0\leq n<\epsilon$, then the module
 $\cM_n$ defined by Proposition~\ref{existence} has the required properties
  and so we may also assume that 
  $\epsilon<\infty$. For $m\geq \epsilon$, write 
 $m=r\epsilon+n$, with $0\leq n<\epsilon$.
 
 As $\cM_n$ is a rank one torsion-free module, 
 $\cM_n\hookrightarrow r^{-1}\cO_S$, for some nonzero $r\in S$; 
 equivalently $\cM_n\hookrightarrow \cO(t)$, for some $t\geq 0$.
 By \cite[Lemma~7.17 and Equation~2.25]{ATV2},
 there exist infinitely many nonisomorphic simple modules 
 $\cN_\alpha\cong  \cO(t)/\cJ_\alpha\in \qgr S$ with 
 $\epsilon(\cN_\alpha ) =\epsilon$. 
 As in the proof of \cite[Theorem~7.3(ii)]{ATV2},
 $\bigcap \cJ_\alpha = 0$ and so
 $\cM_n\cap\bigcap_\alpha ( \cJ_\alpha)=0$.
 We may therefore  chose $\alpha_1,\dots\alpha_r$ such that,
 if  $\cR_m= \cM_n\cap (\bigcap_{i=1}^r \cJ_{\alpha_i})$,
 then 
 $\cM_n/\cR_m \cong \bigoplus_{i=1}^r \cN_{\alpha_i}$. 
 Set $\cM_m=\cR_m$. As in the proof of Proposition~\ref{existence}, 
 $c_1(\cN_{\alpha_i})=0$ and $\chi(\cN_{\alpha_i})=\epsilon$.
 Thus, $c_1(\cM_m)=0$ and $\chi(\cM_m)=1-n-r\epsilon=1-m$ for each $m$.
 Finally,
 \cite[Proposition~7.7]{ATV2} implies that $\cN_\alpha|_E=0$
 and so $(\cM_m)|_E\cong (\cM_n)|_E$ is torsion-free. 
 \end{proof}

We do not know whether there exist line bundles $\cM_m\in \qgr S$ with 
$c_1(\cM_m)=0$ and $\chi(\cM_m)=1-m$, for 
$m\geq \epsilon.$ Such modules clearly do not exist when 
$\qgr S\simeq\coh(\mathbf P^2)$ (where $\epsilon = 1$) and 
 we suspect that they do not exist when $\epsilon>1$.

Combining the last several results gives the following detailed description of 
the spaces  $\M^{ss}_S(1,0,1-n)$
and $(\boldp_S\smallsetminus E)^{[n]}$, and proves Theorem~\ref{mainthm1}
from the introduction.
 
\begin{thm}\label{connectedness1}
Let   $S=S(E,\cL,\sigma)\in \ASthree'(k)$, where $\on{char}(k)=0$.  Then
\begin{enumerate}
\item   $\M^{ss}_S(1,0,1-n)$ is a smooth, projective, fine moduli space
for equivalence  classes of rank one torsion-free modules $\cM\in \qgr S$
with $c_1(\cM)=0$ and $\chi(\cM)=1-n$. Moreover, $\dim \M^{ss}_S(1,0,1-n)=2n$.

\item $(\boldp_S\smallsetminus E)^{[n]}$ is a nonempty open 
subspace of $\M^{ss}_S(1,0,1-n)$.
\item Let $\epsilon=\epsilon(S)$ be defined as in
Proposition~\ref{existence}.
If $n<\epsilon$, then $(\mathbf{P}^2\smallsetminus E)^{[n]}$
parametrizes line bundles with invariants $\{1,0,1-n\}$.
\end{enumerate}
\end{thm}
\begin{proof} Part (1) follows from 
Corollaries~\ref{stable in rank 1}, \ref{smooth} and \ref{dimension}.
Part (2) then follows from Corollaries~\ref{loc free rmk} and 
 \ref{existence2}.

(3)  Let $\M\in \qgr S$ be a torsion-free 
module, with torsion-free restriction to $E$ and invariants $\{1,0,1-n\}$.
If $\M$ is not a line bundle then $M$ is not reflexive, by
Lemma~\ref{vble-reflexive}. Since $\M^{**}/M$  has finite length it must,
 by hypothesis,  have a composition series of 
simple modules $\N_i$ with $\epsilon(\N_i)\geq \epsilon.$ 
Just as in the penultimate paragraph of the proof of
Proposition~\ref{existence}, this forces $\chi(\M^{**})\geq 1-n+\epsilon>1$,
which is impossible. Thus $\M$ is a line bundle.
\end{proof}

When $S$ is a deformation of the polynomial ring $k[x,y,z]$, we can say even 
more about $(\boldp_S\smallsetminus E)^{[n]}$. This next result
 completes the proof of Theorem~\ref{mainthm13}
from the introduction. 
 
\begin{thm}\label{connectedness2}
Let $\cB$ be a smooth curve defined over a field $k$ of characteristic zero and
let $S_\cB=S_{\cB}(E,\cL,\sigma)\in \ASthree'$ be a flat family of algebras 
such that $S_p=k[x,y,z]$ for some point $p\in \cB$.  Set $S=S_b$, for any
point $b\in \cB$. Then:
\begin{enumerate}
\item Both $\M^{ss}_S(1,0,1-n)$
and $(\boldp_S\smallsetminus E)^{[n]}$
are    irreducible and hence connected.
\item  
  $\M^{ss}_S(1,0,1-n)$ is a deformation of the Hilbert scheme 
$(\mathbf{P}^2)^{[n]}$, with its subspace
$(\boldp_S\smallsetminus E)^{[n]}$ being a deformation of 
$(\mathbf{P}^2\smallsetminus E)^{[n]}$.
\end{enumerate}
\end{thm}
\begin{proof}  
The fibre of $\M^{ss}_S(1,0,1-n)$ at the special point 
$p$  is irreducible by \cite[Example~4.5.10]{HLbook}. Thus
both  $\M^{ss}_S(1,0,1-n)$ and its 
open subvariety $(\boldp_S\smallsetminus E)^{[n]}$ are irreducible
by Proposition~\ref{connectedness}, proving (1). Part (2) follows from
Theorem~\ref{smooth family}  the final
assertion of  Corollary~\ref{stable in rank 1}.
\end{proof}

There are no (commutative!) deformations of 
$(\mathbf P^2)^{[1]}=\boldp^2$ and so $\M^{ss}_S(1,0,0)$ must equal $\mathbf
P^2$ in Theorem~\ref{connectedness2}. In fact it is not difficult to prove directly that 
 \begin{equation}\label{case-n-is-one}
 \M^{ss}_S(1,0,0)\cong(\mathbf P^2)^{[1]}=\boldp^2\qquad\text{for all}\quad
  S\in \ASthree.
 \end{equation}
 The proof is left to the interested reader.
 
\subsection{Framed Modules for the Homogenized Weyl algebra}
 As we show in this subsection, the   methods developed in this paper
 give a quick proof that the bijections of
\cite{BW1} and their generalizations in  \cite{KKO} do  come 
from moduli space structures.  This will prove Proposition~\ref{moduli for
others-intro} from the introduction.

 In this subsection we fix $k={\mathbb C}$ and let $U=\mathbb
C\{x,y,z\}/(xy-yx-z^2, z\, \text{central})$ denote the homogenized  Weyl
algebra, as in \eqref{U-defn}.  In this case, rather than  restrict to the
curve $E=\{z^3=0\}$ as we have done previously, we will use the  projection
$U\twoheadrightarrow U/Uz= \mathbb C[x,y]$  to identify $\coh(\pline)$ with a
subcategory of $\qgr U$. Fix   $r\geq 1$ and $\chi\leq 1$. Following
\cite{KKO}, if $R$ is a commutative algebra,  a {\em framed} torsion-free
module in $\qgr U_R$ is a rank $r$ torsion-free  object $\cE$ of $\qgr U_R$
equipped with an isomorphism $\cE|_{\pline} = \theo_{\pline}^r\otimes R$. (This
is a different notion of framing from that considered  for complexes on
page~\pageref{framed-complex}.) A homomorphism $\theta: \cE\to \cF$ between
two  such objects is a   homomorphism in $\qgr U_R$ that induces a  scalar
multiplication   on $\cE|_{\pline}=\theo_{\pline}^r=\cF|_{\pline}$ under the
two framings. Let $\underline{\M}^{\on{fr}}_U(r,0,\chi)$ denote the moduli
functor of framed torsion-free modules in  $\qgr U$ of rank $r$ and  Euler
characteristic~$\chi$.

 Let
$Z={\mathfrak gl}_n\times {\mathfrak gl}_n \times M_{r,n} \times M_{n,r}$
be the vector space of quadruples of matrices of the prescribed sizes and set
\bd
V = \big\{ (b_1, b_2, j, i)\in Z \; 
\big|\; [b_1, b_2] + ij + 2{\mathbf 1}_{n,n} =0
\big\}
\ed
There is a natural action of $\on{GL}(n)$ on $V$.

\begin{prop}\label{moduli for others}
Let $U$ be the homogenized Weyl algebra. Then the quotient $V/\!\!/\on{GL}(n)$
represents the moduli functor $\underline{\M}^{\on{fr}}_U(r,0,1-n)$.
\end{prop}

\begin{proof} The proof amounts to using Theorem~\ref{beilinson5}
to show that the constructions of \cite{BW1, KKO} work in families.

We first use the  construction of \cite[Theorem~6.7]{KKO}
 to define a map $V\rightarrow \underline{\M}^{\on{fr}}_U(r,0,1-n)$.
Given a commutative $\mathbb C$-algebra $R$ and  $(b_1,b_2, j,i)\in V(R)$,   
define a complex
\begin{equation}\label{std in Weyl}
\bK: \; \theo_{U_R}(-1)^n \xrightarrow{A}\theo_{U_R}^n\oplus\theo_{U_R}^n
\oplus\theo_{U_R}^r
\xrightarrow{B} \theo_{U_R}(1)^n
\end{equation}
in $\qgr U_R$, where 
\bd
A= \begin{pmatrix}{\mathbf 1}_{n,n}\\ {\mathbf 0}_{n,n}\\
{\mathbf 0}_{r,n}\end{pmatrix}\cdot x + 
 \begin{pmatrix}{\mathbf 0}_{n,n}\\ {\mathbf 1}_{n,n}\\
{\mathbf 0}_{r,n}\end{pmatrix}\cdot y +
 \begin{pmatrix}b_1\\  b_2\\
j\end{pmatrix}\cdot z
\ed
and\ \
$
B = x\cdot ({\mathbf 0}_{n,n}\;\; {\mathbf 1}_{n,n}\;\; {\mathbf 0}_{n,r}) 
+y\cdot  (-{\mathbf 1}_{n,n}\;\; {\mathbf 0}_{n,n}\;\; {\mathbf 0}_{n,r}),
+z\cdot (-b_2 \;\; b_1\;\; i).
$

By \cite[Theorem~6.7]{KKO}, $\bK\otimes \mathbb C(p)$ is a torsion-free monad for
every $p\in\spec R$ and so $\bK$ is a torsion-free monad in
the sense of Definition~\ref{monad-defn}.
Thus Theorem~\ref{beilinson5} implies that
$\HH^0(\bK)$ is an $R$-flat family of torsion-free objects of $\qgr U_R$.
It follows that $\HH^0(\bK)|_{\pline} = \HH^0(\bK|_{\pline})$ and
so $\HH^0(\bK)$ comes equipped with a framing $\phi$.  We thus obtain
a $\on{GL}(n)$-invariant map $V\rightarrow 
\underline{\M}^{\on{fr}}_U(r,0,1-n)$.

We next show that  the map of functors
$\Psi : V/\on{GL}(n)\rightarrow \underline{\M}^{\on{fr}}_U(r,0,1-n)$
is an \'etale local isomorphism. 
Using Lemma~\ref{exts agree}, it follows from the proof of 
\cite[Theorem~6.7]{KKO} that $\Psi$ is injective, so it suffices to 
prove \'etale local surjectivity.

Observe that if $(\cE,\phi)$ 
is a framed torsion-free object of $\qgr U$,
then
 $\cE$ is $\mu$-semistable. To see this, suppose that
  $\cF\subset\cE$ de-semistabilizes 
 $\cE$. We may enlarge $\cF$ if necessary and  assume that $\cE/\cF$ is
torsion-free.  By Lemma~\ref{restriction of torsion-free}
$\cF|_{\pline} \subset\cE|_{\pline} = \theo_{\pline}^r$ and so
$$c_1(\cF) = c_1(\cF|_{\pline}) \leq 0 = c_1(\cE|_{\pline}) = 
c_1(\cE).$$ Thus, $\cF$ cannot de-semistabilize $\cE$.
The $\mu$-semistability of $\cE$ then implies that 
$\Hom(\cE,\cE|_{\pline})=0$, i.e. that framed torsion-free
objects of $\qgr U$ are rigid. By  Theorem~\ref{CBC}
the same result follows in families.

Lemma~\ref{vanishing cohomology
for ss modules} implies that if $(\cE,\phi)$ is an $R$-flat family of 
framed torsion-free
objects of $\qgr U_R$ then $\cE$ satisfies the Vanishing 
Condition~\ref{thesheaf}.  Theorem~\ref{beilinson5} then implies that
$\cE = \HH^0(\bK)$ for some torsion-free monad $\bK$.  After an \'etale
base change we may assume $\bK$ is a monad of the form \eqref{std in Weyl}.
In order to show that $\Psi$ is surjective,  it remains to show that 
this monad 
$\bK$ is actually isomorphic to a monad coming from
 a quadruple $(b_1,b_2, j,i)\in V(R)$. 
 This is immediate from the linear algebra of \cite[Theorem~6.7]{KKO}.

As in \cite[Theorem~6.7]{KKO}, the action of $\on{GL}(n)$ on $V$ is free,
and so the \'etale slice theorem implies that the map
 $V\rightarrow V/\!\!/\on{GL}(n)$ is a principle $\on{GL}(n)$-bundle.
But, since families of framed objects are rigid,  the fibre product
$V\times_{\underline{\M}^{\on{fr}}} U$ is easily checked to be a
principal $\on{GL}(n)$-bundle for any $\mathbb C$-scheme $U$. Thus
$V\rightarrow  \underline{\M}^{\on{fr}}_U(r,0,1-n)$ is also
a principal $\on{GL}(n)$-bundle and so $V/\!\!/\on{GL}(n)$ represents
$\underline{\M}^{\on{fr}}_U(r,0,1-n)$.
\end{proof}

\section{Poisson Structure  and Higgs Bundles}
\label{current groups} 
Throughout this section we assume that $k={\mathbb C}$ and 
 $S=\skl(E,\cL,\sigma)$; thus $E$ is an elliptic curve.
   Fix integer invariants   $\{r\geq 1,c_1,\chi\}$ and 
write $\M^s=\M^s_S(r,c_1,\chi)$ for the corresponding  smooth 
moduli space, as  defined by Theorem~\ref{projective moduli spaces}.

 The main aim of this 
 section is to prove that $\cM^s$
   is Poisson and  that the moduli space $(\mathbf P_S\smallsetminus E)^{[n]}$
   is symplectic, thereby proving Theorems~\ref{mainthm11} 
   and \ref{poisson-intro} from the introduction.
    We also take a first step towards relating our moduli spaces to integrable
    systems by  describing $\cM^s$ as the 
moduli for meromorphic Higgs bundles with structure group the centrally
extended current group on $E$ of Etingof-Frenkel \cite{EF}.
We achieve all these results by first observing that the  moduli spaces 
of Kronecker complexes on $S$ correspond to  the  moduli spaces for
 their restrictions 
to $E$, which we call  $\sigma$-Kronecker complexes.
This reduces the study of $\cM^s$ to a purely commutative question and 
the desired Poisson structure  then follows  from
 results in  \cite{Polishchuk}.
 
We begin by discussing  $\sigma$-Kronecker complexes.

\begin{defn}\label{sigma-monad-defn} Set $\baL = (\sigma^{-1})^*\cL$.
A {\em family of $\sigma$-Kronecker complexes parametrized by a 
$\mathbb C$-scheme $U$} is a complex on $U\times E$ of the form
\begin{equation}\label{sigma Kronecker complex}
\bK:\;\; V_{-1}\otimes \cL^* \xrightarrow{A} V_0\otimes \theo_E 
\xrightarrow{B} V_1\otimes\baL
\end{equation}
for some vector bundles $V_{-1}$, $V_0$ and $V_1$ on $U$.
By Lemma~\ref{shifts}, a family of $\sigma$  Kronecker complexes
\eqref{sigma Kronecker complex} is the same as a complex of the form 
\begin{equation}\label{sigma Kronecker complex2}
V_{-1}\otimes \theo_{B}(-1) \to V_0\otimes \theo_{B}
\to V_1\otimes \theo_{B}(1),\qquad \text{where } B=S/gS.
\end{equation}
 \end{defn}

We carry over the definitions of the invariants of a $\sigma$-Kronecker 
complex, $\sigma$-monad, etc. exactly as they
appear for Kronecker 
complexes in $\qgr S$ in 
Section~\ref{families of modules and Kronecker complexes}.  
In particular, the definition of (semi)stability for a $\sigma$-Kronecker
complex refers only to $\sigma$-Kronecker subcomplexes.  A {\em $\sigma$-monad} 
is a $\sigma$-Kronecker complex \eqref{sigma Kronecker complex} for
which the map $A$ is injective and the map $B$ is surjective.  
A $\sigma$-Kronecker complex (or, as before, Kronecker complex) 
is {\em framed}\label{framed-defn}
 if it comes equipped with
an isomorphism $V_0 = \theo^n$.

\begin{prop}\label{prop 7.2}
The moduli functors (stacks) for families of (normalized, semistable, stable, framed) Kronecker complexes in $\qgr S$ and 
families of (normalized, semistable, stable, framed) $\sigma$-Kronecker 
complexes on $E$ are isomorphic.
\end{prop}

\begin{proof} Use \eqref{sigma Kronecker complex2} 
as the definition of a  $\sigma$-Kronecker complex and consider the identities 
\begin{equation}\label{prop-7.2.1}
\Hom_{\qgr S}(\theo(n),\theo(n+i)) = \Hom_E(\theo(n)|_E,\theo(n+i)|_E)
\end{equation}
for $i=0, 1,2$.  When $i=1$, the identity gives the equivalence of 
diagrams of the form \eqref{typical Kronecker family} over $S$ and 
\eqref{sigma-monad-defn} over $B$. Then \eqref{prop-7.2.1} for $i=2$
 gives the equivalence 
for the complexes among those diagrams, while  \eqref{prop-7.2.1} for $i=0$ 
provides the equivalence of  maps between these complexes.
\end{proof}

\begin{corollary}\label{semistable gives sigma monad}
If $\bK$ is a family of geometrically semistable $\sigma$-Kronecker complexes,
then $\bK$ is a $\sigma$-monad.
\end{corollary}
\begin{proof}
Suppose that $\spec R$ is the parameter scheme for $\bK$.  
By Proposition~\ref{prop 7.2}, $\bK$ lifts to a family $\widetilde{\bK}$ 
of geometrically semistable Kronecker complexes in $\qgr S_R$.  
For every $p\in\spec R$, we have 
$\bK\otimes \mathbb C(p) = \widetilde{\bK}\otimes \mathbb C(p)|_E$. By
Corollary~\ref{monad families}, $\widetilde{\bK}\otimes \mathbb C(p)$ is therefore a 
torsion-free monad in $\qgr S_{\mathbb C(p)}$.  Surjectivity of $B_{\bK\otimes
\mathbb C(p)}$
then follows from that of $B_{\widetilde{\bK}\otimes \mathbb C(p)}$.
Moreover, $\on{coker}(A_{\widetilde{\bK}\otimes \mathbb C(p)})$
is torsion-free and so Lemma~\ref{restriction of torsion-free}
 implies that $A_{\bK\otimes \mathbb C(p)}$ is injective. 
\end{proof}

\subsection{Construction of Poisson Structure}
\label{construction of poisson structure}

The aim of this subsection is to prove the following result, which 
describes the Poisson structures of our moduli spaces. 

\begin{thm}\label{poisson structure} Let $S=\skl(E,\cL,\sigma)$ be the Sklyanin
algebra defined over $\mathbb C$ and pick integers $\{r\geq 1,c_1,\chi\}$.
 Then:
\begin{enumerate}
\item 
$\M_S^s(r,c_1,\chi)$ admits a Poisson structure. 
\item  Fix  a vector bundle
$H$ on $E$ and let $\M_H$ denote the smooth locus of
the subvariety of $\M_S^s(r,c_1,\chi)$ that 
parametrizes those $\cE$ for which $\cE|_E\cong H$. Then $\M_H$ 
is an open subset of a symplectic leaf of the Poisson space.
\end{enumerate}
\end{thm}
\begin{remark}\label{proving 1.4}
This result proves Theorem~\ref{poisson-intro} from the introduction.
 Combined with Corollary~\ref{monad6}(3) it also proves Theorem~\ref{mainthm11}.
\end{remark}

We begin by explaining the strategy for the proof.
By  Remark~\ref{other-c1} we can and will  assume 
 that $-r<c_1\leq 0$ (this may replace
$H$ by some shift $H(m)$). 
By Proposition~\ref{prop 7.2} it then suffices
 to prove the existence of a
Poisson structure on the corresponding moduli space $\Sigma$ of
 geometrically stable $\sigma$-monads.  
As we will explain, this essentially follows from work of Polishchuk
 \cite{Polishchuk} (see Proposition~\ref{Polish Poiss}),
  although we will first need to set up the appropriate framework to apply his
  results. 

Fix a positive integer $n$ and  let $\M''$ \label{M-dash-defn}
denote the moduli stack
parametrizing data $(\cE_i,\phi_i)$ consisting of $n$-tuples of vector bundles
$\cE_i$  on $E$ and maps $\phi_i: \cE_{i+1}\rightarrow \cE_i$.   There is a
closed substack $\M'$ of $\M''$ parametrizing those $(\cE_{i+1}, \phi_i)$ for
which $\phi_i\circ\phi_{i+1} = 0$ for all $i$; in other words, $\M'$
parametrizes {\em complexes} of vector bundles. If  $(\cF_i)$ is an $n$-tuple
of $E$-vector bundles, let   $\M = \M(\cF_1, \dots, \cF_n)$  denote the 
locally closed substack of $\M'$  parametrizing data $(\cE_i,\phi_i)$ for which
$\cE_i\cong \cF_i$  for all $i$.  As we will show, it follows from
\cite{Polishchuk} that there is  Poisson structure on an open substack of $\M'$
and this  restricts to a Poisson structure on the smooth locus of $\M$.  Since 
$\Sigma$ is  an open substack  of  $\M(\baL\,^{d_1}, \theo_E^n,
(\cL^*)^{d_{-1}})$, in the notation of  \eqref{lin-alg}, this  will  give the
desired  Poisson structure on  $\M^s(r,c_1,\chi)$.

Our first task is to construct   a map $\Psi''$ on (co)tangent spaces of $\M''$,
 for which  we use
following concrete description of those spaces.
Let $(\cE_i,\phi_i)$ be data defining a point of $\M''$.
The discussion of \cite[Section~1]{Polishchuk}, 
shows that the tangent space to
$\M''$ at $(\cE_i,\phi_i)$ is the first hypercohomology ${\mathbf H}^1(C)$ of a complex 
$C = C(\cE_i,\phi_i)$
defined as follows: $C$ is concentrated in degrees $0$ and $1$ 
with  $C^0 = \bigoplus \uEnd(\cE_i)$ and  
$C^1 = \bigoplus \uHom(\cE_{i+1},\cE_i)$.
  The differential $d$ sends
$(e_i)\in C^0$ to $(e_i\phi_i - \phi_i e_{i+1})\in C^1$.

We will actually be interested in the open
substack $\M''_0$\label{M-dash-oh}
of the smooth locus of 
$\M''$ that consists of points $(\cE_i,\phi_i)$ for which
the hypercohomology of $C=C(\cE_i,\phi_i)$ satisfies
\begin{equation}\label{M'' cond} 
{\mathbf H}^2(C) = 0\; \text{and}\; {\mathbf H}^0(C) = {\mathbb C}.
\end{equation}
 This substack parametrizes objects that have only scalar endomorphisms.
The tangent sheaf of $\M''_0$ is a locally free sheaf, the 
fibre of which over $(\cE_i,\phi_i)$ is 
${\mathbf H}^1\big(C(\cE_i,\phi_i)\big)$.   Let $\M'_0$ denote
the smooth locus of $\M'\cap \M''_0$ and similarly let $\M_0$ denote the
smooth locus of $\M\cap \M''_0$. 

We next construct a map $\Psi'': T^*\M''|_{\M'_0}\rightarrow T\M''|_{\M'_0}$. 
Let $C =C(\cE_i,\phi_i)$.
The dual complex $C^{\vee}[-1]$ of  $C$  
 is  the complex
 $\bigoplus\uHom(\cE_i,\cE_{i+1})\xrightarrow{d^*} \bigoplus\uEnd(\cE_i)$,
again concentrated in degrees $0$ and $1$, with
 differential  $d^* : (f_i) \mapsto (\phi_i f_i- f_{i-1}\phi_{i-1})$.
By Serre duality, $T^*_{(\cE_i,\phi_i)}\M''\cong {\mathbf H}^1(C^{\vee}[-1])$.
We define maps
$\psi_k: (C^{\vee}[-1])^k\rightarrow C^k$ by  
$$\psi_0(f_i)= \big( (-1)^{i+1}(\phi_i f_i - f_{i-1}\phi_{i-1})\big)\qquad
\text{and}\qquad \psi_1=0.$$
We   define a second pair of  maps  
$\overline{\psi}_k: C^{\vee}[-1]^k\rightarrow C^k$
by   $$\overline{\psi}_0=0 \qquad
\text{and}\qquad\overline{\psi}_1: (e_i)\mapsto 
 \big((-1)^{i+1}(e_i\phi_i + \phi_ie_{i+1})\big).$$

 The reason for working with $\cM_0'$ rather than $\cM''_0$, is that these
 maps define maps of complexes. Indeed, an elementary calculation, using the
 fact that $\phi_i\phi_{i+1}=0$, gives the following result. (See also 
  \cite[Theorem~2.1]{Polishchuk}.)

\begin{lemma}\label{psi and homotopy}
If $(\cE_i,\phi_i)$ lies in $\M'_0$, then:
\begin{enumerate}
\item $\psi=(\psi_0,\psi_1)$ defines a morphism of complexes
$\psi: C^{\vee}[-1]\rightarrow C$.
\item $\overline{\psi}=(\overline{\psi}_0,\overline{\psi}_1)$
 defines a morphism of complexes
$\overline{\psi}: C^{\vee}[-1]\rightarrow C$.
\item The map $h: \bigoplus \uEnd(\cE_i)\rightarrow \bigoplus \uEnd(\cE_i)$ defined 
by $h(e_i)= \big((-1)^ie_i\big)$ defines a homotopy between $\psi$ and 
$\overline{\psi}$. Hence ${\mathbf H}^1(\psi) =
 {\mathbf H}^1(\overline{\psi})$.\qed
\end{enumerate}
\end{lemma}

As in \cite[Section~2]{Polishchuk}, the maps $\psi$ 
and $\overline{\psi}$ globalize to give a map
$\Psi'': T^*\M''|_{\M'_0}\rightarrow T\M''|_{\M'_0}$\label{Phi-defn}
 that, in the fibres over $(\cE_i,\phi_i)$, is exactly the map
 ${\mathbf H}^1(\psi)$. The proof of 
\cite[Theorem~2.1]{Polishchuk} also shows that $\Psi''$ is skew-symmetric.

We next want to show that $\Psi''$ factors through maps 
\bd
\Psi': T^*\M'_0\rightarrow T\M'_0 \qquad\text{and}\qquad
\Psi: T^*\M_0\rightarrow T\M_0.
\ed
Let $(\cE_i,\phi_i)$ define a point of $\M'_0$ and $C=C(\cE_i,\phi_i)$.
Define a complex $B$ by setting $B^0 = \bigoplus \uEnd(\cE_i)=C^0$
and $B^1 = \bigoplus\uHom(\cE_{i+2},\cE_i)$, with the zero  differential.
A simple calculation shows that there is a map of complexes
$\Xi : C\rightarrow B$, where $\Xi^0$ is the identity and $\Xi^1$
maps $(\psi_i)\in C^1$ to 
$(\phi_i\psi_{i+1} + \psi_i\phi_{i+1})$. 

\begin{lemma}\label{tangents}
\begin{enumerate}
\item $T_{(\cE_i,\phi_i)}\M'_0 = \ker\big[{\mathbf H}^1(C)\rightarrow 
{\mathbf H}^1(B^1)\big]$.
\item  If $(\cE_i,\phi_i)$ determines a point of $\M_0$, then
$T_{(\cE_i,\phi_i)}\M_0 = \ker(\mathbf H^1(\Xi))$. 
\item The  cotangent spaces $T^*_{(\cE_i,\phi_i)}\M'_0$ and 
$T^*_{(\cE_i,\phi_i)}\M_0$
are the cokernels of the respective dual maps.
\end{enumerate}
\end{lemma}
\begin{proof}
Since the differential in $B$ is zero, 
\begin{equation}\label{hypercoh of B}
{\mathbf H}^1(B) = {\mathbf H}^1(B^0)\oplus {\mathbf H}^1(B^1) 
 =  \left[\bigoplus \HH^1(\uEnd(\cE_i))\right] \oplus 
\left[\bigoplus \Hom(\cE_{i+2},\cE_i)\right].
\end{equation}
A cocycle calculation shows that
a class  $c$ corresponding to a first-order deformation 
$(\widetilde{\cE_i},\widetilde{\phi}_i)\in T_{(\cE_i,\phi_i)}\M'_0$ 
lies in $\on{Ker}\big[{\mathbf H}^1(C)\rightarrow\mathbf H^1(B^1)\big]$
 if and only if $\widetilde{\phi}_i\circ\widetilde{\phi}_{i+1}=0$ for all 
$i$. This proves part (1).
  Similarly, a calculation shows
that
$c\in \ker(\mathbf H^1(\Xi))$
if and only if both $\widetilde{\phi}_i\circ\widetilde{\phi}_{i+1}=0$ for all
$i$ and every $\widetilde{\cE}_i$ is a trivial deformation of $\cE_i$, proving 
 (2).
Part (3) then follows by Serre duality.
\end{proof}

\begin{corollary}\label{factor-phi}
$\Psi''|_{\M_0}$
 factors through a map 
$\Psi: T^*\M_0\rightarrow T\M_0$. Similarly,
$\Psi''|_{\M'_0}$ factors through a map
 $\Psi': T^*\M'_0\rightarrow T\M'_0$.   
\end{corollary}
\begin{proof} We will just prove the first assertion, 
since the same proof works in both cases.

By Lemma~\ref{tangents}, to prove that $\Psi''|_{\M'_0}$ 
factors through $T\M_0$ it suffices to show that the composite map 
\begin{equation}\label{Poisson comp}
{\mathbf H}^1(C^{\vee}[-1]) \xrightarrow{\Psi''} {\mathbf H}^1(C)
\xrightarrow{\mathbf H^1(\Xi)}  {\mathbf H}^1(B)
\end{equation}
is zero. 
To prove this we combine
Lemma~\ref{psi and homotopy}(3) with two simple observations:
 first, since $\psi_1 $ is defined
to be zero,
 the composite of \eqref{Poisson comp}
followed by projection onto
 $\bigoplus \Hom(\cE_{i+2},\cE_{2})=\mathbf H^1(B^1)$ is zero. 
Secondly, the homotopic map of complexes $\overline{\psi}$ is
zero in cohomological degree $0$ and so
the composite of \eqref{Poisson comp} followed by
projection onto $\HH^1(\uEnd(\cE_i))$ is zero.  It therefore follows from 
\eqref{hypercoh of B} that \eqref{Poisson comp} is zero.  

To obtain the  factorization through $T^*\M'_0$ and $T^*\M_0$ we use the
description from Lemma~\ref{tangents} in terms of the dual map
${\mathbf H}^1\big((B^1)^{\vee}[-1]\big)\rightarrow 
{\mathbf H}^1(C^{\vee}[-1])$.
Taking the dual of \eqref{Poisson comp} and using that 
$(\Psi'')^{\vee}=-\Psi''$, it follows that 
the factorization through cotangent
spaces is equivalent to the vanishing of the dual to
 \eqref{Poisson comp}.
\end{proof}

We now have the following consequence of Polishchuk's work.

\begin{prop}\label{Polish Poiss}
The map $\Psi'$ is a Poisson structure on  $\M_0'$ while 
 $\Psi$ is a Poisson structure on $\M_0$.
\end{prop}

\begin{proof} We only  prove the first assertion since the second
one follows from the same argument.
This is really a special case of 
\cite[Theorem~6.1]{Polishchuk}, but we need to check that 
this particular result does fit into  Polishchuk's framework. 

We start by reinterpreting the
 map $\psi: C^\vee[-1]\to C$ from Lemma~\ref{psi and homotopy}.   Fix 
a component of $\M''$ and let $n_i = \on{rk}(\cE_i)$ be the corresponding
ranks of the vector bundles. Let $G = \prod \on{GL}_{n_i}({\mathbb C})$ and
$V = \bigoplus \Hom({\mathbb C}^{n_{i+1}}, {\mathbb C}^{n_i})$
with the natural $G$-action.   Then $\M''$ 
can and will be identified with
an open substack of the moduli stack of pairs $(P,s)$ consisting of a 
principal $G$-bundle $P$ and a section $s$ of the associated bundle
$P\times_G V$. Under this identification, the 
pair $(P,s)$ maps to the data $(\cE_i,\phi_i)$
where the sequence of vector bundles $\cE_i$ is given by 
$P\times_G \bigoplus {\mathbb C}^{n_i}$
 and the maps $\phi_i$ are determined by $s$. 

Next, set ${\mathfrak g}_{n_i} = \on{Lie}(\on{GL}_{n_i}({\mathbb C}))$ and
${\mathfrak g} = \on{Lie}(G)$.  Let
 ${\mathfrak t}_{n_i}\in 
\on{Sym}^2({\mathfrak g}_{n_i})^{{\mathfrak g}_{n_i}}$
  denote the  dual of the trace
pairing and set 
${\mathfrak t}=\bigoplus (-1)^{i+1}{\mathfrak t}_{n_i}\in
\on{Sym}^2({\mathfrak g})^{\mathfrak g}$.
By abuse of notation we also use ${\mathfrak t}$ to denote the induced map
${\mathfrak g}^*\rightarrow {\mathfrak g}$.  Associated  to each
$v\in V$ one has the map $d_v: {\mathfrak g} \rightarrow V$ 
defined by $d_v(X) = X\cdot v$ and this induces the map
${\mathfrak t}\circ d_v^*: V^*\rightarrow {\mathfrak g}$.
We therefore obtain 
a map $\psi_0=P\times_G({\mathfrak t}\circ d_v^*):
 P\times_G V^*\rightarrow P\times_G {\mathfrak g}$ 
associated to a pair $(P,s)$. As in \cite{Polishchuk}, 
this is precisely  the map $\psi_0:(C^{\vee}[-1])^0 \rightarrow C^0$
from Lemma~\ref{psi and homotopy}.  
As before, define $\psi_1: (C^{\vee}[-1])^1 \rightarrow C^1$ 
to be the zero map. 

\begin{sublemma}\label{polish-cnd} If $v=(\cE_i,\phi_i)\in \cM'$, then
$d_v\circ {\mathfrak t}\circ d_v^*= 0$.
\end{sublemma}

\begin{proof}[Proof of the sublemma] 
Although the proof is an elementary computation, we give it since
the result is the basic hypothesis for
 \cite[Theorem~6.1]{Polishchuk} (see \cite[p.699]{Polishchuk}).
 Identify 
$V^* =\bigoplus \on{Hom}(\mathbb C^{n_{i}},\mathbb C^{n_{i+1}})$
by $\langle (\alpha_i),(\beta_j)\rangle = \sum_i\on{Tr}(\alpha_i\beta_i)$
for $(\alpha_i)\in V$ and $(\beta_j)\in V^*$.
Fix $(\alpha_i)\in V$ corresponding to a point in $\cM'$. 
Then $d_{(\alpha_i)}(X_j) = (X_i\alpha_i-\alpha_iX_{i+1})$ for 
$(X_j)\in \mathfrak g$.
Thus, $d^*_{(\alpha_i)}$ is defined by 
$$\begin{array}{rl}
&d^*_{(\alpha_i)}(\beta_j)(X_k) =(\beta_j)(d_{(\alpha_i)})(X_k)=
 \langle d_{(\alpha_i)}(X_k)\,,\,(\beta_j)\rangle  \\ 
 \noalign{\vskip 5pt}
&\quad=\sum_i \on{Tr}(X_i\alpha_i\beta_i) 
-\sum_i \on{Tr}(\alpha_{i+1}X_{i+1}\beta_{i}) 
=  \sum_i \on{Tr}(X_i(\alpha_i\beta_i-\beta_{i-1}\alpha_{i-1})).
 \end{array}
 $$ 
This implies that $\mathfrak t \circ d^*_{(\alpha_i)}(\beta_j)
= (-1)^{i+1}(\alpha_i\beta_i+\beta_{i-1}\alpha_{i-1}).$ Finally,
$$d_{(\alpha_i)}\circ\mathfrak t \circ d^*_{(\alpha_i)}(\beta_j)
= (-1)^{i+1}(\alpha_i\beta_i\alpha_i +\beta_{i-1}\alpha_{i-1}\alpha_i
-\alpha_i\alpha_{i+1}\beta_{i+1}-\alpha_i\beta_i\alpha_i).$$
 But  $\alpha_i\alpha_{i+1}=0$      by the definition of $\cM'$. Thus
 $d_{(\alpha_i)}\circ\mathfrak t \circ d^*_{(\alpha_i)}=0$.
\end{proof}

We return to the proof of  the proposition.  As in \cite[p.699]{Polishchuk},
skew-symmetry of the Poisson bracket 
follows from the sublemma, so it remains to prove the Jacobi
identity.  This is equivalent to the equation
\begin{equation}\label{Jacobi}
\Psi'(\omega_1)\cdot \langle\Psi'(\omega_2),\omega_3\rangle
 -\langle [\Psi'(\omega_1),\Psi'(\omega_2)],\omega_3\rangle + 
 \text{cp}(1,2,3) =0,
\end{equation}
where  the $\omega_i$ are local $1$-forms on $\M'_0$
and $\text{cp}(1,2,3)$  denotes cyclic permutations of the indices $\{1,2,3\}$.
By Corollary~\ref{factor-phi}, we may replace 
 $\Psi'$ by $\Psi''|_{\M'_0}$ in \eqref{Jacobi}. 
The fact that \eqref{Jacobi} holds is now
 \cite[Theorem~6.1]{Polishchuk}, except that we have to work 
 with $\cM'$ rather than $\cM''$ in order to ensure that the hypothesis of
 Sublemma~\ref{polish-cnd}
 is valid. However, the proof of \cite[Theorem~6.1]{Polishchuk}
 can be used without change to prove \eqref{Jacobi} and hence the proposition.
\end{proof}

It is now easy to complete the proof of Theorem~\ref{poisson structure}.

\begin{proof}[Proof of Theorem~\ref{poisson structure}]
By Remark~\ref{other-c1} we may assume that $-r<c_1\leq 0$.
By Propositions~\ref{semistable module from complex} and 
\ref{prop 7.2}, it then suffices to
prove the result for $\sigma$-monads rather than
modules. By Lemma~\ref{shifts} and \eqref{lin-alg}, such a  $\sigma$-monad
has
the form
\begin{equation}\label{more-data}
\cE: \;\; \cE_3\xrightarrow{\phi_2} \cE_2\xrightarrow{\phi_1} \cE_1, \quad \text{where}\
\cE_3\cong (\cL^*)^{d_1},\ \cE_2\cong \theo^n\ \text{and}\
\cE_1\cong \baL\,^{d_1}.
\end{equation}
 Thus
$\M=\M\big(\baL\,^{d_1},\theo_E^n,(\cL^*)^{d_{-1}}\big)$
 is the moduli stack of $\sigma$-Kronecker complexes with the specified
invariants and it has an open substack $\cM^s$ that
is the moduli stack parametrizing the corresponding
geometrically stable $\sigma$-monads.
   By Proposition~\ref{Polish Poiss}, in order
to prove  part (1) of the theorem, it suffices to prove:
 \begin{enumerate}
\item[(i)]  $\M^s$ lies in the open
substack of the smooth locus of $\M''$ consisting
of points satisfying \eqref{M'' cond},
\item[(ii)] $\M^s$ lies in the
smooth locus of $\M$ and
\item[(iii)] the
Poisson structure on $\M$ induces one on $\M^s_S(r,c_1,\chi)$.
\end{enumerate}

We first prove (i).
 Let $(\cE_i,\phi_i)$ be the datum of a geometrically stable
$\sigma$-monad $\cE$ on $E$ and set $C=C(\cE_i,\phi_i)$.
If $\on{Bun}$ denotes the moduli stack of triples of vector bundles on $E$
then
$\on{Bun}$ is smooth and
there is a forgetful morphism $\M''\rightarrow \on{Bun}$ sending
$(\cF_i,\phi_i)$ to $(\cF_i)$.  For the given   datum $(\cE_i,\phi_i)$,
we have $\Ext^1(\cE_{i+1},\cE_i)=0$ for $i=1,2$.  Thus
Theorem~\ref{CBC}(3) implies that   $\M''$
is a vector bundle over $\on{Bun}$ in a neighborhood of $(\cE_i)$ and  so
$\M''$ is smooth at $(\cE_i,\phi_i)$.
 By Serre duality and \eqref{more-data},
$${\mathbf H}^2(C)^* = {\mathbf H}^0(C^{\vee}[-1])
=\on{ker}\left[\bigoplus\Hom(\cE_i,\cE_{i+1})\rightarrow
 \bigoplus \End(\cE_i)\right]=0.$$
 On the other hand,
${\mathbf H}^0(C)=
\on{ker}\left[\bigoplus\End(\cE_i)
\xrightarrow{d} \bigoplus\Hom(\cE_{i+1},\cE_i)\right]$
 is the endomorphism ring $\on{End}(\cE)$ of the $\sigma$-monad $\cE$.
Since  $\cE$ is   stable,
Lemma~\ref{stable is simple} and Proposition~\ref{prop 7.2}
imply that  $\on{End}(\cE)=\mathbb C$. Thus (i) holds.

By Remark~\ref{stack-quotient} and Proposition~\ref{prop 7.2},
  $\M^s$ is isomorphic to
  the stack-theoretic quotient $[\N^s/\on{GL}(H)]$
and so it is smooth at $(\cE_i,\phi_i)$, by
Remark~\ref{smooth family-remark}. Thus (ii) holds.

Moreover, $\M^s$ comes equipped with a map
$\alpha:\M^s\rightarrow \M^s(r,c_1,\chi)$.
By the discussion at the beginning of Section~\ref{moduli-smooth},
 the \'etale slice theorem implies that
$\alpha$ is  \'etale locally of the form
$U\times B\cs\rightarrow U$ and so it induces isomorphisms of
tangent and cotangent bundles under pullbacks.  It follows that $\Psi$
also defines a Poisson structure on $\M^s(r,c_1,\chi)$, proving (iii).
Thus part (1) of the theorem is true.

(2) The strategy of the proof is the following. 
 Let $\M_H$ denote the smooth locus
of the locally closed subscheme of $\M_S^s(r,c_1,\chi)$ that parametrizes
those $\cE$ in $\qgr S$ for which $\cE|_E\cong H$.
 For such a module $\cE$, it is awkward to directly relate $T_{\cE}\M_H$
to the map $\Psi$.
  Instead, we construct a complex $D=D(\cE)$ together with a map
$D\rightarrow C$ and use this to identify $\Psi'_{\cE}$ with a map
${\mathbf H}^1(D^{\vee}[-1])\rightarrow {\mathbf H}^1(D)$.  It will then be
easy to prove that $T_{\cE}\M_H$ is exactly the image of this map and
it will  follow that $\Psi'$ induces a nondegenerate map 
$T^*_{\cE}\M_H\rightarrow T_{\cE}\M_H$.  
  
Thus, fix a vector bundle $H$ on $E$ and suppose that 
$\cE$ is a 
$\sigma$-monad corresponding to a point of $\M_H$, as in \eqref{more-data}. Then
$\phi_2$ makes $\cE_3$ a subbundle of $\cE_2$ and $\phi_1$ is
surjective.  
We write 
$\uEnd(\cE)$ for the sheaf of
endomorphisms of the $\sigma$-monad $\cE$.

  Consider the complex
\bd
D: \;\; \bigoplus\uEnd(\cE_i)\xrightarrow{d}\bigoplus \uHom(\cE_{i+1},\cE_i)
\xrightarrow{\Xi^1}\bigoplus\uHom(\cE_{i+2},\cE_i),
\ed
where  $\bigoplus\uEnd(\cE_i)$ lies in degree zero, 
$d$ is the differential in the complex $C$ and $\Xi^1$ 
 is the map defined before Lemma~\ref{tangents}.
There is a natural map of complexes $D\rightarrow C$ given by projection.
 As in \cite[Lemma~3.1]{Polishchuk},
a routine computation shows:

\begin{fact}\label{D cohom}
The cohomology of $D$ satisfies $\HH^i(D)=0$ for $i\neq 0$ 
and $\HH^0(D) =\uEnd(\cE)$.
\end{fact}

A standard \v{C}ech calculation then gives
$T_{(\cE_i,\phi_i)}\M'_0 = \HH^1(\uEnd(\cE))$.  It follows that there is a
canonical isomorphism $T_{(\cE_i,\phi_i)}\M'_0
\xrightarrow{\cong} {\mathbf H}^1(D)$.

The dual complex $D^{\vee}[-1]$ to $D$ has the form
\bd
D^{\vee}[-1]:\;\; \bigoplus\uHom(\cE_i,\cE_{i+2})\rightarrow
\bigoplus\uHom(\cE_i,\cE_{i+1})\rightarrow \bigoplus\uEnd(\cE_i).
\ed
We define a map
$\psi_D: D^{\vee}[-1]\rightarrow D$ by taking $(\psi_D)_i = \psi_i$
for $i=0,1$ (where $\psi$ is the map $C^{\vee}[-1]\rightarrow C$
defined above) and $(\psi_D)_{-1}=0$.  A straightforward computation
shows that
$\psi_D$ is a homomorphism of complexes.
It is clear from the construction that
$\psi: C^{\vee}[-1]\rightarrow C$ factors through $\psi_D$ and it follows
that
\bd
{\mathbf H}^1(\psi_D) = \left(\Psi': T^*_{\cE}\M_0'
\longrightarrow  T_{\cE}\M_0'\right).
\ed

Thus, in order to show that $\Psi'$ factors through 
$T_{\cE}\M_H\subset T_{\cE}\M'_0$ it suffices to show that
$
T_{\cE}\M_H = 
\on{ker}\left[{\mathbf H}^1(D)\rightarrow
 {\mathbf H}^1(\on{Cone}(\psi_D))\right],
$
where $\on{Cone}(\psi_D)$ is the complex
\bd
\bigoplus\uHom(\cE_i,\cE_{i+2})\rightarrow \bigoplus\uHom(\cE_i,\cE_{i+1})
\xrightarrow{(-d^*,\psi_0)} \left[\bigoplus\uEnd(\cE_i)\right]\oplus
\uEnd(\cE),
\ed
with the right hand term in degree $0$.
By construction, $ \on{Cone}(\psi_D)$ is exact except at the right hand term.
Define a map
$
\on{Cone}(\psi_D) \longrightarrow 
\left[\bigoplus\uEnd(\cE_i)\right]\oplus\uEnd(H)
$
by  
\bd
(A_3,A_2,A_1,F) \in\left[\bigoplus\uEnd(\cE_i)\right]\oplus
\uEnd(\cE)\mapsto (A_1+F_1, -A_2+F, A_3+F_3, 
{\mathcal H}^0(F)),
\ed
 where $F_i$ denotes the $i$th graded
component of $F\in \uEnd(\cE)$.  A routine check shows that
this map is a quasi-isomorphism.

It follows that the sequence
\begin{equation}\label{psi'}
{\mathbf H}^1(D^{\vee}[-1])\xrightarrow{\Psi'}
{\mathbf H}^1(D)\xrightarrow{\alpha}
\left[\bigoplus\HH^1(\uEnd(\cE_i))\right]\oplus\HH^1(\uEnd(H))
\end{equation}
is exact.  Here, the map $\alpha$ takes the class of
a first-order deformation
of $\cE$ as a complex to the classes of the
 associated first-order deformations of the $\cE_i$ and of the middle cohomology of
$\cE$.  Therefore, $\on{Im}(\Psi')$
 consists exactly of those first-order deformations
of $(\cE_i,\phi_i)$ such that the induced deformations of the
$\cE_i$ and $H$ are trivial.  This is exactly
$T_{(\cE_i,\phi_i)}\M_H$ and so $\Phi$ does factor through $T_\cE\cM_H$.

The skew-symmetry of $\Psi'$
shows that the Serre dual to \eqref{psi'} is an exact sequence
\bd
{\mathbf H}^1((\on{Cone}(\psi)^{\vee}[-1]) 
\rightarrow {\mathbf H}^1(D^{\vee}[-1]) \xrightarrow{\Psi'} {\mathbf H}^1(D).
\ed  
It follows that $\Psi'$ factors through 
\bd
T^*_{(\cE_i,\phi_i)}\M_H = \on{coker}\left[{\mathbf H}^1((\on{Cone}(\psi)^{\vee}[-1]) 
\rightarrow {\mathbf H}^1(D^{\vee}[-1])\right].
\ed
In particular, $\Psi'$ induces a surjective map 
$T^*_{(\cE_i,\phi_i)}\M_H\rightarrow T_{(\cE_i,\phi_i)}\M_H$, which must
therefore be an isomorphism.  This completes the proof.
\end{proof}

\subsection{From Equivariant Bundles to Equivariant Higgs Bundles}
\label{equivariant higgs section}
The treatment in this section will be informal, since it takes us rather
far afield from the main results of the paper. However, inspired
in part by constructions of integrable particle systems using 
current algebras (see \cite{GN, Nekrasov}), we want to
indicate the relationship between the earlier results of this paper and
a theory of Higgs bundles taking values in a centrally extended
current group  on the elliptic curve \cite{EF}.

Let
$S=\skl(E,\cL,\sigma)$, for $k=\mathbb C$ and (without loss of generality) 
 fix invariants $ r,\dots, d_1$  by
Notation~\ref{invariants-notation}. 
 Fix also 
a vector bundle $H$ on $E$ of rank $r$ and with determinant
$\cL^{\otimes d_{-1}}\otimes {\baL}^{\otimes -d_{1}}$. 
Fix the $C^\infty$ complex 
vector bundle $\cV$ on $E\times{\boldp}^1$, defined as the underlying
$C^{\infty}$ bundle of 
\bd 
\big((\cL^*)^{d_{-1}}\boxtimes \theo_{\pline}\big) 
\oplus \big(H\boxtimes \theo_{\pline}(-1)\big)
\oplus \big(\baL\,^{d_{1}}\boxtimes \theo_{\pline}(-2)\big)
\ \in\ \coh(E\times
\boldp^1).
\ed
Let $\cs$ act on $E\times\pline$ via the trivial action on $E$ and the usual
scaling action on $\pline$.
We give $\cV$ a $\cs$-equivariant structure as follows: 
For $-2\leq i\leq 0$, give $\theo(i)$ the natural $\cs$ structure
in which $\cs$ acts with weight
zero on the fibre at infinity.  Then take the natural $\cs$ action 
this induces on the components of $\cV$.

We regard $\cV$ as a family of equivariant $C^{\infty}$ vector bundles on $E$
parametrized by ${\boldp}^1$ and  will  write $\cV_z$ for  the restriction of
$\cV$ to $E\times\{z\}$. We will also want to think of $\cV$ as a sheaf
living on  ${\boldp}^1$, in which case we write it as  $B_\cV$. Thus, for an
open subset  $U\subset \boldp^1$, the sections  $B_\cV(U)$  consist of the
$C^{\infty}$ sections of $\cV|_{E\times U}$. Let $L_{\cV}$ denote the sheaf 
of (infinite-dimensional)  Lie algebras on $\boldp^1$ whose
sections $L_\cV(U)$ are defined to be  the space of $C^{\infty}$ 
sections of the vector bundle
$\mbox{End}(\cV)|_{E\times U}$.

As in \cite{EF}, the sheaf $L_\cV$ has a sheaf of central 
extensions $\widehat{L}_{\cV}$ that, as an extension 
of sheaves,   has the form
$$
0\rightarrow C^{\infty}(\pline,{\mathbb C}) \rightarrow \widehat{L}_{\cV} 
\rightarrow L_{\cV}\rightarrow 0.
$$ In the
topologically trivial case  $\widehat{L}_{\cV}$  is the central extension 
determined by the cocycle
$ 
\Omega(X,Y) = \int_E\eta\wedge\langle X, \overline{\partial}_t Y\rangle,
$
where $\eta$ is a nonzero holomorphic 1-form on $E$ and 
$\langle \, , \, \rangle$ is the trace form.
For general $\cV$, one needs to
 replace the symbol $\overline{\partial}_t$ by any partial 
$\overline{\partial}$-operator in the $E$ direction.
 For example, this may be done by choosing a 
$\overline{\partial}$-operator for $\cV$ on $E\times \pline$, say 
$\overline{D}: \cV\rightarrow \cV\otimes (T^*)^{0,1}(E\times\pline)$, and
using the splitting $T^*(E\times\pline) \cong p_E^*T^*E 
\oplus p_{\pline}^*T^*\pline$ to project to an 
operator $\overline{\partial}: \cV\rightarrow \cV\otimes p_E^*(T^*)^{0,1}(E)$. 
A useful intuition is to think of $\widehat{L}_{\cV}$  as
the sheaf of sections of a
$C^{\infty}$-bundle of centrally extended current algebras.

Write $E= {\mathbb C}/\Lambda$ and let $t$ denote the coordinate on
${\mathbb C}$ as well as the induced local coordinate on $E$.
Consider a first-order differential operator $D$ on $\cV$ that locally,
 say on an
open subset  $E\times U\subset E\times\boldp^1$ with coordinates $t$ on
$E$ and $z$ on $\boldp^1$,  has the form
 $D =\lambda\cdot d\overline{t}
\cdot \partial/\partial\overline{t} + A\cdot d\overline{t}$,  where $A$ is a
matrix of endomorphisms of $\cV$ with coefficients that are $C^{\infty}$
functions on  $E\times U$ and $\lambda\in{\mathbb C}$. 
Let 
$\widehat{L}_{\cV}^*$ denote the space of such operators, regarded as a
sheaf over  $\mathbf P^1$.
  This bundle is essentially (in a
sense we will not try to make precise) the smooth part of the fibrewise dual of
$\widehat{L}_{\cV}$ (see \cite{EF} for a discussion). We will therefore 
regard $\widehat{L}_{\cV}^*$
with the action of the sheaf of groups of automorphisms of
$\cV$  as the bundle
of coadjoint representations of the centrally extended relative gauge group of
$\cV$ whose associated centrally extended endomorphism algebra is
$\widehat{L}_{\cV}$.  Finally,
 let $(\widehat{L}_{\cV}^*)_\lambda$ denote the bundle 
of hyperplanes in $\widehat{L}_{\cV}^*$ 
consisting of those operators $D$ having the given coefficient $\lambda$ of
$d\overline{t}\cdot\partial/\partial\overline{t}$.

\begin{defn}
A {\em $\overline{\partial}$-operator $\overline{\partial}_B$} on $B_{\cV}$ 
is an operator\label{partial-defn}
$
\overline{\partial}_B: B_{\cV} \to B_{\cV}\otimes 
(T^*)^{0,1}\boldp^1
$
on $\boldp^1$ that can be represented locally on $\boldp^1$ 
in the form 
$\overline{\partial}_B= d\overline{z}\cdot \partial/\partial\overline{z} 
+ B\cdot d\overline{z}$ for some section $B$  of $L_{\cV}$.  

A $C^{\infty}$ {\em meromorphic Higgs field}\label{higgs-defn}
 $\Phi$ on $\boldp^1$ 
is a $C^{\infty}$ section of 
$(T^*)^{1,0}(\boldp^1)\otimes\theo(D)\otimes \widehat{L}_{\cV}^*$,
for some   effective divisor $D$ on $\boldp^1$.  
If one fixes a section $s$ of $(T^*)^{0,1}(\pline)\otimes\theo(D)$, 
then it makes sense to talk of a meromorphic Higgs field lying in
 $s\otimes (\widehat{L}_{\cV}^*)_\lambda$.
\end{defn}

In our case, we will fix the divisor $0+\infty$ on $\pline$, as well as
coadjoint orbits in the fibres of $\widehat{L}_{\cV}$  over $0$ and $\infty$.
By
\cite[Proposition~3.2]{EF}, these orbits 
 correspond to isomorphism  classes of holomorphic structures 
 in  $\cV|_{E\times\{0\}}$ and $\cV|_{E\times\{\infty\}}$, and we choose
these to be the split bundle $(\cL^*)^{d_{-1}}\oplus H\oplus \baL\,^{d_{1}}$
over $E_0$ and the trivial bundle of rank $d_{-1}+r+d_1$ over $E_{\infty}$.
As before, we make these bundles 
 $\cs$-equivariant by using weights $0, -1, -2$ for the factors
of $\cV$ along $E_0$ and weight $0$ along $E_\infty$.
The link between $\sigma$-monads on $E$ and Higgs bundles on
 $\boldp^1$ is then the following correspondence.

\begin{prop}\label{final-prop}
Let $S=\skl(E,\cL,\sigma)$ and keep the above notation.
Then there exists a bijective correspondence between isomorphism classes of 
\begin{enumerate}
\item[(1)] $\sigma$-monads $\bK$ of the form \eqref{sigma Kronecker complex}
that satisfy
$\HH^0(\bK)|_E\cong H$, and
\item[(2)] $\cs$-equivariant meromorphic Higgs pairs 
$(\overline{\partial}_B, \Phi)$ on the $C^{\infty}$-bundle $\cV$ 
for which $\Phi$ lies in $z^{-1}dz\otimes (\widehat{L}_{\cV}^*)_1$ 
with residues at $0$ and $\infty$ in the given coadjoint orbits.
\end{enumerate}
\end{prop}

\begin{proof}
Since this result is tangential to the results of the paper we will
 only outline the proof.

Using a variant of the Rees construction,
one can identify  a $\sigma$-monad $\bK$ 
with a holomorphic structure $D_\cV$  on the equivariant vector
bundle $\cV$  over $E\times\pline$. Now
use the technique of \cite{GM}:  given a
$\overline{\partial}$-operator $D_{\cV}$ on $\cV$ over $E\times\pline$, split
it into components $D_E$ and $D_{\pline}$ taking values in  $\cV\otimes
(T^*)^{0,1}E$, respectively  $\cV\otimes (T^*)^{0,1}\pline$.   The
correspondence of the proposition  then takes $D_{\cV}$ to the pair for which 
$\overline{\partial}_B = D_{\pline}$ and   $\Phi = D_E\cdot z^{-1}dz$.  One
then only needs to check that the gauge groups act in the same way  on the two
sides and that $D_{\cV}\circ D_{\cV} = 0$ is equivalent to the two equations 
$\overline{\partial}_B\circ\overline{\partial}_B=0$ and 
$\overline{\partial}_B(\Phi)= 0$. \end{proof}

\section*{Index of Notation}\label{Index}

Standard definitions concerning $\Qgr S$ can be found 
on pages \pageref{cg-defn}---\pageref{chi-defn}
and these will not be listed in this index.
\begin{multicols}{2}

$\ASthree=\ASthree(k)$,  $\ASthree'$,\hfill\pageref{ASthree-defn}

$A(S)$, \hfill\pageref{A(S)-defn}

$c_1(-)$, the first chern class, \hfill\pageref{chern-defn}, 
\pageref{kronecker-chern}

Canonical filtration, \hfill\pageref{canonical-defn}

\v{C}ech complex ${\mathbf C}^\bullet $, \hfill\pageref{cech-defn}

Diagonal bigraded algebra $\Delta$, \hfill\pageref{delta-defn}

$\epsilon(S)$,\hfill\pageref{epsilon-defn}

Equivalent families, \hfill\pageref{equiv-def}

Equivalent filtration, \hfill\pageref{canonical-defn}
 
Euler characteristic $\chi$,  \hfill\pageref{euler-defn}, \pageref{euler-kronecker}

Family of algebras in $\ASthree$, \hfill\pageref{family-as-defn}

Family of complexes, etc,  \hfill\pageref{families def},
\pageref{typical Kronecker-defn}

Flat module in $\Qgr S_R$, \hfill\pageref{flat-defn}

Framed complex,  \hfill\pageref{framed-complex}

$\Gamma^*(\cM)$, $\Gamma(\cM)$, \hfill\pageref{gamma-fn}

Geometrically (semi)stable, 
\hfill\pageref{stability definition for modules}, \pageref{geom-Kr}

Grassmannians $\Gr_{d_{-1}},\Gr^{d_{1}},\GR_\cB$,\hfill\pageref{Grass-defn}

$\GR=\GR_S$, \hfill\pageref{en-defn}

$\uHom_{S\on{-qgr}}$, \hfill\pageref{m-vee}

$\cK^{ss}$, $\widehat{\cK}^{ss}$, \hfill \pageref{Kronecker functors}

KL-pair,    \hfill \pageref{KL-pair}

Kronecker complex, \hfill\pageref{typical Kronecker complex}

 Hilbert polynomial,   \hfill\pageref{hilbert-defn}, 
\pageref{hilbert-kronecker}

Line bundle in $\Qgr S$, \hfill\pageref{line-ble-defn}

$\cM|_E$,  \hfill\pageref{defn-restriction}

$\cM|_A$, \hfill\pageref{M-to-A-defn}

Maximal subcomplex, \hfill\pageref{max-subcomplex}

Meromorphic Higgs field,  \hfill \pageref{higgs-defn}

Moduli  space (fine or coarse), \hfill\pageref{moduli-defn}

$\M_S^{ss}(r,c_1,\chi)$,  $\M^s_S(r,c_1,\chi)$, \hfill\pageref{moduli-1-defn}

   $\M_{\cB}^s(r,c_1,\chi) $,  \hfill\pageref{M_S-defn}

$\cM'', \cM',\cM=\cM(\cF_1,\dots,\cF_n)$,   \hfill\pageref{M-dash-defn}

$\cM_0'', \cM_0',\cM_0$,   \hfill\pageref{M-dash-oh}

Monad, \hfill\pageref{monad-defn}

$\cN^\vee$,\hfill \pageref{M-vee-defn}

$N_{\cB}$,\hfill\pageref{straight-N-def}

$\cN=\cN_S$,  \hfill\pageref{en-defn}

$\N^{ss}_{\cB}$,  \hfill\pageref{smooth-subdefn}

$\N^s$, $\N^{ss}$, \hfill\pageref{ns-defn}

$\N'$, \hfill\pageref{description of Nss}

Normalized module, \hfill\pageref{normalized-defn}

Normalized complex, \hfill\pageref{norm-complex}

$\Omega^1$, $\widetilde{\Omega}^1$, \hfill\pageref{omega-defn}

$\cP_A$, \hfill\pageref{P-defn}

$(\mathbf P_S\smallsetminus E)^{[n]}$, \hfill\pageref{hilb-defn}

$\Phi'', \Phi', \Phi$,  \hfill\pageref{Phi-defn}
 
Rank of a module, \hfill\pageref{rank-defn}

$\cR(P)$, $\cS(P)$, $ \cV(P)$,\hfill \pageref{cS-defn}

$(\qgr S_R)_{\on{VC} }$,\hfill \pageref{thesheaf}

$S_B$, \hfill\pageref{S-sub-B}

Schematic algebra, \hfill\pageref{schematic-defn2}

(Semi)stable  modules, 
\hfill\pageref{stability definition for modules}

(Semi)stable    complexes, 
\hfill\pageref{semistability definition for complexes}

$\sigma$-Kronecker complex, $\sigma$-monad,  \hfill\pageref{sigma-monad-defn}

Standard complex, \hfill\pageref{standard-defn}

Strongly noetherian ring, \hfill\pageref{strong-defn}

Torsion, torsion-free modules, \hfill\pageref{tf-defn}

Vector bundle in $\Qgr S$, \hfill\pageref{vectorbundledefn}

$\cV_S$, \hfill\pageref{VS-defn}

\end{multicols}


\begin{thebibliography}{EGA33}


\bibitem[Aj]{Aj}  K. Ajitabh, Residue complex for Sklyanin 
algebras of dimension three, {\em Adv. in Math.} {\bf 144} (1999), 137-160. 



 \bibitem[ASh]{ASh}  M. Artin and W. F. Schelter, Graded algebras of global
 dimension 3,
 {\em Adv. in Math.} {\bf 66} (1987), 171-216.

  \bibitem[ASZ]{ASZ} M. Artin, L. W. Small and J. J. Zhang,
  Generic flatness for strongly noetherian algebras,
  {\em J. Algebra} {\bf 221}  (1999), 579-610. 
 
\bibitem[AS]{AS} M. Artin and J. T. Stafford, Noncommutative graded domains with
quadratic growth, {\em Invent. Math.,} {\bf 66} (1995), 231-276.

  \bibitem[ATV1]{ATV1}
M.~Artin, J.~Tate and M.~Van~den Bergh, Some algebras associated
to  automorphisms of elliptic curves,  in {\em The Grothendieck Festschrift},
vol.~1,  Birkh\"auser, Boston, 1990, pp.~33--85.

\bibitem[ATV2]{ATV2}
\bysame, Modules over regular algebras of dimension 3, {\em Invent. Math.}
{\bf  106} (1991), 335--388.

\bibitem[AV]{AVdB}
M.~Artin and M.~Van~den Bergh,  Twisted homogeneous coordinate
rings, {\em J.  Algebra} {\bf 133} (1990), 249--271. 
  
\bibitem[AZ1]{AZ}
M.~Artin and J.~J. Zhang, Noncommutative projective schemes,
{\em Adv. in Math.} {\bf 109}  (1994), no.~2, 228--287.
 
 
 \bibitem[AZ2]{AZ2} \bysame, Abstract Hilbert schemes, 
 {\em Algebr. Represent. Theory} {\bf 4} (2001), 305--394.

 

  \bibitem[BGK1]{BGK}
 V.~Baranovsky, V.~Ginzburg and A.~Kuznetsov,
Quiver varieties and a non-commutative $P^2$,  
{\em Compositio Math.} {\bf 134} (2002), 283--318.

\bibitem[BGK2]{BGK2} \bysame, Wilson's Grassmannian and a
noncommutative quadric, {\em Int. Math. Res. Not.} {\bf 2003:21} (2003), 
1155-1197.

 
 \bibitem[Be]{Beauville} A. Beauville, Vari\'et\'es K\"ahleriennes dont la 
 premi\`ere classe de Chern est nulle.
 {\em J. Differential Geom.} {\bf 18} (1983), 755--782.
 

\bibitem[BJL]{BJL}
D.~ Berenstein, V.~Jejjala and R.~G.~ Leigh, Marginal and
relevant deformations of $N=4$ field theories and non-commutative moduli
spaces of vacua, {\em Nuclear Phys. B} {\bf 589} (2000), 196--248.

\bibitem[BL]{BL}
D.~ Berenstein and R.~G.~ Leigh, Resolution of stringy singularities by
noncommutative algebras,  {\em J. High Energy Phys.,} 2001, no. 6, Paper
30.


  \bibitem[BW1]{BW1}
 Y.~Berest and G.~Wilson, Automorphisms and ideals of the Weyl algebra,
  {\em Math. Ann.} {\bf 318} (2000), 127--147.
 
   \bibitem[BW2]{BW2} \bysame,
    Ideal classes of the Weyl algebra and noncommutative projective 
    geometry.  With an appendix by Michel Van den Bergh,
   {\em Int. Math. Res. Not.}  (2002) {\bf no. 26}, 1347--1396.


 \bibitem[CH]{CH}
 R.~C.~Cannings and M.~P.~Holland,
 Right ideals of rings of differential operators,  {\em J. Algebra} {\bf 167}
  (1994),  116--141
 
 \bibitem[CDS]{CDS}
A.~Connes, M.~Douglas and A.~Schwarz,
Noncommutative geometry and matrix theory: compactification on tori.
{\em J. High Energy Phys.} {\bf 1998}, no. 2, Paper 3, 35 pp.


\bibitem[DV]{dNvdB}
K. de Naeghel and M. van den Bergh, Ideal classes of three dimensional 
Sklyanin algebras,  Preprint, 2002.

 \bibitem[Di]{Di} J.~Dixmier, Sur les al\`ebres de Weyl,
 {\em Bull.~Soc.~Math.~France} {\bf 96} (1968), 209-242.

\bibitem[DN]{DN} M.~R.~Douglas and N.~A.~Nekrasov, 
Noncommutative field theory,
{\em Rev. Modern Phys.},  {\bf 73} (2001), 977--1029.

\bibitem[DL]{Drezet-Le Potier}
J.~M.~Drezet and J.~Le~Potier, Fibr\'es stables et fibr\'es exceptionnels 
sur ${\mathbf P}^2$,
{\em Ann. Scient. \'Ec. Norm. Sup.}, {\bf 18} (1985), 193--244.

\bibitem[Ei]{Ei} D. Eisenbud, \emph{Commutative Algebra with a View Toward
Algebraic Geometry}, Graduate Texts in Math. 150, Springer-Verlag, New York, 
1995.

\bibitem[EF]{EF} P. I. Etingof and I. B.  Frenkel,
Central extensions of current groups in two dimensions, 
{\em Comm. Math. Phys.}, {\bf 165} (1994),  429--444.

\bibitem[FO]{FO2}
B.~L. Feigin and A.~V. Odesskii,
  Vector bundles on an elliptic curve and Sklyanin algebras,
In {\em Topics in quantum groups and finite-type invariants},
{\em Amer. Math. Soc. Transl.} {\bf 185} (1998), 65--84.
   
\bibitem[GM]{GM}
H. Garland and M. K. Murray, Kac-Moody monopoles and periodic 
instantons,  {\em Comm. Math. Phys.}
{\bf 120} (1988), no. 2, 335--351. 

 
\bibitem[Gi]{Ginzburg}
V. Ginzburg, Non-commutative symplectic geometry, quiver varieties,
and operads, {\em Math. Res. Lett.} {\bf 8} (2001), no. 3, 377--400.


\bibitem[GN]{GN}
A. Gorsky and N. Nekrasov, Elliptic Calogero-Moser system from
two-dimensional current algebra.  arXiv:hep-th/9401021.

\bibitem[EGA]{EGA}
A. Grothendieck, \'{E}l\'{e}ments de g\'eom\'etrie
alg\'ebrique, Chapters III and IV, 
{\em Inst. Hautes \'Etudes Sci. Publ. Math.}, {\bf 11} (1961), {\bf 17} (1963),
 {\bf 20} (1964), {\bf 24} (1965),
{\bf 28} (1966), {\bf 32} (1967).


\bibitem[Ha]{H}
R. Hartshorne, {\em Algebraic geometry}, Springer Verlag, Berlin, 1977. 

\bibitem[HL]{HLbook} D.~Huybrechts and M.~Lehn, 
 {\em The geometry of moduli spaces of sheaves},
 Friedr. Vieweg \& Son, Braunschweig, 1997.

\bibitem[KKO]{KKO}
A.~Kapustin, A.~Kuznetsov and D.~Orlov,
Noncommutative instantons and twistor transform,
 {\em Comm.~Math.~Phys.} {\bf 221}
(2001), 385--432.

\bibitem[LB1]{Le} L. LeBruyn,
Moduli spaces for right ideals of the Weyl algebra,
{\em J. Algebra}, {\bf  172}  (1995), 32--48.

 
\bibitem[LB2]{LeBruynbook}
\bysame, {\em Noncommutative geometry at $n$}, in preparation.
 

\bibitem[LP1]{Le Potier} J. Le Potier, \`A propos de la construction
 de l'espace
de modules des faisceaux semi-stables sur le plan projectif, {\em Bull. Soc.
Math. France} {\bf 122} (1994), 363--369.

\bibitem[LP2]{LP2} \bysame, {\em Lectures on vector bundles},  
CUP, Cambridge,
1997.

 
\bibitem[Ma]{Ma}
 Yu.~I.~Manin, {\em Quantum groups and noncommutative geometry},
  Universit\'e de Montr\'eal, Centre de Recherches Math\'ematiques,
  Montr\'eal, 1988.

\bibitem[MR]{MR}
J.~C. McConnell and J.~C. Robson, {\em Noncommutative {N}oetherian
rings}, John  Wiley {\&} Sons, New York, 1987.


\bibitem[Na]{Nakajima} H. Nakajima, Instantons on ALE spaces, quiver varieties,
 and Kac-Moody algebras,
{\em Duke Math. J.} {\bf 76} (1994),  365--416.

\bibitem[NV]{NV}
C.~N\u ast\u asescu and F.~Van~Oystaeyen, {\em Graded ring theory},
North Holland, Amsterdam,  1982.

\bibitem[Ne]{Nekrasov}
N. Nekrasov, Infinite-dimensional algebras, many-body systems and
gauge theories, in {\em Moscow seminars in mathematical physics}, 
{\em Amer. Math. Soc. Transl. Ser. 2} {\bf 191} (1999), 263--299.


\bibitem[NSc]{Nekrasov-Schwarz}
N. Nekrasov and A. Schwarz, Instantons on 
noncommutative ${\mathbf R}^4$ and 
$(2,0)$ superconformal six-dimensional theory,
 {\em Comm. Math. Phys.} {\bf 198} (1998),
 689--703.

\bibitem[Od]{Od}  A.~V.~Odesskii, Elliptic algebras, {\em 
Russian Math. Surveys}, {\bf 75} (2002), 1127-1162.
  
  \bibitem[OSS]{OSS} C. Okonek, M. Schneider and H. Spindler, 
  {\em Vector bundles on complex projective spaces},
  Progress in Math. 3, Birkh\"auser,
  Boston, 1980.
  
\bibitem[Pl]{Polishchuk} A. Polishchuk, Poisson structures and
birational morphisms associated with bundles on elliptic curves,
{\em Internat. Math. Res. Notices} {\bf 1998}, no. 13, 683--703.

 \bibitem[Pp]{Po} N. Popescu, {\em Abelian categories with applications to rings
  and modules}, Academic Press, London, 1973.

\bibitem[RSS]{RSS} R. Resco, L.W. Small
and J. T. Stafford, Krull and global dimension of semiprime noetherian PI
rings, {\em Trans. Amer. Math. Soc.} {\bf 274} (1982), 285-296.

\bibitem[Rg]{rog}
D. Rogalski, Examples of generic noncommutative surfaces,
{\em Adv. in Math.}, to appear. 
  arXiv:math.RA/0203180.
  
\bibitem[Ru]{Rudakov}
A. Rudakov, Stability for an abelian category, {\em J. Algebra} {\bf 197} 
(1997), 231--245.

\bibitem[SW]{Seiberg-Witten}
N. Seiberg and E. Witten, String theory and noncommutative
geometry, {\em J. High Energy Phys.} {\bf 1999}, no. 9, Paper 32, 93 pp.


\bibitem[Se]{Seshadri}
 C.~S.~Seshadri, Geometric reductivity over arbitrary base,
 {\em Adv. in Math.} {\bf 26} (1977), 225--274.

\bibitem[Si]{Simpson}
C. Simpson, Moduli of representations of the fundamental group of
a smooth projective variety I,  {\em Publ. Math. IHES},
{\bf 79} (1994), 47--129.

  \bibitem[Sm]{Sm} S.~P.~Smith, Some finite-dimensional algebras related to 
  elliptic curves,  in {\em Representation theory of algebras 
  and related topics
(Mexico City, 1994)}, 315--348, CMS Conf. Proc., 
19, Amer. Math. Soc., Providence, RI, 1996.


\bibitem[St1]{St2} J.~T.~Stafford, 
Stably free  projective right ideals, {\em Compositio Math.}, {\bf 54}
(1985), 63-78.

\bibitem[St2]{Stafford} \bysame,   Noncommutative projective geometry,
{\em Proceedings of the ICM, Vol. II (Beijing, 2002)}, 93--103.




\bibitem[SV]{SVdB} J.~T.~ Stafford and M.~Van den Bergh, 
Noncommutative curves and noncommutative surfaces,
{\em Bull.~Amer.~Math.~Soc.} {\bf 38} (2001), 171-216.

\bibitem[TV]{TVdB}
J.~Tate and M.~Van~den Bergh, Homological properties of {S}klyanin
  algebras, {\em Invent. Math.} {\bf 124} (1996), 619--647.

\bibitem[VW1]{VW2} F. Van Oystaeyen and L. Willaert, 
Cohomology of schematic algebras, {\em J. Algebra} {\bf 185} (1996), 74-84. 

\bibitem[VW2]{VW3} \bysame, 
Examples and quantum sections of schematic algebras,  {\em J. Pure Appl.
Algebra} {\bf 120} (1997), 195--211.

 \bibitem[Ve]{Ve}
A.~B. Verevkin, On a non-commutative analogue of the category of coherent
  sheaves on a projective scheme, {\em Amer. Math. Soc. Transl.} {\bf 151}
   (1992), 41--53.

\bibitem[Wi]{Wilson}  G. Wilson, Collisions of Calogero-Moser particles
and an adelic Grassmannian,  {\em Invent. Math.} {\bf 133} (1998), 1--41.

\bibitem[YZ]{YZ}
A.~Yekutieli and J.~J. Zhang, Serre duality for noncommutative
projective  schemes, {\em Proc. Amer. Math. Soc.} {\bf 125} (1997),  697--707.

\bibitem[Zh]{Zh2} J. J. Zhang, Connected 
graded Gorenstein rings with enough normal
elements, {\em J. Algebra}, {\bf 189}
   (1997),  390--405

\end{thebibliography}
\end{document}